\documentclass[twoside,pdftex]{amsart}

\usepackage[latin1]{inputenc}
\usepackage{color}
\usepackage{url}
\usepackage{graphicx}
\usepackage{pstricks}
\usepackage[bookmarks=false,pdfborder={0 0 0.05}]{hyperref}
\usepackage{amssymb}

\definecolor{sagecolor}{rgb}{.78, .36, .04}
\newcommand{\sageex}[1]{\vspace{0.5em} {\small \tt #1} \vspace{0.5em}}
\def\sageret#1{\begin{center} #1 \end{center}}

\def\sagepromt{{\normalsize \color{sagecolor} sage: }}
\def\sagedots{{\normalsize \color{sagecolor}....: }}
\def\sage{{\tt sage}}
\def\sagecomb{{\tt sage-combinat}}

\def\True{{\tt True}}
\def\False{{\tt False}}
\def\None{{\tt None}}

\newtheorem{theorem}{Theorem}[section]
\newtheorem{algorithm}[theorem]{Algorithm}
\newtheorem{proposition}[theorem]{Proposition}
\newtheorem{corollary}[theorem]{Corollary}

\theoremstyle{definition}
\newtheorem{definition}[theorem]{Definition}

\newtheorem{conjecture}[theorem]{Conjecture}
\newtheorem{remark}[theorem]{Remark}  

\begin{document}
  \title[A compendium on the cluster algebra and quiver package in {\large \sage}]{A compendium on the cluster algebra\\ and quiver package in {\LARGE \sage}}
  \author{Gregg Musiker}
  \address{School of Mathematics, University of Minnesota, Minneapolis, MN 55455, USA}
  \email{musiker@math.umn.edu}
  \urladdr{http://www.math.umn.edu/$\sim$musiker}

  \author{Christian Stump}
  \address{LaCIM, Universit\'e du Qu\'ebec \`a Montr\'eal, Montr\'eal (Qu\'ebec), Canada}
  \email{christian.stump@lacim.ca}
  \urladdr{http://homepage.univie.ac.at/christian.stump/}

  \date{\today}

  \maketitle

\setcounter{section}{-1}

  \tableofcontents

\section{Preface}

The idea for a cluster algebra and quiver package in the open-source computer algebra system \sage\ was born during the \sage\ days 20.5 which were held at the Fields Institute in May 2010. The purpose of this package is to provide a platform to work with cluster algebras in graduate courses and to further develop the theory by working on examples, by gathering data, and by exhibiting and testing conjectures. In this compendium, we include the relevant theory to introduce the reader to cluster algebras assuming no prior background; this exposition has been written to be accessible to an interested undergraduate.\\

Version 1.0 of the software package and this compendium is the result of many discussions on mathematical background and on implementation algorithms, and of many, many hours of coding. It is part of the \sagecomb\ project \cite{sagecomb}.

For more information on \sage, in particular on a detailed description how to install the program, we refer the reader to \url{http://www.sagemath.org} \cite{sage}; for more on the \sagecomb\ project, see \url{http://wiki.sagemath.org/combinat}.  Throughout this compendium, we include examples that the user can run in the {\tt sage-notebook} or \sage\ command line. The package provides as well an interactive mode for the {\tt sage-notebook} as shown in Figure~\ref{fig:interactivemode} at the end of Section \ref{sec:quivers}. We will close with a detailed description of the data structures and methods in this package. We follow the usual \sage\ convention of indexing all lists starting at zero.\\

{\bf Currently, installing the \sagecomb\ queue is a necessary requirement for working with the cluster algebra and quiver package.} In order to install the \sagecomb\ queue, you have to, after installing \sage, run the command
\sageex{

  \# ./sage -combinat install

}
\noindent on the {\tt unix} command line. Once the \sagecomb\ branch is created, one can use the command

\sageex{
  \# ./sage -combinat update
}

\noindent to update to the latest version of the \sagecomb\ queue, or one can use the command

\sageex{
  \# ./sage -combinat upgrade
}

\noindent to update to the latest version of the \sagecomb\ queue and to upgrade \sage\ to its newest version.  For more detailed explanations, please visit the \sagecomb\ wiki page.  Installing the \sagecomb\ queue will eventually become obsolete after the project has gone through testing and reviewing processes, which might take time due to the involvedness of the algorithms (especially on mutation type detections). \\

This current version 1.0 should not be considered a complete, unchangeable, totally stable version. We will keep working on this project by fixing bugs, improving algorithms, and by adding functionalities. So it might be a good idea to update the \sagecomb\ queue once in a while, especially if you have encountered a problem.

We anticipate this ongoing project being improved with feedback from users.  We are very interested in getting any type of feedback: on ways in which the package has been used, on features people like or that could be done better, or requests for new functionalities.  If you are interested in helping us make improvements or further develop this package, we would be happy to have you involved.\\

Several other people have also worked on software for computations involving cluster algebras and quiver representations.  Links to these are available at Fomin's Cluster Algebra Portal {\tt http://www.math.lsa.umich.edu/$\sim$fomin/cluster.html}.  This software includes work of Chapoton \cite{ChapMutations}, Dupont-P\'erotin \cite{DP-QME}, Keller \cite{Keller}, and L. Williams.

\vspace{0.5em}

\noindent {\bf Acknowledgements.}~~ 
We thank Franco Saliola and S\'ebastien Labb\'e for help on details of the \sage\ development process. We also thank Florian Block, Hugh Thomas, and Leandro Vendramin for several contributions in the early stages of this project.  We thank William Stein, Florent Hivert, Nicolas Thi\'ery, and all of the developers of \sage\ and \sagecomb\ for their continued work on this open-source mathematical software. Finally, we like to thank the Fields Institute for its hospitality during the \sage\ days 20.5 in May 2010 where this project was initiated.  We also thank Bernhard Keller for a careful reading and numerous helpful edits to this guide.

\section{Introduction} \label{intro}

Cluster algebras, invented by Fomin and Zelevinsky \cite{ClustI}, are certain commutative algebras which are isomorphic to subalgebras of the field of rational functions.  Each cluster algebra has a distinguished set of generators called cluster variables; this set is a union of overlapping algebraically independent finite subsets called clusters, which together have the structure of a simplicial complex.  The clusters are related to each other by binomial exchange relations. In the past ten years, such algebras have been found to be related to a number of other topics such as quiver representations, tropical geometry, canonical bases of semisimple algebraic groups, total positivity, generalized associahedra, Poisson geometry, and Teichm\"{u}ller theory.

Usually, when one defines an \emph{algebra} $\mathcal{A}$, one describes it by writing down the \emph{generators} and \emph{relations} of $\mathcal{A}$.  Instead, when working with a \emph{cluster algebra}, only a finite set of generators are provided at first, along with combinatorial data that allows one to algebraically construct the rest of the generators by applying a sequence of exchange rules.  With this definition in mind, a \emph{seed} for a cluster algebra $\mathcal{A}$ is a pair $(\mathbf{x}, B)$, where $\mathbf{x}$ denotes the \emph{initial cluster},  and $B$ denotes an \emph{exchange matrix} (or $B$-matrix)\footnote{Technically, this is the definition for a seed of a cluster algebra of geometric type.  We give a more general definition of cluster algebra seeds in the next section.}.  Here, the cluster $\mathbf{x}$ consists of exchangeable generators, known as \emph{cluster variables} and non-exchangeable generators, known as \emph{coefficients} or \emph{frozen variables}.

One of the simplest families of cluster algebras are those which are coefficient-free and rank two.  Such algebras are parametrized by two positive integers $(b,c)$, and the associated cluster algebra $\mathcal{A}(b,c)$ is defined to be the algebra generated by the set $\{x_n\}_{n \in \mathbb{Z}}$, where for $n \not \in \{0,1\}$, 
$$x_n = 
\frac{x_{n-1}^b+1}{x_{n-2}}  \mbox{ if $n$ is even, and }
\frac{x_{n-1}^c+1}{x_{n-2}}  \mbox{ if $n$ is odd}.
$$
These are implemented in \sage, for example (letting $b=2$, and $c=3$) as 

\sageex{
\sagepromt S23 = ClusterSeed(['R2',[2,3],2]); S23

\sageret{A seed for a cluster algebra of rank 2 of type ['R2',[2,3],2]}
}
\noindent Here, {\tt 'R2'} refers to \lq\lq rank $2$\rq\rq, {\tt [2,3]} gives the parameters. For an explanation of the final $2$, we refer to Section~\ref{QMT}. Notice that if instead we let $b=1$ and $c=1$, we obtain 

\sageex{
\sagepromt S11 = ClusterSeed(['R2',[1,1],2]); S11

\sageret{A seed for a cluster algebra of rank 2 of type ['A',2]}
}
\noindent We will see more examples of this phenomenon in a moment, but the point is that when $(b,c)=(1,1)$, the associated cluster algebra is of ``type $A_2$''.

Let us keep working with the cluster seed $S11$ at the moment.  We can see the $B$-matrix and initial cluster corresponding to this seed quite easily.

\sageex{
\sagepromt S11.cluster()

$$\left[x_{0}, x_{1}\right]$$

\sagepromt S11.b\_matrix()

$$\left(\begin{array}{rr}
0 & 1 \\
-1 & 0
\end{array}\right)$$
}

\noindent Using this data, it is possible to construct the other generators of $\mathcal{A}(1,1)$ by applying a sequence of exchanges.  We define mutation in general down below.  For now, let us mention that if we start with the initial cluster $\left[x_0, x_1\right]$, and mutate in the $0$th direction, we replace the $x_0$ with $x_2$, defined as $x_2 = \frac{x_1 + 1}{x_0}$.  This gives us a new seed, whose cluster is $[x_2, x_0]$.

\sageex{
\sagepromt S11.mutate(0); S11.cluster()

$$\left[\frac{x_{1} + 1}{x_{0}}, x_{1}\right]$$
}

\noindent The exchange matrix of this new seed is simply $-B$. 

\sageex{
\sagepromt S11.b\_matrix()

$$\left(\begin{array}{rr}
0 & -1 \\
1 & 0
\end{array}\right)$$
}

We can continue this procedure, and now mutate in the $1$st direction, letting $x_3 = \frac{x_2+1}{x_1}$ replace $x_1$.

\sageex{
\sagepromt S11.mutate(1); S11.cluster()

$$\left[\frac{x_{1} + 1}{x_{0}}, \frac{x_{0} + x_{1} + 1}{x_{0} x_{1}}\right]$$

\sagepromt S11.b\_matrix()
$$\left(\begin{array}{rr}
0 & 1 \\
-1 & 0
\end{array}\right)$$
}

\noindent Notice that after this mutation, the exchange matrix of the obtained seed is $B$ again.  Consequently, we can iterate this procedure, applying the {\tt mutate} command over and over.  If we want to do this more efficiently, we can as well call {\tt mutate} with a list of indices to apply from left to right.

\sageex{
\sagepromt S11.mutate([0,1,0,1])
}

\noindent If we are not only interested in the final seed, we can instead use the procedure {\tt mutation\_sequence}.  Before doing that, we reset the cluster to the initial sequence of variables (in the initial seed).  

\sageex{
\sagepromt S11.reset\_cluster();

\sagepromt S11.mutation\_sequence([0,1,0,1,0],return\_output='matrix')

$$\left[\left(\begin{array}{rr}
0 & 1 \\
-1 & 0
\end{array}\right), \left(\begin{array}{rr}
0 & -1 \\
1 & 0
\end{array}\right), \left(\begin{array}{rr}
0 & 1 \\
-1 & 0
\end{array}\right), \left(\begin{array}{rr}
0 & -1 \\
1 & 0
\end{array}\right), \left(\begin{array}{rr}
0 & 1 \\
-1 & 0
\end{array}\right), \left(\begin{array}{rr}
0 & -1 \\
1 & 0
\end{array}\right)\right]$$

\sagepromt S11.reset\_cluster();

\sagepromt S11.mutation\_sequence([0,1,0,1,0],return\_output='var')

$$\left[\frac{x_{1} + 1}{x_{0}}, \frac{x_{0} + x_{1} + 1}{x_{0} x_{1}}, \frac{x_{0} + 1}{x_{1}}, x_{0}, x_{1}\right]$$
}

Here, the first command returns the sequence of exchange matrices obtained from this sequence of mutations, including the initial one.  Notice, the sequence is exactly $[B, -B, B, -B, B, -B]$.  The second command returns the list of cluster variables encountered as these exchanges occur.  In the rank two case, this list is equivalent to $[x_2,x_3,x_4,x_5,x_6]$ corresponding to the $(b,c)=(1,1)$-sequence $\{x_n\}$ referred to above.  

Notice, that we have already found an interesting pattern, that is after five exchanges, we have arrived back essentially\footnote{To be precise, this seed uses matrix $-B$ (equivalently $B^T$) instead of $B$, but these seeds are the same \lq\lq up to equivalence\rq\rq, see Remark~\ref{rem:equivalence}.} at the same seed with which we started.  This is therefore known as a cluster algebra of \emph{finite type} and \emph{finite mutation type}.  Both of these concepts will be described in more detail below.

For our next example, we look at the $(b,c)=(2,2)$ case, again a rank two cluster algebra.

\sageex{
\sagepromt S22 = ClusterSeed(['R2',[2,2],2]); S22

\sageret{A seed for a cluster algebra of rank 2 of type ['A',[1,1],1]}
}

\noindent Here again, notice that this specific rank two cluster algebra is recognized.  In this case, this is our notation for a cluster algebra of affine type $\widetilde A_{1,1}$.  We again, run the procedure {\tt mutation\_sequence}, and obtain the following:

\sageex{
\sagepromt S22.mutation\_sequence([0,1,0,1,0],return\_output='var')\\

$\begin{array}{c}

\hspace{-15em}\left[\frac{x_{1}^{2} + 1}{x_{0}}, \frac{x_{1}^{4} + x_{0}^{2} + 2 x_{1}^{2} + 1}{x_{0}^{2} x_{1}}, \frac{x_{1}^{6} + x_{0}^{4} + 2 x_{0}^{2} x_{1}^{2} + 3 x_{1}^{4} + 2 x_{0}^{2} + 3 x_{1}^{2} + 1}{x_{0}^{3} x_{1}^{2}}, \right.\\[15pt]
\hspace{-14em}\frac{x_{1}^{8} + x_{0}^{6} + 2 x_{0}^{4} x_{1}^{2} + 3 x_{0}^{2} x_{1}^{4} + 4 x_{1}^{6} + 3 x_{0}^{4} + 6 x_{0}^{2} x_{1}^{2} + 6 x_{1}^{4} + 3 x_{0}^{2} + 4 x_{1}^{2} + 1}{x_{0}^{4} x_{1}^{3}}, \\[15pt]
\hspace{-1em}\left.\frac{x_{1}^{10} + x_{0}^{8} + 2 x_{0}^{6} x_{1}^{2} + 3 x_{0}^{4} x_{1}^{4} + 4 x_{0}^{2} x_{1}^{6} + 5 x_{1}^{8} + 4 x_{0}^{6} + 9 x_{0}^{4} x_{1}^{2} + 12 x_{0}^{2} x_{1}^{4} + 10 x_{1}^{6} + 6 x_{0}^{4} + 12 x_{0}^{2} x_{1}^{2} + 10 x_{1}^{4} + 4 x_{0}^{2} + 5 x_{1}^{2} + 1}{x_{0}^{5} x_{1}^{4}}\right]
\end{array}$
}

\noindent Unlike the previous case, the cluster variables appear to be getting more and more complicated, and that pattern continues.  To understand these expressions better, we plug in the value $1$ for $x_0$ and $x_1$.

\sageex{
\sagepromt [cv.subs({S.x(0):1,S.x(1):1}) for cv in ms]

\sageret{[2, 5, 13, 34, 89]}
}
From this data, one might conjecture, and it is in fact true, that the sequence $$\{x_n : x_{n}x_{n-2} = x_{n-1}^2+1 \mathrm{~and~}x_0=x_1=1\}$$ is precisely the Fibonacci numbers with even index.

It is also clear that the cluster variables $x_n$'s obtained by an instance of the $(b,c)$-sequence are rational functions in the indeterminates $x_0$ and $x_1$.  More surprisingly, in spite of the divisions appearing, all such $x_n$'s are actually \emph{Laurent} polynomials, i.e. in the ring $\mathbb{Z}[x_0^{\pm 1},x_1^{\pm 1}]$.  This is actually a special case of one of the first major results in the theory of cluster algebras.

\begin{theorem} [Laurent Phenomenon \cite{ClustI,Laurent}] \label{th:Laurent}
Given any cluster algebra $\mathcal{A}$, which is parameterized by a choice of exchange pattern, a choice of coefficients 
(whose group ring is given as $\mathbb{Z}\mathbb{P}$) and a choice of initial cluster $\{x_0,x_1,x_2,\dots, x_{n-1}\}$ of generators, then all other generators, i.e. cluster variables, are Laurent polynomials in the ring $\mathbb{Z}\mathbb{P}[x_0^{\pm 1}, x_1^{\pm}, \dots, x_{n-1}^{\pm 1}]$. 
\end{theorem}
\noindent In the same paper in which Fomin and Zelevinsky prove this Laurent phenomenon, they made the following \emph{positivity conjecture}.

\begin{conjecture} [Positivity Conjecture] Given a cluster algebra $\mathcal{A}$ with an arbitrary exchange pattern, choice of coefficients $\mathbb{P}$, and an arbitrary initial cluster $\{x_0,x_1,\dots, x_{n-1}\}$, then every generator of $\mathcal{A}$ can be written in
$$\mathbb{Z}_{\geq 0}\mathbb{P}[x_0^{\pm 1}, x_1^{\pm 1}, \dots, x_{n-1}^{\pm 1}].$$
In other words, the Laurent expansions for cluster variables can be written using \emph{positive} coefficients.
\end{conjecture}

Positivity of the coefficients is significant, as it is conjecturally related to total-positivity properties of dual canonical bases \cite{Bruhat, Total, Canonical}.
Nonetheless, this conjecture is still open despite nearly a decade of work by many researchers proving it for certain families of cluster algebras.  Much of this work \cite{CalderoZel, CalderoK, CR, CP, cube, DiF1, Dup2, FQ, MusPropp, MusSch, MSW, Nak, markoff, S2, ST, Sherman, SpeyerOct, Z} has been accomplished by exploration of examples, either by hand or by computer.  As patterns to the Laurent polynomial expansions of cluster variables were noticed, the positivity conjecture and explicit formulas have been proven for more and more cases.  This software provides further tools for such explorations.

\section{What is a cluster algebra?}

In this section, we give a more general and complete definition of cluster algebras, and in the next one, we describe the connection between cluster algebras and quivers.  We say that a cluster algebra $\mathcal{A}$ is of \emph{rank} $n$ if $\mathcal{A}$ is subalgebra of an \emph{ambient field} $\mathbb{F}$ isomorphic to a field of rational functions in $n$ variables.  Algebras are typically defined by \emph{generators} and \emph{relations}, but in the case of cluster algebras, instead of being handed all the generators at once, you are instead handed a distinguished set of $n$ of them along with a constructive algorithm that can be used to obtain a complete set of generators.  Note, that in general, a cluster algebra is infinitely-generated, however, any element of this distinguished generating set can be reached in finite time.

This distinguished generating set is called the set of \emph{cluster variables}, the first $n$ of which are known as the \emph{initial} cluster variables.  Any set of algebraically independent cluster variables of maximal size has the same cardinality, namely $n$, and these $n$-subsets are known as \emph{clusters}.  Pairs of clusters $\mathbf{x}$, $\mathbf{x'}$ whose intersection is of size $(n-1)$ are related to one another by a binomial exchange relation of the form
$$
\mathbf{x'} = \left(\mathbf{x} - \{x_k\}\right) \cup x_k' \hspace{2em}\mathrm{~where~}\hspace{2em} x_k x_k' = p^+ M^+ + p^- M^-.
$$

In the second equation, $p^+$ and $p^-$ belong to a \emph{coefficient semifield} $\mathbb{P}$ (a semifield is a field missing additive inverses and is subtraction-free), and $M^+$, $M^-$ are monomials in the elements of $\mathbf{x}- \{x\}$ which share no common factor.  We let $\oplus$ denote addition in the semifield $\mathbb{P}$, and multiplication is denoted in the usual way.

\begin{definition}[Skew-symmetrizable matrices]
  An $n$-by-$n$ matrix $B$ is called \emph{skew-symmetrizable} if there exists a diagonal integer matrix $D$ with strictly positive entries on the diagonal such that $DB$ is skew-symmetric.
\end{definition}
There is an algorithmic way to determine whether a matrix is skew-symmetrizable, and to find the diagonal matrix $D$, see Section~\ref{sec:skew-sym}.

\begin{definition} [Labeled Seed for a cluster algebra]
A \emph{labeled seed} for a cluster algebra $\mathcal{A} = \mathcal{A}(\mathbf{x},\mathbf{y},B)$ is a triple $(\mathbf{x}, \mathbf{y}, B)$ where

\begin{itemize}
 \item $\mathbf{x} = \{x_0,x_1,\dots, x_{n-1}\}$ is a cluster of $n$ algebraically independent elements of ambient field $\mathcal{F}$,

 \item $\mathbf{y} = \{y_0,y_1,\dots, y_{n-1}\}$ is an $n$-tuple of coefficients, elements of the semifield $\mathbb{P}$, and 

 \item $B$ is an $n$-by-$n$ matrix that is skew-symmetrizable.

\end{itemize}

\end{definition}

A labeled seed can be \emph{mutated} into another labeled seed $(\mathbf{x'},\mathbf{y'},B')$ and all other clusters of $\mathcal{A}$, hence all other cluster variables, can be reached by applying a sequence of such mutations.

\begin{definition} [Mutation of labeled seeds]
If $\mathcal{A}$ is a cluster algebra of rank $n$ and $(\mathbf{x},\mathbf{y},B)$ is a labeled seed of $\mathcal{A}$, then for any $k \in \{0,1,\dots, n-1\}$, there exists another labeled seed $\mu_k(\mathbf{x},\mathbf{y},B) = (\mathbf{x'},\mathbf{y'},B') = \left(\mu_k(\mathbf{x}),\mu_k(\mathbf{y}),\mu_k(B)\right)$ defined as follows:

The cluster $\mathbf{x'}= \{x_0,x_1,\dots, \widehat{x_k}, \dots, x_{n-1}\} \cup \{x_k'\}$ where
$$x_k' = \left(y_k \prod_{b_{ik}>0} x_i^{b_{ik}} + \prod_{b_{ik}<0} x_i^{-b_{ik}}\right)~\bigg/~(y_k \oplus 1)x_k;$$
the coefficient tuple $\mathbf{y'} = (y_0',y_1', \dots, y_{n-1}')$ is given by 
$$y_j' = 
\begin{cases} 
y_j~ y_k^{\max(b_{kj},0)}(y_k \oplus 1)^{-b_{kj}} \mathrm{~~if~~}j\not=k, \\
1/y_k \mathrm{~~if~~}j=k
\end{cases};
$$
and the matrix $B' = \left[b_{ij}'\right]$ is given by 
$$b'_{ij} =
\begin{cases}
-b_{ij} & \text{if $i=k$ or $j=k$,} \\
b_{ij} & \text{if $b_{ik}b_{kj} \leq 0$,} \\
b_{ij} + b_{ik}b_{kj} & \text{if $b_{ik}, b_{kj} > 0$, or} \\
b_{ij} - b_{ik}b_{kj} & \text{if $b_{ik}, b_{kj} < 0$}
\end{cases}.$$
We say that $\mu_k(\mathbf{x},\mathbf{y},B)$ is the mutation in the $k$th direction.
\end{definition}

The following important observation ensures that mutation of a labeled seed is again a labeled seed.

\begin{proposition} [Proposition 4.5 of \cite{ClustI}]
If $B$ is a skew-symmetrizable matrix, then so is $\mu_k(B)$ for $0 \leq k \leq n-1$.
\end{proposition}

Another helpful fact about mutation is that it is an involution, i.e. for any $0 \leq k \leq n-1$, $\mu_k(\mu_k(\mathbf{x},\mathbf{y},B)) = (\mathbf{x},\mathbf{y},B)$.

\begin{definition} [Cluster Algebras of Geometric Type]
A cluster algebra of \emph{geometric type} is one where the coefficient semifield $\mathbb{P}$ is a special choice which leads to a great simplification of the above formulas.  In this case,
$$\mathbb{P}= Trop(u_0,u_1,u_2,\dots, u_{m-1})$$
where the addition $\oplus$ in $Trop(u_0,u_1,u_2,\dots, u_{m-1})$ is defined as $$\prod_j u_j^{a_j} \oplus \prod_j u_j^{b_j} = \prod_{j} u_j^{\min(a_j,b_j)}.$$  In particular, note that if $y_i = \prod_{j} u_j^{a_j}$ where all $a_j \geq 0$, then $y_i \oplus 1 = 1$.  With $\mathbb{P}$ given as this tropical semifield, the group ring $\mathbb{ZP}$ is simply the ring of Laurent polynomials $\mathbb{Z}[u_0^{\pm 1}, u_1^{\pm 1}, \dots, u_{m-1}^{\pm 1}]$.
\end{definition}

\begin{remark}
Letting $\mathbb{P} = Trop(u_0,u_1,\dots, u_{m-1})$, a labeled seed for a cluster algebra of geometric type is simply given as a pair $(\mathbf{x},B)$, as opposed to a triple $(\mathbf{x},\mathbf{y},B)$, where $\mathbf{x} = \{x_0,x_1,\dots, x_{n-1},u_0,u_1,\dots, u_{m-1}\} = \{x_0,x_1,\dots, x_{m+n-1}\}$ and $B$ is an $(n+m)$-by-$n$ matrix whose top $n$-by-$n$ portion is skew-symmetrizable.  This notation agrees with that of Section \ref{intro}.

By abuse of notation, we refer to the entire set of $(n+m)$ variables as a cluster, Only the first $n$ cluster variables are \emph{exchangeable}, the last $m$ of them are known as \emph{frozen variables} and appear in every single cluster.  These frozen variables encode the same data as the coefficients did in the more general case described above.  The exchange rules for mutation instead look like the following: 

$$x_k'x_k = \prod_{b_{ik}>0} x_i^{b_{ik}} + \prod_{b_{ik}<0} x_i^{-b_{ik}}$$ and the mutation rule for the $B$-matrix is unchanged except that we must mutate entires in the last $m$ rows appropriately as well.  This mutation of the last $m$ rows exactly agrees with the mutation of coefficients $\mathbf{y}$ in the general definition.  In particular, if we let $y_j = \prod_{0 \leq i\leq m} u_i^{b_{i+n,j}}$ for $0\leq j \leq n-1$, then we can recover the coefficient tuple $\mathbf{y}$ from the second halves of $\mathbf{x}$ and $B$. 
\end{remark}

\begin{remark}
Since cluster algebras of \emph{geometric type} are sufficient for many applications and all of the computations currently possible in the cluster algebra package, we henceforth discuss the theory in terms of cluster algebras only of geometric type.  We shall say that $\mathcal{A} = \mathcal{A}(\mathbf{x},B)$ is a cluster algebra of rank $n$ (with $m$ coefficients) if it is a subalgebra of an \emph{ambient field} $\mathbb{F}$ isomorphic to a field of rational functions in $(n+m)$ variables, $m$ of which are \emph{frozen}.  This is because the cluster algebra $\mathcal{A}$ is a subalgebra of $\mathbb{ZP}[x_0^{\pm 1}, \dots, x_{n-1}^{\pm 1}]$, and if $\mathbb{ZP} = \mathbb{Z}[u_0^{\pm 1}, \dots, u_{m-1}^{\pm 1}]$, then $\mathcal{A}$ can be thought of as a subalgebra of $\mathbb{Z}[x_0^{\pm 1}, \dots, x_{n-1}^{\pm 1}, u_0^{\pm 1}, \dots, u_{m-1}^{\pm 1}]$, where the $u_i^{\pm 1}$'s are simply extra generators of $\mathcal{A}$ in addition to the set of exchangeable cluster variables.

{\bf Note:} We abuse notation and often refer to the frozen variables as \emph{coefficients}; and we always denote frozen variables as $y_0$ through $y_{m-1}$ rather than the $x_{n+i}$ or $u_i$ notations used above. 
\end{remark}

We close this section with some examples and more information on some basic commands.

\sageex{
\sagepromt B3 = matrix([[0,1,0],[-1,0,-1],[0,1,0]]);

\sagepromt S3 = ClusterSeed(B3); S3

\sageret{A seed for a cluster algebra of rank 3}
}

Notice that unlike the earlier examples, the description of the seed does not include the type.  This is because the input was only the matrix, and \sage\ will not attempt to recognize the type unless it is asked for by the user or by a method.

\sageex{
\sagepromt S3.cluster()

$$[x_0, x_1, x_2]$$

\sagepromt S3.mutate(0)

\sagepromt S3.b\_matrix()

$$B' = \left(\begin{array}{rrr}
0 & -1 & 0 \\
1 & 0 & -1 \\
0 & 1 & 0
\end{array}\right)
$$

\sagepromt S3.cluster()

$$\left[\frac{x_{1} + 1}{x_{0}}, x_{1}, x_{2}\right]$$
}

We have therefore obtained a new labeled seed $(\mathbf{x'},B')$ by mutating in the $0$th direction.  Note that by default, {\tt S3.mutate(0)} acted on and changed the object {\tt S3} in place.  There is an option to leave {\tt S3} alone and just return the new object as a new output. If this behavior is desired, the command would be 

\sageex{
\sagepromt S3new = S3.mutate(0,inplace=False)
}

\noindent Since mutation is an involution, if we mutate again in the $0$th direction, we would recover the original labeled seed.  So we instead mutate in a different direction.

\sageex{
\sagepromt S3.mutate(1)

\sagepromt S3.b\_matrix()

$$B'' = 
\left(\begin{array}{rrr}
0 & 1 & -1 \\
-1 & 0 & 1 \\
1 & -1 & 0
\end{array}\right)
$$

\sagepromt S3.cluster()

$$\left[\frac{x_{1} + 1}{x_{0}}, \frac{x_{0} x_{2} + x_{1} + 1}{x_{0} x_{1}}, x_{2}\right]$$
}

Let us explain why the second element (the element $x_1'$) of this cluster is now 
$\frac{x_{0} x_{2} + x_{1} + 1}{x_{0} x_{1}}$.  This came from the exchange relation 
$$x_1 x_1' = x_2 + x_0',$$ which we read off of the second column of the exchange matrix $B' = \mu_0(B)$.  Here $x_0' = \frac{x_1+1}{x_0}$ and so we obtain the desired Laurent polynomial in terms of the initial cluster variables $x_0$, $x_1$, and $x_2$ by plugging in for $x_0'$ and simplifying.

For one more example of the exchange relation, let us now mutate in the $0$th direction again.  This corresponds to reading the first column of $B'' = \mu_1(\mu_0(B))$ which gives us the exchange relation $x_0'' = \frac{x_2 + x_1'}{x_0'}$.  Plugging in the relevant Laurent polynomials for $x_0'$ and $x_1'$, and dividing, we get a surprising cancellation and $x_0''$ is a Laurent polynomial:
$$x_0'' = \left(x_2 + \frac{x_0 x_2 + x_1 + 1}{x_0 x_1}\right)\bigg/ \left(\frac{x_1+1}{x_0} \right)= 
\frac{x_0 x_1 x_2 + x_0 x_2 + x_ 1 + 1}{x_1(x_1 + 1)} = 
\frac{x_{0} x_{2} + 1}{x_{1}}.$$

\noindent Using \sage, we see

\sageex{
\sagepromt S3.mutate(0)

\sagepromt S3.b\_matrix()

$$
\left(\begin{array}{rrr}
0 & -1 & 1 \\
1 & 0 & 0 \\
-1 & 0 & 0
\end{array}\right)
$$

\sagepromt S3.cluster()

$$\left[\frac{x_{0} x_{2} + 1}{x_{1}}, \frac{x_{0} x_{2} + x_{1} + 1}{x_{0} x_{1}}, x_{2}\right]
$$
}

\noindent One last note about mutation, we can compress the above steps as the command 

\sageex{
\sagepromt S3 = ClusterSeed(B3); S3.mutate([0,1,0]).
}

\noindent In other words, if a list is used at the input to {\tt S3.mutate}, then the seed is mutated to a new seed by applying the sequence of mutations in the same order as given by the list.

At this point, {\tt S3} is a labeled seed with matrix $B''$ and cluster $\mathbf{x''}$ as given.  However, since a labeled seed is a choice of both an exchange matrix and a cluster, we also have methods to change the cluster.  The first one is 

\sageex{
\sagepromt S3.reset\_cluster()
}

\noindent This command resets the cluster to the initial cluster $[x_0,x_1,\dots, x_{n-1}]$ while leaving the exchange matrix alone.  After running {\tt S3.reset\_cluster()}, we compute the exchange matrix and cluster, and obtain:

\sageex{
\sagepromt S3.b\_matrix()

$$
\left(\begin{array}{rrr}
0 & -1 & 1 \\
1 & 0 & 0 \\
-1 & 0 & 0
\end{array}\right)
$$

\sagepromt S3.cluster()

$$\left[x_{0}, x_{1}, x_{2}\right]
$$
}

\noindent A related command is {\tt S3.set\_cluster()} which lets the user set the initial cluster to be whatever they like.  Note that in \sage, arbitrary expressions in terms of indeterminates are not defined.  However, integers (or even rational numbers) are fair to be plugged in.  Additionally, if a rational function in terms of $x_0$ through $x_{n-1}$ is desired, this can be accomplished by the commands {\tt S3.x(0)} through {\tt S3.x(n-1)}.   

\sageex{
\sagepromt S3.set\_cluster([7,11,13]); S3.b\_matrix()

$$
\left(\begin{array}{rrr}
0 & -1 & 1 \\
1 & 0 & 0 \\
-1 & 0 & 0
\end{array}\right)
$$

\sagepromt S3.cluster()

$$\left[7,11,13\right]
$$

\sagepromt S3.mutate([0,1,2,0]); S3.cluster()

$$\left[8/11,115/77,192/1001\right]
$$
}

\noindent Note that at first glance, this might seem to falsify the Laurent Phenomenon, but it is actually allowed because all cluster variables are supposed to be Laurent polynomials in terms of the initial cluster variables.  Since the integers $7$, $11$, and $13$ are initial cluster variables, they are allowed to appear in the denominator.

\sageex{
\sagepromt S3.b\_matrix()

$$
\left(\begin{array}{rrr}
0 & 1 & 0 \\
-1 & 0 & 1 \\
0 & -1 & 0
\end{array}\right)
$$

\sagepromt S3.set\_cluster([S3.x(0)+S3.x(1),S3.x(1)\^~\hspace{-0.5em}2,S3.x(0)/S3.x(2)])

\sagepromt S3.cluster()

$$\left[x_{0} + x_{1}, x_{1}^{2}, \frac{x_{0}}{x_{2}}\right]$$

\sagepromt S3.mutate([0,1,0,2,0])

\sagepromt S3.cluster()

$$
\left[\frac{x_{1}^{2} + 1}{x_{0} + x_{1}}, \frac{x_{0} x_{1}^{2} + x_{0} x_{2} + x_{1} x_{2} + x_{0}}{x_{0} x_{1}^{2} x_{2} + x_{1}^{3} x_{2}}, \frac{x_{0} x_{1}^{2} x_{2} + x_{1}^{3} x_{2} + x_{0} x_{1}^{2} + x_{0} x_{2} + x_{1} x_{2} + x_{0}}{x_{0}^{2} x_{1}^{2} + x_{0} x_{1}^{3}}\right]
$$
}

\noindent Again, these are Laurent polynomials in terms of the \emph{initial} cluster variables obtained after setting them in this way.

\begin{definition} [Principal coefficients]
An important cluster algebra of geometric type is one \emph{with principal coefficients}.  In this case, the initial exchange matrix $B$ is $2n$-by-$n$ and where the last $n$ rows of this matrix is a rank $n$ identity matrix.  
\end{definition}

Cluster algebras with principal coefficients are fundamental, because as explained in \cite{ClustIV} by Fomin and Zelevinsky, the formula for cluster variables in a cluster algebra with general coefficients (including those not of geometric type) can be described as a simple algebraic transformation of the formulas obtained for cluster variables with principal coefficients.  See Theorem 3.7 of \cite{ClustIV} for more details.  In future versions of this package, it is anticipated that working with \emph{F-polynomials} and \emph{g-vectors} will be easier with relevant methods included. 

A cluster algebra with principal coefficients can be constructed rather simply by the command {\tt S3.principal\_extension()}.  Before demonstration, let us reset the cluster:

\sageex{
\sagepromt S3.b\_matrix()

$$\left(\begin{array}{rrr}
0 & 1 & 0 \\
-1 & 0 & 1 \\
0 & -1 & 0
\end{array}\right)
$$

\sagepromt S3.reset\_cluster(); S3.cluster()

$$\left[x_0,x_1,x_2\right]$$
}

\noindent Now, we demonstrate working with principal coefficients.

\sageex{
\sagepromt SP3 = S3.principal\_extension(); SP3

\sageret{A seed for a cluster algebra of rank 3 with 3 frozen variables}

\sagepromt SP3.b\_matrix()

$$\left(\begin{array}{rrr}
0 & 1 & 0 \\
-1 & 0 & 1 \\
0 & -1 & 0 \\
1 & 0 & 0 \\
0 & 1 & 0 \\
0 & 0 & 1
\end{array}\right)
$$

\sagepromt SP3.cluster()

$$
\left[x_{0}, x_{1}, x_{2}, y_{0}, y_{1}, y_{2}\right]
$$
}

\noindent Notice, that unlike the {\tt mutate} command, which is a verb, {\tt S3} is unaffected by the operation {\tt SP3 = S3.principal\_extension()}.

\sageex{
\sagepromt S3.b\_matrix()

$$\left(\begin{array}{rrr}
0 & 1 & 0 \\
-1 & 0 & 1 \\
0 & -1 & 0
\end{array}\right)
$$
}

\noindent Let us try an example of mutating in this cluster algebra with principal coefficients.  Here we use the command {\tt SP3.mutation\_sequence()} instead with optional arguments {\tt return\_output} which is either set to {\tt 'matrix'} or {\tt 'var'}.

\sageex{
\sagepromt SP3.mutation\_sequence([0,1,0,2],return\_output='matrix')

\vspace{10pt}
{\small $\begin{array}{c}

\hspace{4em}\left[\left(\begin{array}{rrr}
0 & 1 & 0 \\
-1 & 0 & 1 \\
0 & -1 & 0 \\
1 & 0 & 0 \\
0 & 1 & 0 \\
0 & 0 & 1
\end{array}\right), \left(\begin{array}{rrr}
0 & -1 & 0 \\
1 & 0 & 1 \\
0 & -1 & 0 \\
-1 & 1 & 0 \\
0 & 1 & 0 \\
0 & 0 & 1
\end{array}\right), \left(\begin{array}{rrr}
0 & 1 & 0 \\
-1 & 0 & -1 \\
0 & 1 & 0 \\
0 & -1 & 1 \\
1 & -1 & 1 \\
0 & 0 & 1
\end{array}\right),\right.\\[40pt]
\hspace{1em}\left.\left(\begin{array}{rrr}
0 & -1 & 0 \\
1 & 0 & -1 \\
0 & 1 & 0 \\
0 & -1 & 1 \\
-1 & 0 & 1 \\
0 & 0 & 1
\end{array}\right), \left(\begin{array}{rrr}
0 & -1 & 0 \\
1 & 0 & 1 \\
0 & -1 & 0 \\
0 & 0 & -1 \\
-1 & 1 & -1 \\
0 & 1 & -1
\end{array}\right)\right]
\end{array}$
}
\vspace{15pt}

\sagepromt SP3.mutation\_sequence([0,1,0,2],return\_output='var')

$$\left[\frac{x_{1} + y_{0}}{x_{0}}, \frac{x_{0} y_{0} y_{1} + x_{1} x_{2} + x_{2} y_{0}}{x_{0} x_{1}}, \frac{x_{0} y_{1} + x_{2}}{x_{1}}, \frac{x_{0} x_{1} y_{0} y_{1} y_{2} + x_{0} y_{0} y_{1} + x_{1} x_{2} + x_{2} y_{0}}{x_{0} x_{1} x_{2}}\right]$$
}

\noindent A few words about this procedure.  

\begin{enumerate}

\item The command {\tt SP3.mutation\_sequence()} does not affect the object {\tt SP3}, only returns the results of mutating in this order.  If one wants actual seeds to work with rather than simply an output of matrices or cluster variables, one should use the option {\tt return\_output='seed'} (or omit this optional parameter since this is the default setting). 
\sageex{

\sagepromt seeds3 = SP3.mutation\_sequence([0,1,0,2]); seeds3

\sageret{
  \hspace*{-5em}[A seed for a cluster algebra of rank 3 with 3 frozen variables, \\
  \hspace*{-4em}A seed for a cluster algebra of rank 3 with 3 frozen variables, \\
  \hspace*{-3em}A seed for a cluster algebra of rank 3 with 3 frozen variables, \\
  \hspace*{-2em}A seed for a cluster algebra of rank 3 with 3 frozen variables, \\
  \hspace*{-1em}A seed for a cluster algebra of rank 3 with 3 frozen variables]
}}

\item With the optional parameters for returning output, the other options are {\tt 'matrix'} or {\tt 'var'}.  The option {\tt matrix} is self-explanatory.  The option {\tt var} outputs the new cluster variable at each step.  The rest of the cluster variables in the associated clusters are suppressed, since otherwise a lot of redundant information would be printed or saved.
\end{enumerate}

To return the rank (i.e. the number of exchangeable variables or columns in the exchange matrix $B$), one can simply use the command {\tt SP3.n()}.
The set of exchangeable variables can also be obtained by {\tt SP3.exchangeable\_variables()}.  
To return the number of coefficients or frozen variables (also equal to the number of rows minus the number of columns in $B$), we use the command {\tt SP3.m()}.  The frozen variables are obtained by {\tt SP3.frozen\_variables()}.

Not surprisingly, if we mutate {\tt SP3} in place with the same sequence, it equals the last seed returned in the sequence.  

\sageex{
\sagepromt SP3.mutate([0,1,0,2]); SP3.cluster()

$$\left[\frac{x_{0} y_{1} + x_{2}}{x_{1}}, \frac{x_{0} y_{0} y_{1} + x_{1} x_{2} + x_{2} y_{0}}{x_{0} x_{1}}, \frac{x_{0} x_{1} y_{0} y_{1} y_{2} + x_{0} y_{0} y_{1} + x_{1} x_{2} + x_{2} y_{0}}{x_{0} x_{1} x_{2}}, y_{0}, y_{1}, y_{2}\right]$$

\sagepromt SP3 == seeds3[len(seeds3)-1]

\sageret{True}
}

\noindent Notice that it is because of \sage's indexing starting at zero that the last seed is indexed by {\tt len(seeds3)-1}, where {\tt len} stands for ``length''. One can also access the last entry by {\tt seeds3[-1]}.

It is also rather simple to strip off the frozen variables and obtain the coefficient-free cluster algebra by the command {\tt SP3.principal\_restriction()}.  Like the command {\tt S3.principal\_extension()}, this does not change the object in place, and only returns a new object where only the top half of the matrix and the first $n$ cluster variables are kept.  

\vspace{0.5em}

{\bf Warning:} Since we did non-trivial mutations before restricting, we have a labeled seed where $B$ is a $3$-by-$3$ matrix but the associated cluster is Laurent polynomials in terms of $y_0$, $y_1$, and $y_2$, as well as $x_0$, $x_1$, and $x_2$.

\sageex{
\sagepromt SPR3 = SP3.principal\_restriction(); SPR3.cluster()

$$\left[\frac{x_{0} y_{1} + x_{2}}{x_{1}}, \frac{x_{0} y_{0} y_{1} + x_{1} x_{2} + x_{2} y_{0}}{x_{0} x_{1}}, \frac{x_{0} x_{1} y_{0} y_{1} y_{2} + x_{0} y_{0} y_{1} + x_{1} x_{2} + x_{2} y_{0}}{x_{0} x_{1} x_{2}}\right]$$
}

\noindent This causes a problem when we try to mutate:

\sageex{
\sagepromt SPR3.mutate(0); SPR3.cluster()

$$\left[\frac{x_{0} y_{0} y_{1} + x_{0} x_{1} + x_{1} x_{2} + x_{2} y_{0}}{x_{0}^{2} y_{1} + x_{0} x_{2}}, \frac{x_{0} y_{0} y_{1} + x_{1} x_{2} + x_{2} y_{0}}{x_{0} x_{1}},\right.$$
$$\left. \hspace*{100pt}\frac{x_{0} x_{1} y_{0} y_{1} y_{2} + x_{0} y_{0} y_{1} + x_{1} x_{2} + x_{2} y_{0}}{x_{0} x_{1} x_{2}}\right]
$$
}

\noindent since the exchange relation should have involved $y_i$'s but was truncated, and so we did not get the expected cancellation.  Compare with

\sageex{
\sagepromt SP3.mutate(0); SP3.cluster()

$$\left[\frac{x_{1} + y_{0}}{x_{0}}, \frac{x_{0} y_{0} y_{1} + x_{1} x_{2} + x_{2} y_{0}}{x_{0} x_{1}}, \frac{x_{0} x_{1} y_{0} y_{1} y_{2} + x_{0} y_{0} y_{1} + x_{1} x_{2} + x_{2} y_{0}}{x_{0} x_{1} x_{2}}, y_{0}, y_{1}, y_{2}\right]
$$
}

\noindent Thus it is {\bf recommended} to run {\tt SPR3.reset\_cluster()} or {\tt SPR3.set\_cluster()} before actually restricting to the exchangeable cluster variables by running {\tt SPR3 = SP3.principal\_restriction()}.

\section{Using quivers as cluster algebra seeds} \label{sec:quivers}

In this section, we introduce a second way to input a cluster algebra seed.  This uses the language of \emph{quivers}, which is a fancy way of saying a directed (or oriented) graph.  The term \emph{quiver} originates in representation theory, where it was introduced by P. Gabriel at the beginning of the seventies.   Gabriel wanted to emphasize the difference between the representation-theoretic and the graph-theoretic aspects of one and the same notion.  For a quick introduction to quiver representations, please see references such as Section 5 of \cite{Keller2}.  An in-depth treatment is given, for example in the book by Assem, Simson, and Skowronski \cite{ASS}.  The theory of quivers is an important one in representation theory, where fundamental questions come from studying the path algebra associated to such a directed graph.    

For our purposes, we mostly use the quivers for bookkeeping purposes and thinking of them simply as directed graphs will be sufficient for most of our applications.  In this package, a class of objects has been included as a placeholder for future development.  For example, it is planned that in future versions of this package, some of the methods for quiver representations, as in preparation by Franco Saliola, will be available from this class as well.  In the meantime, we will define what we need from quiver theory and describe the methods available in the current package as relevant to cluster algebra theory. 

\begin{definition}
A \emph{quiver} $Q$ is a directed graph.  We will only work with quivers on a finite number of vertices and which contains no loops ($1$-cycles) or $2$-cycles.  However, we do allow our quivers to have multiple edges between a pair of vertices, but since there are no $2$-cycles, this means that all edges between two vertices must have the same direction.  In general there is no restriction against oriented cycles on $\geq 3$ vertices.
\end{definition}

\begin{definition} [Constructing an exchange matrix from a quiver]
Given a quiver $Q$ on vertices $v_0$, $v_1$, through $v_{n-1}$, we let $\pm b_{ij}$ denote the number of edges between $v_i$ and $v_j$.  We let this number be positive if the edges are oriented from $v_i$ to $v_j$ and negative otherwise.  We construct $B_Q = \left[b_{ij}\right]$ as the associated $n$-by-$n$ matrix.  
\end{definition}

\begin{definition} [Constructing a \emph{pair-weighted} quiver from an exchange matrix]
To get a quiver $Q_B$ from an $(m+n)$-by-$n$ exchange matrix $B$ is the reverse of the above construction, however, there are two nuances to emphasize.

\begin{enumerate}

\item For a cluster algebra seed to correspond to a quiver, the corresponding matrix $B$ must satisfy $b_{ij} = - b_{ji}$ for all pairs $0 \leq i,j\leq n-1$.  In other words, the top $n$-by-$n$ portion of $B$ must be skew-symmetric, not just skew-symmetrizable.  Since cluster algebras for non-skew-symmetric seeds are also quite prevalent in the literature, our package works with a slight generalization of quivers, which we call \emph{pair-weighted quivers}.  

We do not allow parallel edges in such quivers, and instead, we label each directed edge as an ordered pair $[b_{ij}, b_{ji}]$ such that the associated edge is oriented from $v_i$ to $v_j$.  Consequently, the first entry of each such pair is necessarily positive and the second is negative, but the direction of the edge must also be recorded.  In the case that $b_{ji} = -b_{ij}$, i.e. the case of parallel edges, this label is simplified to be simply the positive number $b_{ij}$.  We also omit the label $b_{ij}=1$ when displaying graphics to make pictures easier to view.  

Note, that this notation differs from that in places such as \cite{ClustII} or \cite{FeSTuII}, but is necessarily for precise computations.  Our notation is inspired by Dlab-Ringel \cite{DR1976} and Dupont-P\'erotin \cite{DupontPerotin}.

\item If $m > 0$, i.e. the matrix has more rows than columns, and for any $b_{ij}$ where $i \geq n$, there is no $b_{ji}$ in the matrix and so we do not have to worry about checking skew-symmetry for such entries.  However, such vertices $v_i$ correspond to a frozen variable and so we designate these vertices accordingly as ``frozen vertices'' to remind the user not to mutate or apply exchanges at such vertices.

\end{enumerate}

\end{definition}

\begin{remark}\label{rem:equivalence}
This immediate connection between quivers and exchange matrices explains why we often consider exchange matrices up to simultaneous row and column permutations: two quivers are considered to be isomorphic if they are isomorphic as unlabeled digraphs, and this corresponds to considering exchange matrices up to simultaneous row and column permutations. The isomorphism reflects the fact that, as the cluster of an initial cluster seed $(\{x_1,\ldots,x_n\},B)$ is invariant under permuting the variable indices, the cluster algebra does not depend on the ordering of the vertices in the corresponding quiver.
\end{remark}

\begin{definition} [Quiver Mutation]
While a quiver $Q$ can be mutated in any of the $n$ directions by constructing the associated exchange matrix $B_Q$, applying $\mu_k$ and then pulling back to the quiver $Q_{\mu_k(B)} = \mu_k(Q)$, there is also a three step process that allows for a a nice visual description of quiver mutation (in the case of skew-symmetric $B$'s).

\begin{enumerate}
 \item Reverse the direction of every oriented edge incident to vertex $v_k$.  Call the resulting quiver $Q'$.

 \item For any $2$-path $v_i \to v_k \to v_j$ that went through $v_k$ in the original quiver $Q$, add a directed edge $v_i \to v_j$ in $Q'$.  In other words, for any pair of vertices, $\{v_i, v_j\}$, if there are $b_{ik}$ parallel edges from $v_i$ to $v_k$ and $b_{kj}$ parallel edges from $v_k$ to $v_j$, then in $Q'$, we add $b_{ik}b_{kj}$ directed edges between $v_i$ and $v_j$.  

 \item In step 2, a $2$-cycle may have been created, so the last step is to pair off and erase any such anti-parallel edges.

\end{enumerate}

It is an easy exercise to see that the definition of matrix mutation $\mu_k(B)$ given in the previous section agrees with mutation of the quiver $Q_B$ at vertex $v_k$.  In the case of a pair-weighted quiver, it is easiest to mutate the associated matrix and then pull-back to a pair-weighted quiver.  
\end{definition}

\begin{definition}[Mutation-equivalence]
   Two quivers $Q_1, Q_2$ are said to be \emph{mutation-equivalent} if one can be obtained from the other by a finite sequence of mutations, i.e., if there exists a finite sequence $i_1,\ldots,i_k$ such that $\mu_{i_k} \circ \dots \circ \mu_{i_1}(Q_1) = Q_2$. The collection of all quivers mutation-equivalent to a given quiver $Q$ is called \emph{mutation class} of $Q$.
\end{definition}

We now describe the numerous ways that a quiver can be constructed in our package.  Firstly, a quiver can be constructed directly from an exchange matrix, or from a cluster seed in multiple ways. 

\sageex{
\sagepromt B3 = matrix([[0,1,0],[-1,0,-1],[0,1,0]])

\sagepromt S3 = ClusterSeed(B3)

\sagepromt Q1 = Quiver(B3)

\sagepromt Q2 = Quiver(S3)

\sagepromt Q3 = S3.quiver()

\sagepromt Q1 == Q2; Q2 == Q3; Q1

\sageret{True \hspace{3em} True \hspace{3em} Quiver on 3 vertices}
}

\noindent There are other possible constructors, such as from a directed graph:

\sageex{
\sagepromt dg = DiGraph()

\sagepromt dg.add\_edges([[0,1],[2,1]])

\sagepromt Q4 = Quiver(dg)

\sagepromt Q1 == Q4

\sageret{True}
}

\noindent {\bf Warning:} If one uses the digraph constructor, one must follow the conventions for that constructor as a \sage\ object, in particular, digraphs do not allow multiple edges by default.  For example, to get a quiver with parallel edges, one might be tempted to type 

\sageex{
\sagepromt dg = DiGraph()

\sagepromt dg.allows\_multiple\_edges()

\sageret{False}

\sagepromt dg.add\_edges([[0,1],[2,1],[2,1]])

\sagepromt dg.edges(labels=False)

\sageret{[(0, 1), (2, 1)]}

\sagepromt Q5 = Quiver(dg)
}

\noindent However, if one then asks 

\sageex{
\sagepromt Q1 == Q5

\sageret{True}
}

\noindent as the multiple copies of edge $v_2 \to v_1$ are ignored.  Instead, one should use the construction

\sageex{
\sagepromt dg = DiGraph()

\sagepromt dg.add\_edges([[0,1,1],[2,1,2]])

\sagepromt Q6 = Quiver(dg)

\sagepromt Q1 == Q6; Q6.digraph().edges()

\sageret{False \hspace{3em} [(0,1,(2,-2)), (2,1,(1,-1))]}
}

\noindent Note that all quivers are actually implemented as pair-weighted quivers, i.e. as a labeled digraph where multiple edges correspond to a pair $(b, -b)$ where $b \geq 2$.  The program automatically converts the user's input with a single number indicating the edge label to a pair.  A user can even label some edges as a single number, leave some edges unlabeled (as a single edge with pair-weight $(1,-1)$), and other edges as pairs; and the program will interpret this correctly.  As mentioned above, multiple copies of an edge are ignored.  More precisely, if they are given as labeled edges, then the label assigned is the one given to the last copy of the edge included.

\sageex{
\sagepromt dg = DiGraph()

\sagepromt dg.add\_edges([[0,1,(1,-1)],[2,1,2]])

\sagepromt Q8 = Quiver(dg)

\sagepromt Q6 == Q8

\sageret{True}
}

\noindent A quiver can also be constructed more quickly by having \sage\ do the intermediate work of constructing the digraph for you.  Just simply type 

\sageex{
\sagepromt Q9 = Quiver([[0,1,1],[2,1,2]])
}

\noindent or any of the analogous constructions described above for encoding the edges of a digraph (although again one should include edge labels instead of multiple copies of edges).

\sageex{
\sagepromt Q6 == Q9

\sageret{True}
}

\noindent You can also get a copy of a quiver already defined by a command such as

\sageex{
\sagepromt Q10 = Quiver(Q9)

\sagepromt Q10 == Q9

\sageret{True}

\sagepromt Q10.mutate(0)

\sagepromt Q10 == Q9

\sageret{False}
}

\noindent We did not emphasize it above, but a similar technique allows one to get a copy of a cluster seed.  There is one other technique that can be used to construct a quiver, or for that matter a cluster seed, but requires knowledge of quiver mutation types, and we leave the description of this construction to the next section. 

We now introduce some of the possible methods our package contains for working with quivers along with associated examples.  Most importantly, to get a picture of a quiver, we use the {\tt show} command.

\sageex{
\sagepromt Q1.show()

\sageret{\includegraphics[width=2in]{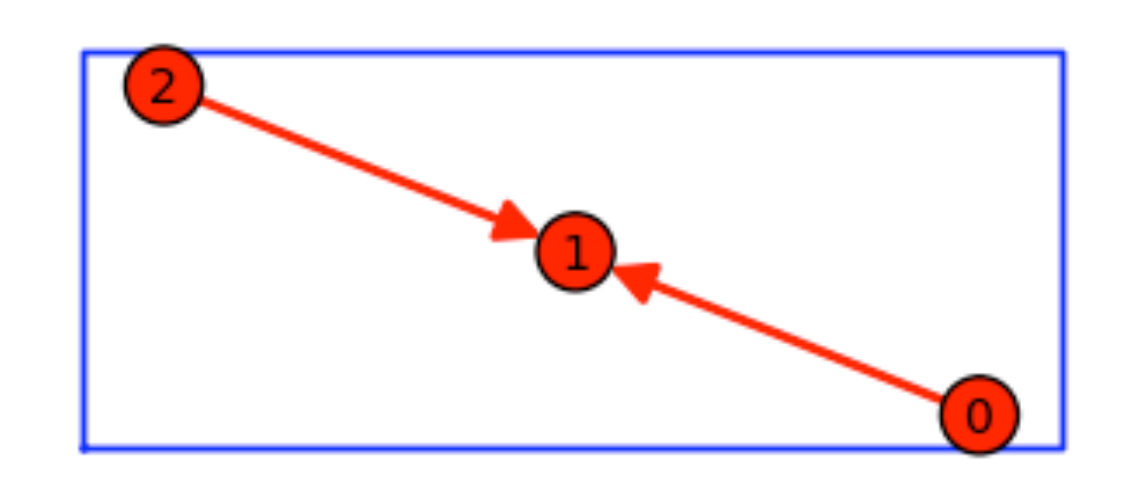}}

\sagepromt Q1.mutate(1); Q1.show()

\sageret{\includegraphics[width=2in]{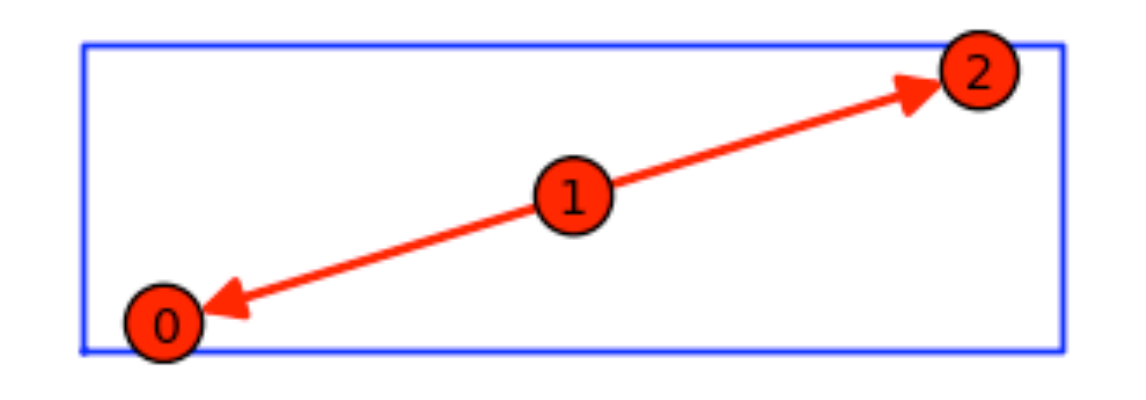}}
}

\noindent One quirk about the method {\tt show} is that the graphic obtained can be a little random.  If the placement of the vertices or the drawing of the graph is not optimal, it is recommended the user try running the {\tt show} command again until the quiver renders in a more visually pleasing way.

Note that just as before, the command {\tt Q1.mutate()} changes the object in place.

\sageex{
\sagepromt Q1.b\_matrix()

$$
\left(\begin{array}{rrr}
0 & -1 & 0 \\
1 & 0 & 1 \\
0 & -1 & 0
\end{array}\right)
$$

\sagepromt Q1.mutate(1); Q1.b\_matrix()

$$
\left(\begin{array}{rrr}
0 & 1 & 0 \\
-1 & 0 & -1 \\
0 & 1 & 0
\end{array}\right)
$$
}

\noindent A nice way to visualize a sequence of quiver mutations is the use of the command {\tt mutation\_sequence} with the optional parameter {\tt show\_sequence=True}.  Unlike the {\tt show} command, the quivers always render in the same circular way using this procedure so it is easier to compare vertices to one another.

\sageex{
\sagepromt Q1.mutation\_sequence([0,1,2,0,1,0],show\_sequence=True)

\sageret{\hspace{-12em}[Quiver on 3 vertices,\\ \hspace{-9em}Quiver on 3 vertices,\\ \hspace{-6em}Quiver on 3 vertices,\\ \hspace{-3em}Quiver on 3 vertices,\\ \hspace{0em}Quiver on 3 vertices,\\ \hspace{3em}Quiver on 3 vertices,\\ \hspace{6em}Quiver on 3 vertices]}

\sageret{\includegraphics[width=5in]{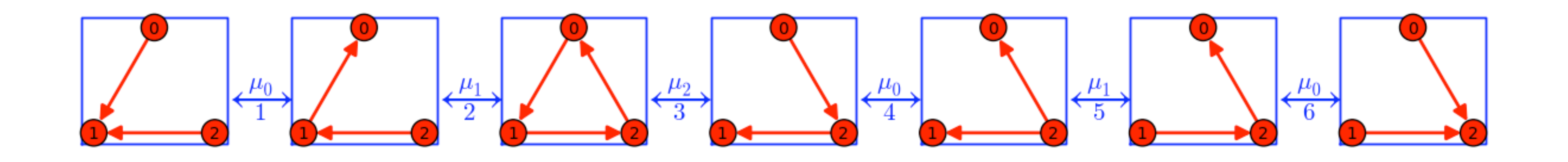}}
}

\noindent Note that, here we are using the {\tt mutation\_sequence} on the quivers (rather than seeds), so the optional argument of {\tt return\_output} is not allowed.  However,  we can construct the associated cluster seed quite easily and then methods for viewing the associated quiver are still accessible, along with the other commands for cluster seeds.

\sageex{
\sagepromt NewS = ClusterSeed(Q1); Q2 = NewS.quiver()

\sagepromt Q2 == Q1; NewS.show()

\sageret{True}

\sageret{\includegraphics[width=1.5in]{tmp_1.pdf}}

\sagepromt NewS.mutation\_sequence([0,1,2,0,1,0],\\
\hspace*{6em}show\_sequence=True,return\_output='matrix')

\vspace*{10pt}
\centering
{\tiny
$\begin{array}{c} \hspace{-2em} \left[\left(\begin{array}{rrr}
0 & 1 & 0 \\
-1 & 0 & -1 \\
0 & 1 & 0
\end{array}\right)\right., \left(\begin{array}{rrr}
0 & -1 & 0 \\
1 & 0 & -1 \\
0 & 1 & 0
\end{array}\right), \left(\begin{array}{rrr}
0 & 1 & -1 \\
-1 & 0 & 1 \\
1 & -1 & 0
\end{array}\right), \left(\begin{array}{rrr}
0 & 0 & 1 \\
0 & 0 & -1 \\
-1 & 1 & 0
\end{array}\right), \\[15pt] \hspace*{15em} \left[\left(\begin{array}{rrr}
0 & 0 & -1 \\
0 & 0 & -1 \\
1 & 1 & 0
\end{array}\right), \left(\begin{array}{rrr}
0 & 0 & -1 \\
0 & 0 & 1 \\
1 & -1 & 0
\end{array}\right), \left(\begin{array}{rrr}
0 & 0 & 1 \\
0 & 0 & 1 \\
-1 & -1 & 0
\end{array}\right)\right]\end{array}$
}

\sageret{\includegraphics[width=5in]{tmp_3.pdf}}
}
\noindent Another instructive command is {\tt DiGraph} which lets the user construct the associated labeled directed graph encoding the quiver.  Since {\tt DiGraph} is already a class of objects, this allows the user access to a variety of other methods.  One can then reconstruct a quiver with the altered directed graph, {\tt dg} whenever desired, using the techniques described above, i.e. {\tt Quiver(dg)}.

\sageex{
\sagepromt Quivs = Q1.mutation\_sequence([0,1,2,0,1,0])

\sagepromt [Q.digraph().edges() for Q in Quivs]

\sageret{
  \hspace{-12em}[[(0, 1, (1, -1)), (2, 1, (1, -1))],\\
  \hspace{-9em}[(1, 0, (1, -1)), (2, 1, (1, -1))],\\
  \hspace{3em}[(0, 1, (1, -1)), (1, 2, (1, -1)), (2, 0, (1, -1))],\\
  \hspace{-3em}[(0, 2, (1, -1)), (2, 1, (1, -1))],\\
  \hspace{0em}[(2, 0, (1, -1)), (2, 1, (1, -1))],\\
  \hspace{3em}[(1, 2, (1, -1)), (2, 0, (1, -1))],\\
  \hspace{6em}[(0, 2, (1, -1)), (1, 2, (1, -1))]]
}
}

\noindent Thus far, the examples included have been skew-symmetric and coefficient-free.  We close this section with some examples which require pair-weighted quivers and frozen vertices.

\sageex{
\sagepromt B = matrix([\\
\hspace*{6em}(0, 1, 0, 0, 0, 1), \\
\hspace*{7em}(-1, 0, -1, 0, 0, 0), \\
\hspace*{8em}(0, 1, 0, 1, 0, 0), \\
\hspace*{9em}(0, 0, -1, 0, -2, 0), \\
\hspace*{10em}(0, 0, 0, 1, 0, 0), \\
\hspace*{11em}(-2, 0, 0, 0, 0, 0)])

\sagepromt S = ClusterSeed(B); S.show()

\sageret{\includegraphics[width=3in]{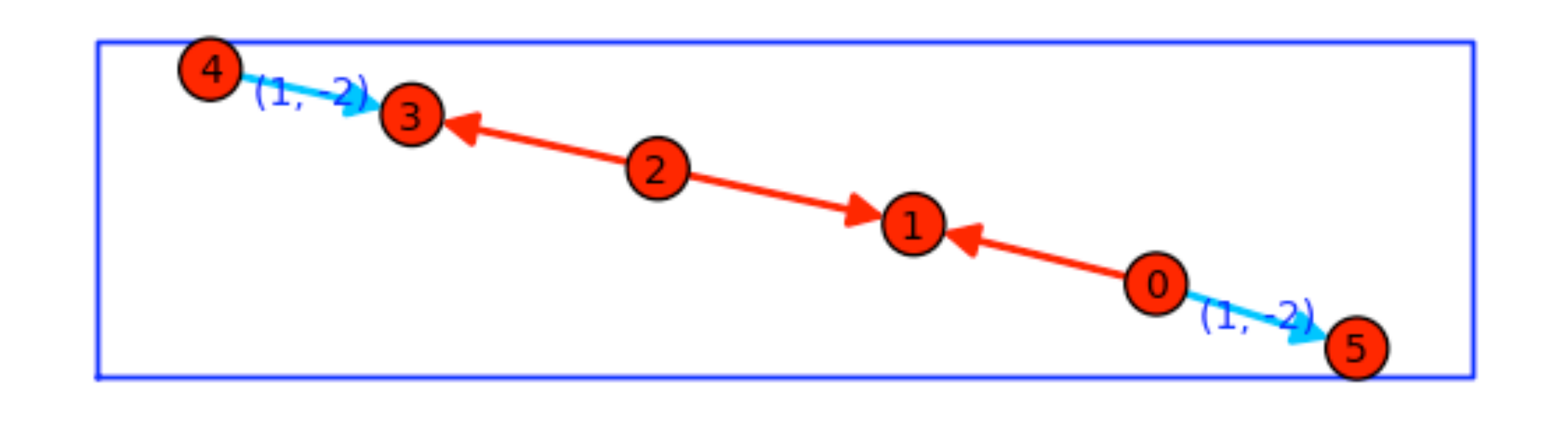}}

\sagepromt S.mutation\_sequence([1,2,0], show\_sequence=True, \\
\hspace*{6em}return\_output='matrix')

{\centering
\vspace*{5pt}
{\tiny
$\begin{array}{c}\hspace{-4em}
\left[\left(\begin{array}{rrrrrr}
0 & 1 & 0 & 0 & 0 & 1 \\
-1 & 0 & -1 & 0 & 0 & 0 \\
0 & 1 & 0 & 1 & 0 & 0 \\
0 & 0 & -1 & 0 & -2 & 0 \\
0 & 0 & 0 & 1 & 0 & 0 \\
-2 & 0 & 0 & 0 & 0 & 0
\end{array}\right)\right., \left(\begin{array}{rrrrrr}
0 & -1 & 0 & 0 & 0 & 1 \\
1 & 0 & 1 & 0 & 0 & 0 \\
0 & -1 & 0 & 1 & 0 & 0 \\
0 & 0 & -1 & 0 & -2 & 0 \\
0 & 0 & 0 & 1 & 0 & 0 \\
-2 & 0 & 0 & 0 & 0 & 0
\end{array}\right), \\[20pt] \hspace*{6em}\left(\begin{array}{rrrrrr}
0 & -1 & 0 & 0 & 0 & 1 \\
1 & 0 & -1 & 1 & 0 & 0 \\
0 & 1 & 0 & -1 & 0 & 0 \\
0 & -1 & 1 & 0 & -2 & 0 \\
0 & 0 & 0 & 1 & 0 & 0 \\
-2 & 0 & 0 & 0 & 0 & 0
\end{array}\right), \left. \left(\begin{array}{rrrrrr}
0 & 1 & 0 & 0 & 0 & -1 \\
-1 & 0 & -1 & 1 & 0 & 1 \\
0 & 1 & 0 & -1 & 0 & 0 \\
0 & -1 & 1 & 0 & -2 & 0 \\
0 & 0 & 0 & 1 & 0 & 0 \\
2 & -2 & 0 & 0 & 0 & 0
\end{array}\right)\right]
\end{array}$}

\sageret{\includegraphics[width=5in]{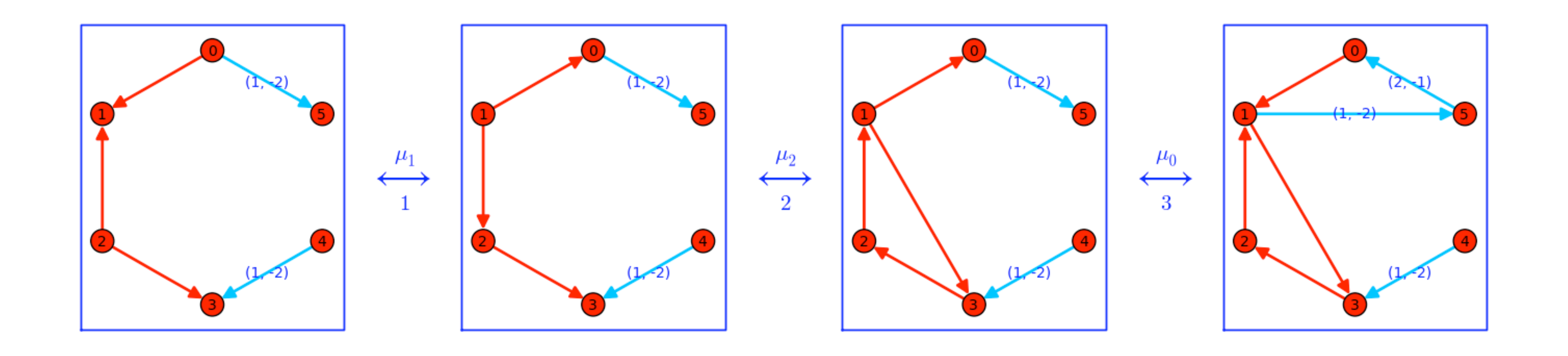}}
}

\sagepromt Q = Quiver(S)

\sagepromt Q2 = Q.principal\_extension(); Q2; Q2.show()

\sageret{Quiver on 6 vertices with 6 frozen vertices}

\sageret{\includegraphics[width=2.5in]{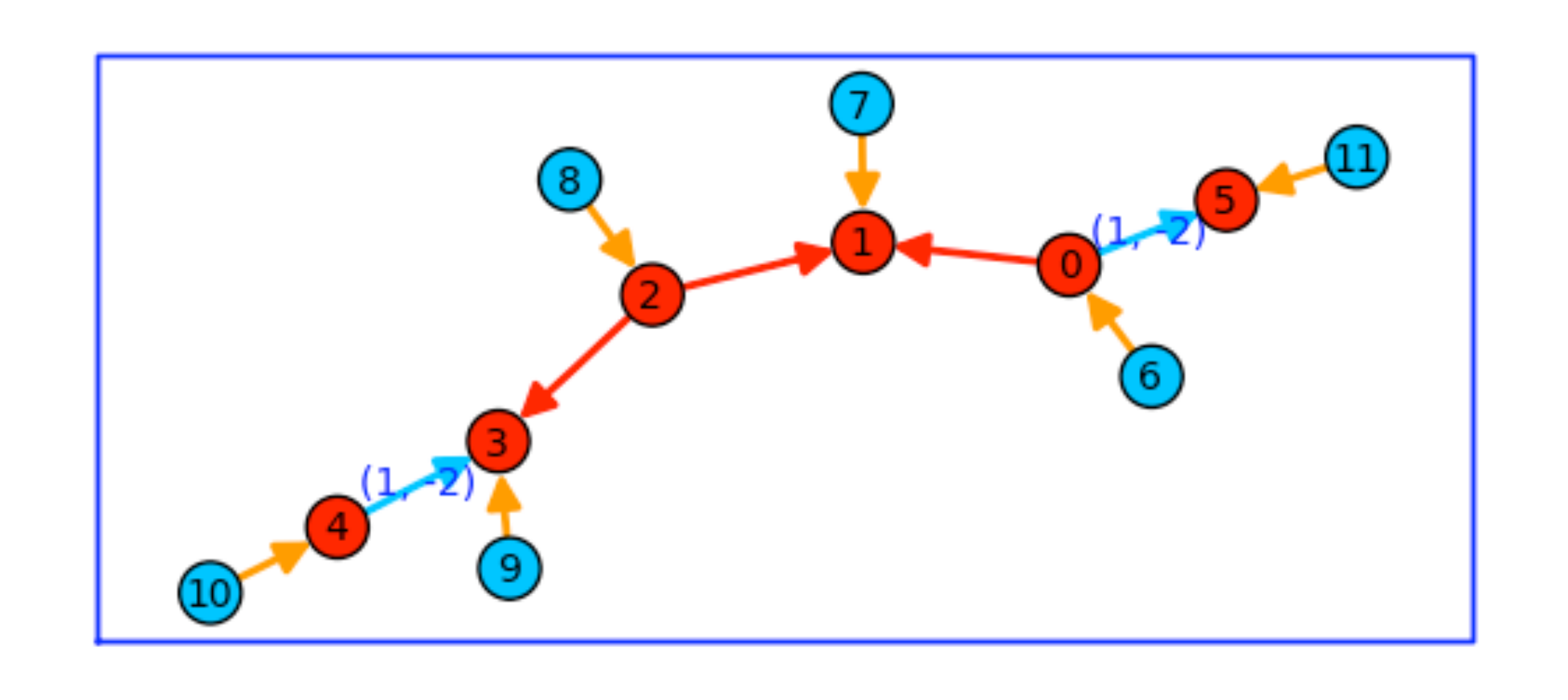}}

\sagepromt Q2.mutation\_sequence([1,2,0],show\_sequence=True)

\sageret{\includegraphics[width=5in]{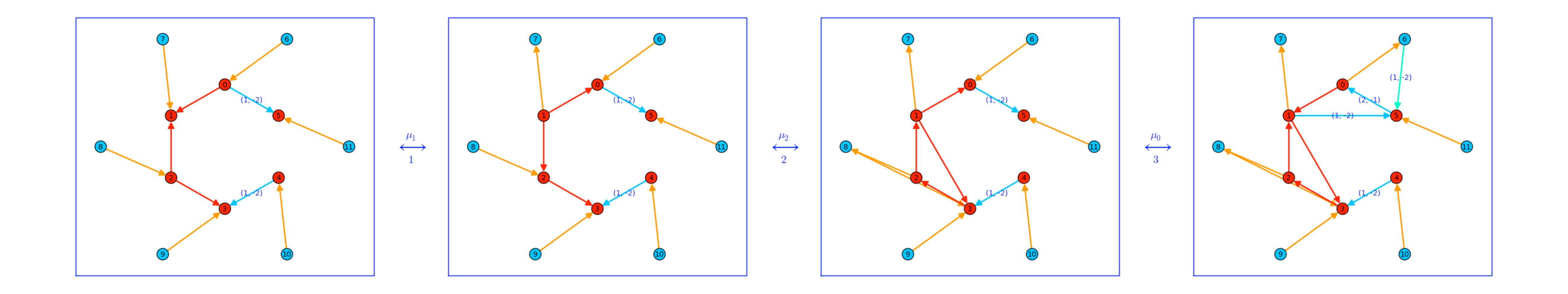}}
}

\noindent If we instead produce a quiver by first producing the {\tt principal\_extension} of the cluster seed, and then constructing a quiver from it, we obtain an equal quiver as a result.

\sageex{
\sagepromt S2 = S.principal\_extension()

\sagepromt Q3 = Quiver(S2); Q2 == Q3

\sageret{True}
}

\begin{figure}
    \centering
    \includegraphics[width=5in]{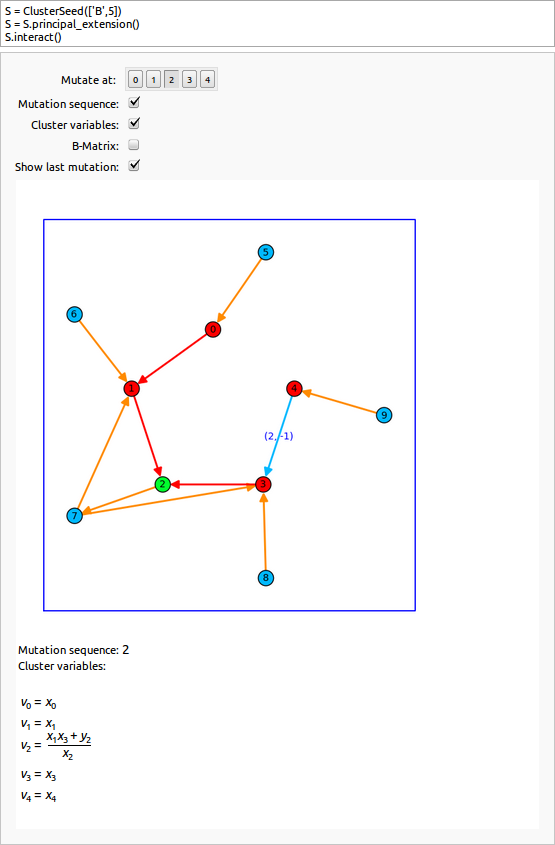}
    \caption{The interactive mode of the cluster package in the {\tt sage-notebook}.}
    \label{fig:interactivemode}
  \end{figure}

Another way to work with quivers and cluster seeds is through the interactive mode available through the {\tt sage-notebook}.  This involves a command such as {\tt S.interact()} or {\tt Q.interact()}, as shown in Figure~\ref{fig:interactivemode}.

  \section{Finite type and finite mutation type classifications} \label{sec:typeclassification}

So far we have described how a cluster algebra seed can be constructed from a skew-symmetrizable matrix or from a quiver.  The last construction that we wish to discuss utilizes the notion of \emph{quiver mutation types}.  Before, we delve more into the specifics of this discussion, we begin with a few theoretical preliminaries.

Two natural questions that one can ask about a cluster algebra (or its seed) once the initial definitions have been given are the following:

\begin{itemize}

\item[{}] \noindent \hspace{-3.5em} Given a cluster algebra $\mathcal{A} = \mathcal{A}(\mathbf{x}_0,B_0)$, with initial seed $(\mathbf{x}_0,B_0)$,
 
\item[$\bullet$] are there a finite number of generators (cluster variables) $x$ for $\mathcal{A}$ as we take the union of all clusters $\mathbf{x}$ as we mutate?

\item[$\bullet$] are there a finite number of exchange matrices $B$ for $\mathcal{A}$ as we mutate into different seeds?
\end{itemize}

\begin{definition}
If there are a finite number of cluster variables for $\mathcal{A}$, we say that $\mathcal{A}$ is of \emph{finite type}.
\end{definition}

\begin{definition}
If there are a finite number of exchange matrices for $\mathcal{A}$, we say that $\mathcal{A}$ is of \emph{finite mutation type}.
\end{definition}

An important theorem that greatly simplifies our notation for geometric type is the following theorem by Gekhtman, Shapiro, and Vainshtein:

\begin{theorem} [Theorem 7.4 of \cite{GSV2}]
If $\mathcal{A}$ is a cluster algebra (i) of geometric type, or (ii) has nondegenerate exchange matrix, and $(\mathbf{x},B)$, $(\mathbf{x'},B')$ are two seeds for $\mathcal{A}$, such that cluster $\mathbf{x'}$ is simply the permutation $\sigma$ of cluster $\mathbf{x}$, then the exchange matrices $B'$ and $B$ must also be the same, up to simultaneous permutation of its rows and columns by the same $\sigma$.  In particular, the \emph{cluster determines the seed} in the above cases. 
\end{theorem}

From this theorem, it is clear that any cluster algebra of finite type must have a finite number of clusters, hence a finite number of seeds and exchange matrices.

\begin{corollary}[Finite type implies finite mutation type]
  A cluster algebra of finite type \emph{is also} of finite mutation type.
\end{corollary}
However, the converse is false, the simplest counter-example being the rank two example $\mathcal{A}(2,2)$ discussed in the Introduction.

Classifying cluster algebras of finite type was one of the first natural questions about cluster algebras, and led Fomin and Zelevinsky to the following beautiful theorem.

\begin{theorem} [Theorem 1.5 of \cite{ClustII}] \label{th:finitetypeclassification}
The following three conditions about a cluster algebra $\mathcal{A} = \mathcal{A}(\mathbf{x}_0,B_0)$ are equivalent:

\begin{itemize}
\item Cluster algebra $\mathcal{A}$ is of finite type.
\item In every seed $(\mathbf{x},B)$ that is mutation-equivalent to $(\mathbf{x}_0,B_0)$, the exchange matrix $B$ satisfies $|b_{ij} b_{ji}| \leq 3$ for all pairs $1 \leq i, j \leq n$.
\item There exists a mutation-equivalent seed $(\mathbf{x}_1,B_1)$ such that the exchange matrix $B_1$ is a \emph{skew-symmetric version} of a Cartan matrix of a finite-dimensional Lie algebra\footnote{
Given a Cartan matrix $A$, we make a skew-symmetric $B_A$ by replacing the $2$'s on the diagonal with $0$'s, and picking a bipartite coloring of the Dynkin diagram associated to $A$ so that $b_{ij} = |a_{ij}|$ if directed edge $v_i \to v_j$ would go from white to black, and $b_{ij} = - |a_{ij}|$ otherwise, see Section~\ref{sec:associahedra}.}.
\end{itemize}
In particular, cluster algebras of finite type are given by the same Cartan-Killing classification as that describing Lie algebras via Dynkin diagrams: $$A_n, ~B_n,~ C_n,~ D_n,~ E_6,~ E_7,~ E_8,~ F_4,~ \mathrm{and}~ G_2.$$  
\end{theorem}

Given a cluster algebra seed $S$ for $\mathcal{A}$, it therefore makes sense to ask whether or not $S$ is mutation-equivalent to a seed $(\mathbf{x},B)$ where the exchange matrix $B$ is a skew-symmetric version of the Cartan matrix of type $A_n$ (resp. $B_n$, $C_n$, $D_n$, $E_6$, $E_7$, $E_8$, $F_4$, or $G_2$).  If so, we call $\mathcal{A}$ a cluster algebra of mutation type $A_n$ (resp. $B_n$, $C_n$, $D_n$, $E_6$, $E_7$, $E_8$, $F_4$, or $G_2$).  We also call all exchange matrices and the corresponding quivers of such a cluster algebra of mutation type $A_n$ (resp. $B_n$, $C_n$, $D_n$, $E_6$, $E_7$, $E_8$, $F_4$, or $G_2$).

Our program has algorithms for identifying mutation types of exchange matrices and quivers.  In the cases of the exceptional types, $E_6$, $E_7$, $E_8$, $F_4$ and $G_2$, it is sufficient to hard-code a catalog of the mutation classes. This is done to avoid recomputing the mutation class whenever checking a mutation type. In classical types however, the parameter $n$ can be any positive integer, and we instead utilize theoretical results of \cite{SchifflerA} (type $A_n$), \cite{Stump-pre} (types $B_n$ and $C_n$), 
and \cite{Vatne} (type $D_n$) to identify them for any rank $n$. 

Recall that a quiver (resp. pair-weighted quiver) encodes the same information as a skew-symmetric (resp. skew-symmetrizable) matrix.  To avoid duplication of data types, we have introduced a new class of objects known as quiver mutation types.  Note that these can be implemented with or without brackets.

\sageex{
\sagepromt QM1 = QuiverMutationType(['A',5])

\sagepromt QM2 = QuiverMutationType('A',5); QM1 == QM2

\sageret{True}

\sagepromt QM1

\sageret{['A', 5]}

\sagepromt type(QM1)
}

\sageret{$\langle$class 'sage.all\_cmdline.QuiverMutationType\_Irreducible'$\rangle$}

\sageex{
\sagepromt QM1.b\_matrix()

$$
\left(\begin{array}{rrrrr}
0 & 1 & 0 & 0 & 0 \\
-1 & 0 & -1 & 0 & 0 \\
0 & 1 & 0 & 1 & 0 \\
0 & 0 & -1 & 0 & -1 \\
0 & 0 & 0 & 1 & 0
\end{array}\right)
$$

\sagepromt Quiv = QM1.standard\_quiver(); Quiv

\sageret{Quiver on 5 vertices of type ['A', 5]}

\sagepromt Quiv.show()

\sageret{\includegraphics[width=3in]{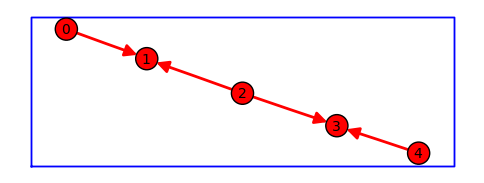}}

\sagepromt QM2 = QuiverMutationType('BC',6,1); QM2

\sageret{['BC',5,1]}

\sagepromt QM2.b\_matrix()

$$\left(\begin{array}{rrrrrrr}
0 & 1 & 0 & 0 & 0 & 0 & 1 \\
-1 & 0 & -1 & 0 & 0 & 0 & 0 \\
0 & 1 & 0 & 1 & 0 & 0 & 0 \\
0 & 0 & -1 & 0 & -1 & 0 & 0 \\
0 & 0 & 0 & 1 & 0 & 2 & 0 \\
0 & 0 & 0 & 0 & -1 & 0 & 0 \\
-2 & 0 & 0 & 0 & 0 & 0 & 0
\end{array}\right)
$$

\sagepromt QM2.standard\_quiver().show()

\sageret{\includegraphics[width=4in]{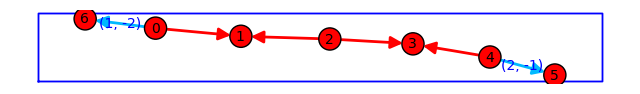}}

}

\noindent Each quiver mutation type has a number of attributes and methods associated to it.  We already saw an example of two key methods: {\tt b\_matrix} and {\tt standard\_quiver}, i.e.
each quiver mutation type object encodes a specific canonical exchange matrix and the associated pair-weighted quiver.  This characterizes only one representative out of the relevant possible mutation-class, but it is enough data to determine the appropriate cluster algebra seed up to mutation-equivalence.  We hard-coded these representatives so that the associated quiver is an oriented Dynkin diagram such that each vertex is a sink or source.  For future reference, such a quiver and seed is known as \emph{bipartite}.

More generally, each of these representative quivers are \emph{trees} and \emph{acyclic}.  Because of results from representation theory and otherwise, there are a number of results in cluster algebra theory that hold when the associated quiver is bipartite (resp. a tree or acyclic), but the result is incorrect, or a proof is unknown when the quiver lacks the relevant property.  Here are some examples:

\begin{theorem} \cite{Nak}
If a cluster algebra $\mathcal{A}$ is given by a seed that is mutation-equivalent to one which is skew-symmetric and bipartite, then all cluster variables of $\mathcal{A}$ have positive expansions as Laurent polynomials
\footnote{By theorems of Fan Qin \cite{FQ} and an updated version of \cite{Nak}, positivity has been proven for all skew-symmetric acyclic seeds.}.  
\end{theorem}

\begin{theorem} [Proposition 9.2 in \cite{ClustII}] \label{th:tree} If $Q$ is a quiver that is a tree as an undirected graph then $Q$ is mutation-equivalent to any $Q'$ where $Q'$ has the same underlying undirected graph as $Q$ but the edges of $Q'$ are oriented arbitrarily\footnote{Note: this list of mutation-equivalent quivers is not exhaustive, for example a quiver of type $A_3$ is both mutation-equivalent to any orientation of a path on three vertices; or to an oriented triangle.}.
\end{theorem}    

\begin{theorem} [Corollary 1.21 in \cite{ClustIII}] \label{th:exchangerelation} Let $\mathcal{A} = \mathcal{A}(\mathbf{x}, B)$ be a cluster algebra where $B$ corresponds to an acyclic seed.  Let $x_i'$ denote the unique element in cluster $\mu_i(\mathbf{x})$ which is not contained in $\mathbf{x}$.  Then we have the following:
\begin{itemize}
\item[$\bullet$] $\mathcal{A}$ is finitely generated by the set $\chi = \{x_1, x_1', \dots, x_n, x_n'\}$, 
\item[$\bullet$] The standard monomials (those not containing the factor $x_i x_i'$ for any $i\in \{0,1,\dots, n-1\}$) in $\chi$ form a $\mathbb{ZP}$-basis of $\mathcal{A}$, and 
\item[$\bullet$] The binomial exchange relations involving $x_i x_i'$ on the left-hand-sides generate the ideal of relations among the generators $\chi$.
\end{itemize}
\end{theorem}

Because of the importance of these properties, and other related ones, there are methods to check whether a given cluster seed, quiver, or quiver mutation type satisfies them:

\sageex{\sageret{is\_finite(), is\_mutation\_finite(), is\_bipartite(), is\_acyclic(),\dots}}

\noindent There are a few other checks that we have not explained yet, but we will provide an annotated list of all of the checkable properties in Section \ref{glossary}.

\sageex{
\sagepromt QM1.properties()

\sageret{
['A', 5] has rank 5 and the following properties:\\
	- irreducible:       True\\
	- mutation finite:   True\\
	- simply-laced:      True\\
	- skew-symmetric:    True\\
	- finite:            True\\
	- affine:            False\\
	- elliptic:          False
}

\sagepromt QM2.properties()

\sageret{
['BC', 6, 1] has rank 7 and the following properties:\\
	- irreducible:       True\\
	- mutation finite:   True\\
	- simply-laced:      False\\
	- skew-symmetric:    False\\
	- finite:            False\\
	- affine:            True\\
	- elliptic:          False
}
}

\noindent Most importantly, our program allows the user to construct a cluster seed or quiver by using a quiver mutation type. The associated quiver is the standard quiver that is hard-coded as a representative for each type; and the associated cluster seed is obtained from this choice of quiver.

\sageex{
\sagepromt ClusterSeed(['A',5])

\sageret{
A seed for a cluster algebra of rank 5 of type ['A', 5]
}

\sagepromt ClusterSeed(['BC',6,1])

\sageret{
A seed for a cluster algebra of rank 7 of type ['BC', 6, 1]
}

\sagepromt Quiver(['A',5])

\sageret{
Quiver on 5 vertices of type ['A', 5]
}

\sagepromt Quiver(['BC',6,1])

\sageret{
Quiver on 7 vertices of type ['BC', 6, 1]
}
}

\subsection{Finite mutation type classification}

We now describe theoretical results regarding the classification of cluster algebras of finite mutation type.  Again, we use the notation of pair-weighted quivers so our descriptions of some of the results will differ slightly from the work of Felikson-Shapiro-Tumarkin \cite{FeSTuII}.  Our story begins however with Felikson-Shaprio-Tumarkin's first paper \cite{FeSTu} which classified \emph{skew-symmetric} cluster algebras of finite mutation type.

\begin{theorem} [Theorem 6.1 of \cite{FeSTu}]\label{th:skew-symmetric}
The following two conditions about a cluster algebra $\mathcal{A} = \mathcal{A}(\mathbf{x}_0,B_0)$ with skew-symmetric $B_0$ are equivalent:
\begin{itemize}
 \item $\mathcal{A}$ is of finite mutation type,
 \item $\mathcal{A}$ has one of the following properties:
  \begin{enumerate}
   \item $\mathcal{A}$ is of rank $2$,
   \item $\mathcal{A}$ is associated to a \emph{cluster algebra corresponding to a surface}, or
   \item $\mathcal{A}$ is one of $11$ exceptional types $E_6$, $E_7$, $E_8$, affine $\tilde{E}_6$, $\tilde E_7$, $\tilde E_8$, elliptic $\tilde E_6^{(1)}$, $\tilde E_7^{(1)}$, $\tilde E_8^{(1)}$, or one of two other types $X_6$ and $X_7$, which were found by Derksen and Owen \cite{DerkOwen}.
  \end{enumerate}
\end{itemize}
\end{theorem} 

Rank two cluster algebras were already described in the introduction, and are clearly mutation-finite since mutation of such an exchange matrix $B$ simply leads to $-B$.

Describing cluster algebras of surfaces is beyond the scope of this compendium, however it is planned that future installments of this software will handle such cluster algebras and their description will be spelled out at that time.  Please see Fomin, Shaprio, and D. Thurston's papers \cite{FST, FT} for a description or \cite{MSW} where Schiffler, Williams, and the first author prove positivity of Laurent expansions for such cluster algebras.  Nonetheless, we mention here that cluster algebras corresponding to \emph{polygons with $0$, $1$, or $2$ punctures, or to an annulus}, can also be described as the skew-symmetric types $A_n$, $D_n$, $\tilde{D}_n$, or $\tilde{A}_{r,s}$, respectively.  The first two cases are of finite type and the second two are of affine type.  Any other finite or affine type is of exceptional type or is not skew-symmetric.  We illustrate corresponding representative quivers in the next section. 

We have met some of the eleven exceptional types before, the types $E_6$, $E_7$, and $E_8$ are of finite type and thus of finite mutation type.  We give representative quivers for the remaining eight in the next section.  The affine types $\tilde{E}_6$, $\tilde{E}_7$, and $\tilde{E}_8$ each have a bipartitely oriented tree as a quiver representative; however the other five have no acyclic representatives.    

\subsection{Skew-symmetrizable cluster algebra seeds of finite mutation type} \label{sec:Kacclassification}

In cutting edge work this summer \cite{FeSTuII}, Felikson-Shapiro-Tumarkin generalized their previous work to a classification including mutation-finite weighted quivers that are not skew-symmetric.  

\begin{theorem} [Theorems~2.8 and 5.13 of \cite{FeSTuII}]
The following three conditions about a cluster algebra $\mathcal{A} = \mathcal{A}(\mathbf{x}_0,B_0)$ with skew-symmetrizable $B_0$ are equivalent:
\begin{itemize}
 \item $\mathcal{A}$ is of finite mutation type,
 \item  In every seed $(\mathbf{x},B)$ that is mutation-equivalent to $(\mathbf{x}_0,B_0)$, the exchange matrix $B$ satisfies $|b_{ij} b_{ji}| \leq 4$ for all pairs $1 \leq i, j \leq n$.
 \item $\mathcal{A}$ has one of the following properties:
  \begin{enumerate}
   \item $\mathcal{A}$ is of rank $2$,
   \item $\mathcal{A}$ is decomposable into \emph{blocks}, as described in \cite{FeSTuII}, or
   \item $\mathcal{A}$ is one of the $11$ exceptional types in Theorem~\ref{th:skew-symmetric} or one of the $7$ exceptional types $\tilde G_2, F_4, \tilde F_4, V_4, W_4, Y_4,$ and $Z_6$.
  \end{enumerate}
\end{itemize}
\end{theorem} 

\begin{remark}
One can get from our notation of pair-weighted quivers to the notion of weighted quivers in \cite{FeSTuII} by the following: if an edge of our quiver has the pair-weight $[b, -c]$, then the corresponding weight in their notation is $bc$.  While their notation has several advantages and simplifies the statements of certain theorems, for computations it obscures the differences between different mutation classes.  For example, cluster algebras of types $B_n$ and $C_n$ would have the same weighted quivers.  Even though these cluster algebras give rise to the same cluster complexes (i.e. the clique complex induced by the graph whose vertices are seeds and whose edges are mutations), the Laurent expansions of cluster variables are quite different in these two cases.
\end{remark}  

\noindent To illustrate this example we introduce two new commands.  See Section \ref{sec:varclass} for details on the associated algorithms:

1) Given a cluster algebra of finite mutation type, we can use the command {\tt b\_matrix\_class} to obtain a list of all the exchange matrices that are mutation-equivalent to a given initial seed.  To avoid extraneous duplication, we only output one matrix up to simultaneous permutation of rows and columns.  

For example, in the $B_3$ versus $C_3$ cases, notice that the list of exchange matrices in the respective mutation classes are negative transposes of one another\footnote{This would be clearer if we included all mutation-equivalent matrices rather than just those up to permutation, which could be accomplished by {\tt S3.b\_matrix\_class(up\_to\_equivalence=False)}.  In particular the last matrices in both of these lists are negative transposes of each other if we also swap the first and second rows/columns.}.

\sageex{
\sagepromt S3 = ClusterSeed(['B',3]); S3.b\_matrix\_class()

\vspace*{10pt}
$\begin{array}{c}\left[\left(\begin{array}{rrr}
0 & 0 & 1 \\
0 & 0 & 2 \\
-1 & -1 & 0
\end{array}\right) \right., \left(\begin{array}{rrr}
0 & 0 & 1 \\
0 & 0 & -2 \\
-1 & 1 & 0
\end{array}\right), \left(\begin{array}{rrr}
0 & 1 & 1 \\
-2 & 0 & 0 \\
-1 & 0 & 0
\end{array}\right), \\[20pt] \hspace{15em}\left(\begin{array}{rrr}
0 & 2 & 0 \\
-1 & 0 & 1 \\
0 & -1 & 0
\end{array}\right), \left. \left(\begin{array}{rrr}
0 & -1 & 1 \\
2 & 0 & -2 \\
-1 & 1 & 0
\end{array}\right)\right]
\end{array}$

\sagepromt S4 = ClusterSeed(['C',3]); S4.b\_matrix\_class()

$\begin{array}{c}\left[\left(\begin{array}{rrr}
0 & 0 & 1 \\
0 & 0 & 1 \\
-1 & -2 & 0
\end{array}\right) \right., \left(\begin{array}{rrr}
0 & 0 & 1 \\
0 & 0 & -1 \\
-1 & 2 & 0
\end{array}\right), \left(\begin{array}{rrr}
0 & 2 & 1 \\
-1 & 0 & 0 \\
-1 & 0 & 0
\end{array}\right), \\[20pt] \hspace{15em}\left(\begin{array}{rrr}
0 & 1 & 0 \\
-2 & 0 & 1 \\
0 & -1 & 0
\end{array}\right), \left. \left(\begin{array}{rrr}
0 & 1 & -1 \\
-2 & 0 & 1 \\
2 & -1 & 0
\end{array}\right)\right]
\end{array}$
}
\sageex{

\sagepromt S3.show(); S4.show()

\sageret{\includegraphics[width=2in]{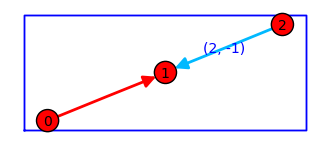}
\includegraphics[width=1.5in]{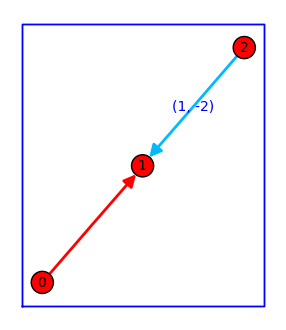}}

\sagepromt S3.quiver().digraph().edges()

\sageret{[(0, 1, (1, -1)), (2, 1, (2, -1))]}

\sagepromt S4.quiver().digraph().edges()

\sageret{[(0, 1, (1, -1)), (2, 1, (1, -2))]}
}

\noindent There is an analogous command that works for cluster algebras of \emph{finite type}:

2) The command {\tt variable\_class} will output the list of all cluster variables obtained as one mutates through all mutation-equivalent seeds.  

\sageex{
\sagepromt S3.variable\_class()

$$\hspace*{-45pt}\bigg[ x_0, x_1, x_2, \frac{x_{1} + 1}{x_{0}}, \frac{x_{0} x_{2}^{2} + 1}{x_{1}}, \frac{x_{1} + 1}{x_{2}}, \frac{x_{0} x_{2}^{2} + x_{1} + 1}{x_{0} x_{1}}, \frac{x_{0} x_{2}^{2} + x_{1} + 1}{x_{1} x_{2}},$$
$$\frac{x_{0} x_{2}^{2} + x_{1}^{2} + 2 x_{1} + 1}{x_{0} x_{1} x_{2}}, \frac{x_{0} x_{2}^{2} + x_{1}^{2} + 2 x_{1} + 1}{x_{1} x_{2}^{2}}, \frac{x_{1}^{3} + x_{0} x_{2}^{2} + 3 x_{1}^{2} + 3 x_{1} + 1}{x_{0} x_{1} x_{2}^{2}},$$
$$\hspace*{140pt}\frac{x_{0}^{2} x_{2}^{4} + 3 x_{0} x_{1} x_{2}^{2} + x_{1}^{3} + 2 x_{0} x_{2}^{2} + 3 x_{1}^{2} + 3 x_{1} + 1}{x_{0} x_{1}^{2} x_{2}^{2}}\bigg]$$

\sagepromt S4.variable\_class()

$$\hspace*{-40pt}\bigg[x_{0}, x_{1}, x_{2}, \frac{x_{1} + 1}{x_{0}}, \frac{x_{0} x_{2} + 1}{x_{1}}, \frac{x_{1}^{2} + 1}{x_{2}}, \frac{x_{0} x_{2} + x_{1} + 1}{x_{0} x_{1}}, \frac{x_{1}^{2} + x_{0} x_{2} + 1}{x_{1} x_{2}},$$
$$\hspace*{-45pt}\frac{x_{1}^{3} + x_{1}^{2} + x_{0} x_{2} + x_{1} + 1}{x_{0} x_{1} x_{2}}, \frac{x_{0}^{2} x_{2}^{2} + x_{1}^{2} + 2 x_{0} x_{2} + 1}{x_{1}^{2} x_{2}},$$
$$\frac{x_{0}^{2} x_{2}^{2} + x_{1}^{3} + x_{0} x_{1} x_{2} + x_{1}^{2} + 2 x_{0} x_{2} + x_{1} + 1}{x_{0} x_{1}^{2} x_{2}},$$
$$\hspace*{95pt}\frac{x_{1}^{4} + x_{0}^{2} x_{2}^{2} + 2 x_{1}^{3} + 2 x_{0} x_{1} x_{2} + 2 x_{1}^{2} + 2 x_{0} x_{2} + 2 x_{1} + 1}{x_{0}^{2} x_{1}^{2} x_{2}}\bigg]$$
}

\noindent In conclusion, even though the quivers of type $B_3$ and $C_3$ look quite similar and they have the same cluster complex, the Laurent polynomials are quite different.  For example, the bipartite seed for a cluster algebra of type $B_3$ leads to cluster variables whose numerators have degree $6$, while the numerators are only of degree at most $4$ in the case of $C_3$.  Similar phenomena happen for other dual cluster algebras, e.g. types $B_n$ versus $C_n$ for $n\geq 3$, or pairs of seeds: $(\mathbf{x},B)$ and $(\mathbf{x}, B^T)$.  Here and below, we adapt the term \lq\lq dual\rq\rq\ from the notion for Kac-Moody algebras.

Nuances like these make the non-skew-symmetric cases more difficult to analyze.  Nonetheless, using the classification (via folding of skew-symmetric quivers) appearing in \cite{FeSTuII}, it has been possible to include descriptions of mutation classes for those classes that correspond to a non-simply laced Dynkin diagram of finite or affine type, as well as the weighted quivers listed as exceptional cases in \cite{FeSTuII}.  For the classification of non-simply laced affine Dynkin diagrams, we use the tables of Kac \cite[pgs. 53-55]{Kac}.  However, the notation here is not explicit enough either as a number of cluster algebra mutation classes are again collapsed together.  We therefore follow notation of Dupont-P\'erotin \cite{DupontPerotin} instead.  The Dupont-P\'erotin notation specifies a quiver by indicating what the two ends look like, where the choices are that of a Dynkin diagram of type $B$, $C$ or $D$.  We say more about this notation in the next section. Since many users might be more familiar with the Kac-Moody notation, through careful coercing, if a user inputs a typical Kac-Moody type, it is recognized and translated into the appropriate notation that our software uses.

\sageex{
\sagepromt QuiverMutationType('C',2)

\sageret{['B',2]}

\sagepromt QuiverMutationType('B',4,1)

\sageret{['BD',4,1]}

\sagepromt QuiverMutationType('C',4,1)

\sageret{['BC',4,1]}

\sagepromt QuiverMutationType('A',2,2)

\sageret{['BC',1,1]}

\sagepromt QuiverMutationType('A',4,2)

\sageret{['BC', 2, 1]}

\sagepromt QuiverMutationType('A',5,2)

\sageret{['CD', 3, 1]}

\sagepromt QuiverMutationType('A',6,2)

\sageret{['BC', 3, 1]}

\sagepromt QuiverMutationType('A',7,2)

\sageret{['CD', 4, 1]}

\sagepromt QuiverMutationType('D',5,1)

\sageret{['D',5,1]}

\sagepromt QuiverMutationType('D',5,2)

\sageret{['CC',5,1]}

\sagepromt QuiverMutationType('D',4,3)

\sageret{['G',2,-1]}

\sagepromt QuiverMutationType('E',6,1)

\sageret{['E',6,1]}

\sagepromt QuiverMutationType('E',6,2)

\sageret{['F',4,-1]}

\sagepromt QuiverMutationType('F',4,1)

\sageret{['F',4,1]}
}

As for the finite types, our program has algorithms for identifying exchange matrices of affine types. In affine type $\tilde A_n$, we have a similar coercion issue in the case of simply-laced affine $\tilde{A}_{r,s}$ types where two parameters (rather than one parameter) is required to specify a mutation-equivalence type.  This example is special because it is the only finite or affine type with a Dynkin diagram which is not a tree.  Instead its Dynkin diagram is a cycle on $n$ vertices, and here quivers $Q_1$ and $Q_2$ are only mutation-equivalent if they have the same number of edges oriented clockwise and the same number of edges oriented counter-clockwise.  Actually, if all arrows are reversed, it is also the same type.  The mutation classes of types $\tilde{A}_{r,s}$ can be classified using theoretical results in \cite{Bastian}.

\sageex{
\sagepromt Qu = Quiver(['A',[2,3],1]); Qu

\sageret{Quiver on 5 vertices of type ['A', [2, 3], 1]}

\sagepromt Qu.show()

\sagepromt Quiver(['A',[4,1],1]).show()

\sagepromt Quiver(['A',[3,3],1]).show()

\sageret{\includegraphics[width=1.5in]{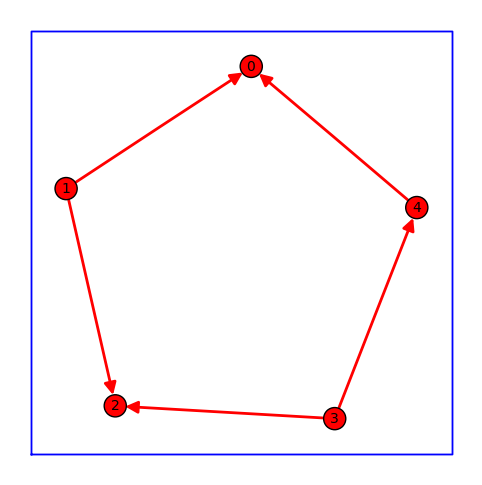}
\includegraphics[width=1.5in]{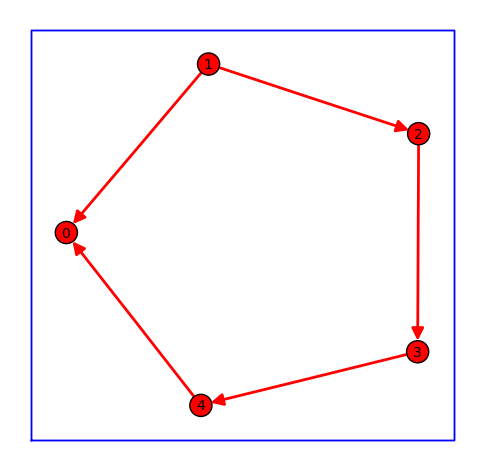}
\includegraphics[width=1.5in]{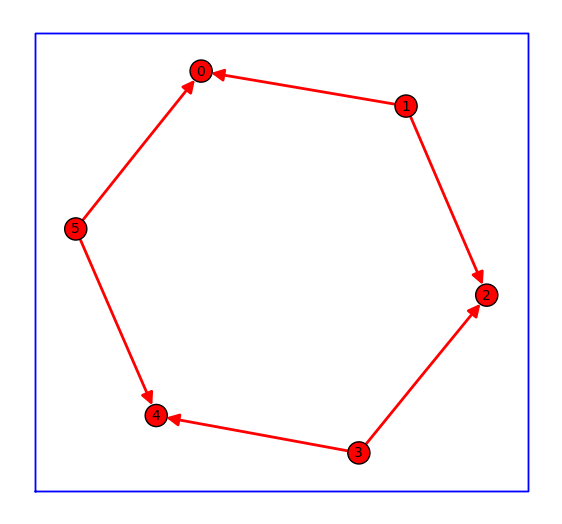}}
}

\noindent Notice also that the representative quiver for an affine $\tilde{A}_{r,s}$-type is made as bipartite as possible and that  mutation type {\tt ['A',[r,s],1]} is coerced into type {\tt ['A',[s,r],1]} when $s < r$.

The remaining affine types can be found in Section~\ref{QMT} and are classified using results in \cite{Hen2009} and \cite{Stump-pre}.

Beside the described coercions, we also include some basic coercions such as letting type $D_2$ coerce into type $A_1 \times A_1$, $D_3$ coerce into $A_3$, $C_2$ coerce into $B_2$, small rank two examples $\mathcal{A}(b,c)$ coerce into $A_2$, $B_2$, $G_2$, and $\tilde{A}_{1,1}$, and $\tilde{BC}_1$ for $(b,c) = (1,1), (1,2), (1,3), (2,2)$, and $(1,4)$, respectively.  Here, $\tilde{BC}_1$ simply means the type {\tt ['BC',1,1]} which is a denegenerate version of the {\tt ['BC',n,1]} family of Dynkin diagrams used above.  More technical details can be found in Section \ref{QMT}, including other families of types and more coercions.

\subsection{Class sizes of finite and affine quiver mutation types}\label{sec:classsizes}

In this section, we discuss the sizes of mutation classes of finite and affine types. Those results and conjectures are used to compute the size of mutation classes without explicitly computing the class. The \emph{class size} of a cluster seed or quiver is defined to be the number of exchange matrices or quivers which are mutation-equivalent to the given cluster seed or quiver, respectively. Here, we consider seeds and quivers up to isomorphism.
\begin{theorem}[Class sizes of finite types]
  The number of exchange matrices or quivers of finite
  \begin{itemize}
      \item type $A_n$ \cite{Tor2008} is given by
      $$\frac{1}{n+3}\left[ \frac{1}{n+1}\binom{2n}{n} + \binom{n+1}{(n+1)/2} + \binom{2n/3}{n/3}\right],$$
      where the second term is omitted if $(n+1)/2$ is not integral and the third term if $n/3$ is not integral.

      \item type $B_n$ or of type $C_n$ \cite{Stump-pre} is given by
      $$\frac{1}{n+1}\binom{2n}{n}.$$
    
      \item type $D_n$ \cite{BT2009} is for $n=4$ given by $6$, and for $n \geq 5$, it is given by
      $$\sum_{d|n} \frac{\phi(n/d)}{2n} \binom{2d}{d}.$$
      
      \item types $E_6,E_7,E_8,F_4,$ and $G_2$ are given by $67,416,1574,15,$ and $2$.
  \end{itemize}
\end{theorem}
\begin{theorem}[\cite{BPRS}]
  The number of exchange matrices or quivers of affine type $\tilde{A}_{r,s}$ is given by
      $$
        \begin{cases}
          \frac{1}{2} \sum\limits_{k|r,k|s}\frac{\phi(k)}{r+s}\binom{2r/k}{r/k}\binom{2s/k}{s/k}& \text{if } r \neq s,\\ \\
          \frac{1}{2} \left(\frac{1}{2} \binom{2r}{r} + \sum\limits_{k|r} \frac{\phi(k)}{4r}
          \binom{2r/k}{r/k}^2 \right) & \text{if } r=s.
        \end{cases}
      $$
      where $\phi(k)$ is Euler's totient function, i.e., the number of $1\leq d \leq k$ coprime to $k$.
\end{theorem}
\begin{conjecture}[\cite{Stump-pre}]
  The number of exchange matrices or quivers of affine
  \begin{itemize}
      \item type $\tilde{BB}_n$ or of type $\tilde{CC}_n$ is given by
      $$\binom{2n-1}{n-1} + \binom{n-1}{n/2-1}$$
      where the second term is omitted if $n$ is odd.\\
    
      \item type $\tilde D_n$ is for $n=4$ given by $9$, and for $n \geq 5$, it is given by
      $$2\binom{2n}{n} + \binom{n}{n/2},$$
      where the second term is omitted if $n$ is odd.\\

      \item type $\tilde{BC}_n$ is given by
      $$\binom{2n}{n}.$$
    
      \item type $\tilde{BD}_n$ or of type $\tilde{CD}_n$ is given by
      $$2\binom{2(n-1)}{n-1}.$$
    \end{itemize}
\end{conjecture}
\begin{theorem}
  The number of exchange matrices or quivers of
  \begin{itemize}
   \item affine types $\tilde E_6,\tilde E_7,\tilde E_8,\tilde F_4,$ and $\tilde G_2$ are given by $132,1080,7560,60,$ and $6$.
   \item elliptic types $\tilde E_6^{(1)},\tilde E_7^{(1)},$ and $\tilde E_8^{(1)}$ are given by $49,506,$ and $5739$.
   \item the other exceptional mutation-finite types $V_4,W_4,X_6,X_7,Y_6,$ and $Z_6$ are given by $7,2,5,2,90,$ and $35$.
  \end{itemize}
\end{theorem}

\subsection{Algorithms for computing mutation classes}
\label{sec:varclass}

The four commands
\sageret{{\tt mutation\_class}, \quad {\tt b\_matrix\_class}, \quad {\tt cluster\_class}, \quad {\tt variable\_class}}
each utilize the auxiliary command obtained by adding {\tt \_iter}, which constructs an iterator that will run through all the objects in the corresponding mutation class.  For quivers, there is only the method {\tt mutation\_class}. The first three methods are directly derived from {\tt mutation\_class\_iter}, we therefore begin by describing how this method works.

Note first that {\tt mutation\_class\_iter} is, as the name already indicates, an \emph{iterator}. This means that the next element is only computed if the iterator is asked to do so. Here is an example. One might be interested if there exists a seed or quiver in a given infinite mutation class having a certain property. Of course, we cannot test all elements, but we can construct the iterator and then let the computer run through the elements, constructing one after the other, and checking this property. If the program finds an element having the property, one could halt the process and return the element, together with all mutations applied to the initial element. If the computer keeps running, you might (or might not) get convinced that such an element does not exist.

The command {\tt mutation\_class\_iter} has five (resp. six) optional arguments if it is acting on a cluster seed (resp. quiver).  The additional optional argument for quivers is {\tt data\_type} which is initially set to {\tt `quiver'} but can also be allowed to be {\tt matrix}, {\tt digraph}, {\tt dig6}, or {\tt path}.  This argument does not appear in the cluster seed since the data type is assumed to be a cluster seed here.

The second optional argument is {\tt depth}, which is set to be `infinity` by default and instructs how large a ``ball'' in the mutation-equivalence-class around the initial input is supposed to be constructed. If the cluster algebra is of finite type (resp. finite mutation type) however then a depth of infinity will eventually construct the entire mutation class, when the original input is a cluster seed (resp. quiver).  

Another optional argument is {\tt show\_depth}, which allows the user to print extra information of the actual depth, the number of constructed seeds or quivers, and the elapsed time. It is set to be {\tt False} by default.  The argument {\tt up\_to\_equivalence} works differently depending on whether the input is a cluster seed or a quiver.  In the default case \True, cluster seeds are considered up to simultaneous row and column permutations and quivers are considered unlabeled; see Remark~\ref{rem:equivalence}. Otherwise, equivalence of seeds and quivers are not considered.  

\sageex{
\sagepromt S = ClusterSeed(['A',2]);  

\sagepromt S.cluster\_class() 

$$\left[\left[x_{0}, x_{1}\right], \left[x_{0}, \frac{x_{0} + 1}{x_{1}}\right], \left[\frac{x_{1} + 1}{x_{0}}, x_{1}\right], \left[\frac{x_{0} + x_{1} + 1}{x_{0} x_{1}}, \frac{x_{0} + 1}{x_{1}}\right], \left[\frac{x_{1} + 1}{x_{0}}, \frac{x_{0} + x_{1} + 1}{x_{0} x_{1}}\right]\right]$$

\sagepromt S.cluster\_class(up\_to\_equivalence=False) 

$$\bigg[\left[x_{0}, x_{1}\right], 
\left[x_{0}, \frac{x_{0} + 1}{x_{1}}\right], 
\left[\frac{x_{1} + 1}{x_{0}}, x_{1}\right], 
\left[\frac{x_{1} + 1}{x_{0}}, \frac{x_{0} + x_{1} + 1}{x_{0} x_{1}}\right], 
\left[\frac{x_{0} + x_{1} + 1}{x_{0} x_{1}}, \frac{x_{0} + 1}{x_{1}}\right],$$
$$\left[\frac{x_{0} + x_{1} + 1}{x_{0} x_{1}}, \frac{x_{1} + 1}{x_{0}}\right], \left[\frac{x_{0} + 1}{x_{1}}, \frac{x_{0} + x_{1} + 1}{x_{0} x_{1}}\right], \left[x_{1}, \frac{x_{1} + 1}{x_{0}}\right], \left[\frac{x_{0} + 1}{x_{1}}, x_{0}\right], \left[x_{1}, x_{0}\right]
\bigg]$$
}

\noindent The argument {\tt sink\_source} is set to be {\tt False} by default, but if set to {\tt True}, then only mutations at sinks and sources are performed.  This option is helpful for working with bipartite seeds or studying the BGP reflection functors on quiver representations.

Finally, the last argument {\tt return\_paths}, again {\tt False} by default, will keep track of the shortest mutation sequence that can be used to produce a given seed  (or quiver) from the initial one.  This data can be accessed by other commands and then utilized for future work.  Note that such a sequence is not unique so accessing this shortest sequence during different computational sessions might not give the same result but for most purposes a single example of the mutation sequence between two seeds is sufficient data.

With this iterator, one can then call {\tt mutation\_class} which will output the associated list of seeds or quivers in the mutation class.  However, since this output cannot be infinite, the argument {\tt depth} cannot be {\tt infinity} unless the input is of finite (resp. finite mutation type).  The data associated to the optional arguments is also returned at this time.  The commands {\tt b\_matrix\_class} and {\tt cluster\_class}, which each can only be performed on a cluster seed, work analogously. The algorithm for {\tt variable\_class}, which again only works on a cluster seed, requires a little more explanation.

The procedure for {\tt variable\_class\_iter} starts by running through an iterator for the mutation class and by yielding all found cluster variables.  However, since the set of cluster variables is dwarfed by the number of clusters, this search-based algorithm is quite slow.  

On the other hand, if we are in the lucky situation that the initial cluster is bipartite, then we can use \cite[Theorem~8.8]{ClustIV} to efficiently compute the variable class.

\begin{theorem} [Theorem 8.8 of \cite{ClustIV}]
Suppose that an exchange matrix $B$ is bipartite, and its Cartan counterpart $A=A(B)$ is indecomposable.

1) If $A$ is of finite type, then the corresponding bipartite belt (see Definition \ref{BipBelt}) has the following periodicity property: the labeled seeds $\Sigma_m$ and $\Sigma_{m+2(h+2)}$ are equal to each other for all $m \in \mathbb{Z}$.  Here, $h$ is the Coxeter number of the corresponding Cartan matrix $A$.

2) If $A$ is of infinite type, then all of the elements $x_{i;2m}$, denoting the $n$ cluster variables of $\Sigma_{2m} = (\{x_{1;2m},x_{2;2m},\dots, x_{n;2m}\},B)$ as $m$ ranges over the integers are distinct Laurent polynomials in the initial data.
\end{theorem}

Note that in this theorem, the \emph{Cartan counterpart} of $B$ (see Section~\ref{sec:associahedra}) is the (generalized) \emph{Cartan matrix} $A=A(B) = (a_{ij})$ defined by $$a_{ij} = \begin{cases} 2 &\mathrm{~if~} i = j \\ -|b_{ij}| &\mathrm{~if~} i \not = j \end{cases}.$$

\begin{definition} \label{BipBelt}
We use $\Sigma_0 = (\mathbf{x_0}, B)$ to denote an initial bipartite seed and let $\mu_+$ (resp. $\mu_-$) denote the concatenation of all mutations at sources (sinks) of the quiver $Q(B)$\footnote{Since sources and sinks are not adjacent, the factors of $\mu_+$ (resp. $\mu_-$) commute with one another, hence why $\mu_+$ and $\mu_-$ are well-defined.}.  Observe that $\mu_+(B) = \mu_-(B) = -B$. 

Define the associated \emph{bipartite belt} to be the seeds $\Sigma_m = (\mathbf{x_m}, (-1)^m B)$ for $m \in \mathbb{Z}$, defined recursively by 
$$\Sigma_r = \begin{cases} \mu_+ (\Sigma_{r-1}) &\mathrm{~if~}r\mathrm{~is~odd} \\  
\mu_- (\Sigma_{r-1}) &\mathrm{~if~}r\mathrm{~is~even}.
\end{cases}$$
\end{definition}

As a consequence, given an initial bipartite seed $(\mathbf{x},B)$, it is sufficient to mutate all vertices labeling sinks in $Q(B)$ followed by mutating all vertices labeling sources in $Q(B)$, and iterate.  We will get no repeats in this list and thus the most efficient way to obtain all cluster variables in the case of a finite type cluster algebra\footnote{\label{bibfoot} If the cluster algebra is of infinite type, one can also mutate along the bipartite belt to efficiently generate a large list of cluster variables but \emph{not all cluster variables are reachable} in this way.}.

Our algorithm thus first checks if the initial seed is bipartite for this reason.  If not, it proceeds as above trying to mutate in all directions.  

It is a difficult computational problem to find a mutation sequence, if one exists, from an initial non-bipartite seed to a bipartite one, so it is not computationally feasible to use the shortcut if we do not have a bipartite seed at hand.  However, since our proceeding is doing a search through all seeds mutation-equivalent to the initial one anyway as its default behavior if we get lucky and find a bipartite seed, the program can record this path and take advantage of this find.  

In the case that the search algorithm finds a bipartite seed, the algorithm then does the following procedure instead:

1) Starts over at the initial seed.

2) Mutates along the recorded path to get to the bipartite seed $\Sigma_0$.

3) Mutate along the bipartite belt the appropriate distance from there in both directions (i.e. applying $\mu_+$ first or $\mu_-$ first).

In step (3) the appropriate distance is either the period $2(h+2)$ in the case of a cluster algebra of finite type or the {\tt depth} chosen beforehand by the user.  Note well that the meaning of {\tt depth} is actually different here, as the algorithm will no longer spread out in all directions.  Instead, the argument {\tt depth} now instructs the computer how many iterations of the bipartite belt to use.  The program will actually output the cluster variables found on the way to the bipartite seed $\Sigma_0$, as well as all cluster variables in the seeds $\{\Sigma_m: m\in \mathbb{Z}, |m| \leq \mathtt{depth}\}$.

Since in the case of infinite type, not all cluster variables can be reached by using the bipartite belt, for example even cluster variables lying in clusters two mutations away from the bipartite seed might not be reachable (see the bipartite $\tilde{A}_{2,2}$ example below), the optional argument {\tt ignore\_bipartite\_belt=False} is included.  If set to be {\tt True}, the original (albeit slower) algorithm of mutating in all directions out to a certain depth is utilized even if a bipartite seed is found.

\sageex{

\sagepromt S = ClusterSeed(['A',[2,2],1]); S.b\_matrix(); S.is\_bipartite()

\sageret{ $\left(\begin{array}{rrrr}
0 & -1 & 0 & -1 \\
1 & 0 & 1 & 0 \\
0 & -1 & 0 & -1 \\
1 & 0 & 1 & 0
\end{array}\right)$
\qquad
{\tt True}
}

\sagepromt S.variable\_class(depth=1)

\sageret{\hspace{-170pt}Found a bipartite seed - \\ \hspace{50pt}constructing the variable class into its bipartite belt.}

$$\hspace{-50pt}\bigg[x_{0}, x_{1}, x_{2}, x_{3}, \frac{x_{1} x_{3} + 1}{x_{0}}, \frac{x_{0} x_{2} + 1}{x_{1}}, \frac{x_{1} x_{3} + 1}{x_{2}}, \frac{x_{0} x_{2} + 1}{x_{3}},$$
$$\frac{x_{1}^{2} x_{3}^{2} + x_{0} x_{2} + 2 x_{1} x_{3} + 1}{x_{0} x_{1} x_{2}}, \frac{x_{0}^{2} x_{2}^{2} + 2 x_{0} x_{2} + x_{1} x_{3} + 1}{x_{0} x_{1} x_{3}},$$ 
$$\hspace{50pt}\frac{x_{1}^{2} x_{3}^{2} + x_{0} x_{2} + 2 x_{1} x_{3} + 1}{x_{0} x_{2} x_{3}}, \frac{x_{0}^{2} x_{2}^{2} + 2 x_{0} x_{2} + x_{1} x_{3} + 1}{x_{1} x_{2} x_{3}}\bigg]$$

}  

\noindent If we look at the output from {\tt S.variable\_class(depth=2)} or higher depth, we will see that the denominators grow larger and larger but no denominator of $x_0x_1$ appears.  Compare this output with the examples below.

\sageex{

\sagepromt S.mutate([0,1]); S.cluster()

$$\left[\frac{x_{1} x_{3} + 1}{x_{0}}, \frac{x_{0} x_{2} + x_{1} x_{3} + 1}{x_{0} x_{1}}, x_{2}, x_{3}\right]
$$

\sagepromt S.variable\_class(depth=2, ignore\_bipartite\_belt=True)

$$\hspace{-120pt}\bigg[x_{0}, x_{1}, x_{2}, x_{3}, \frac{x_{1} x_{3} + 1}{x_{0}}, \frac{x_{0} x_{2} + 1}{x_{1}}, \frac{x_{1} x_{3} + 1}{x_{2}},$$
$$\hspace{-50pt}\frac{x_{0} x_{2} + x_{1} x_{3} + 1}{x_{0} x_{1}}, \frac{x_{0} x_{2} + x_{1} x_{3} + 1}{x_{0} x_{3}}, \frac{x_{0} x_{2} + x_{1} x_{3} + 1}{x_{1} x_{2}},$$
$$\frac{x_{1}^{2} x_{3}^{2} + x_{0} x_{2} + 2 x_{1} x_{3} + 1}{x_{0} x_{1} x_{2}}, \frac{x_{0}^{2} x_{2}^{2} + 2 x_{0} x_{2} + x_{1} x_{3} + 1}{x_{0} x_{1} x_{3}},$$
$$\hspace{60pt}\frac{x_{1}^{3} x_{3}^{3} + x_{0}^{2} x_{2}^{2} + 2 x_{0} x_{1} x_{2} x_{3} + 3 x_{1}^{2} x_{3}^{2} + 2 x_{0} x_{2} + 3 x_{1} x_{3} + 1}{x_{0}^{2} x_{1} x_{2} x_{3}}\bigg]$$

}

\section{Associahedra and the cluster complex} \label{sec:associahedra}

Before looking at associahedra, the cluster complex and their implementations, we need to start with some basic background on root systems for (generalized) Cartan matrices. For further details, we refer to \cite{Hum1972,Kac}.

\begin{definition}[Generalized Cartan matrix]
  An $n \times n$-matrix $A = (a_{ij})$ with integer entries is called a \emph{generalized Cartan matrix} if
  \begin{itemize}
    \item $a_{ii} = 2$,
    \item $a_{ij} < 0$ for $i \neq j$,
    \item $A$ is symmetrizable, i.e., there exists a diagonal matrix $D$ with positive entries such that $DA$ is symmetric.
  \end{itemize}
  A generalized Cartan matrix is called \emph{of finite type} if $DA$ is positive definite, and \emph{of affine type} if $DA$ is positive semi-definite.
\end{definition}
Recalling the definition of $B$-matrices, we see that we can associate a generalized Cartan matrix to every $B$-matrix (see \cite[(1.6)]{ClustII}). The terms \emph{finite} and \emph{affine} come from their connections to finite and affine \emph{Lie algebras}. Indecomposable generalized Cartan matrices of finite type (resp. of affine type) classify Lie algebras of finite type (resp. of affine type).

A \emph{realization} of a Cartan matrix $A$ (of finite type) is a (rational, real, or complex) vector space $V$ with distinguished basis $\Delta = \{ \alpha_i : 0 \leq i < n \}$, and with dual space $V^*$ with distinguished basis $\Delta^\vee = \{ \alpha^\vee_i : 0 \leq i < n \}$, together with the pairing $\langle \alpha^\vee_i, \alpha_j \rangle = a_{ij}$. For $\beta \in V$ (resp. $\beta^\vee \in V^\vee$), we write $[ \beta, \alpha_i ]$ (resp. $[ \beta^\vee, \alpha_i^\vee ]$) for the coefficient of $\alpha_i$ in $\beta$ (resp. $\alpha_i^\vee$ in $\beta^\vee$).

Define a \emph{reflection} on $V$ by
$$s_i(\alpha_j) = \alpha_j - a_{ij} \alpha_i,$$
and define moreover, the \emph{Weyl group} by $W = \langle s_i : 0 \leq i < n \rangle \leq \operatorname{Aut}(V)$ and the \emph{root system} by
$$\Phi = \big\{ \omega(\alpha) : \omega \in W, \alpha \in \Delta \big\}.$$

It can be shown that $\Phi$ can be written as $\Phi^+ \cup \Phi^-$ where
$$\Phi^+ = \{ \beta \in \Phi : [ \beta, \alpha ] > 0 \text{ for all } \alpha \in \Delta \},$$
and $\Phi^- = \{-\beta: \beta \in \Phi^+$. The elements in $\Phi$ are called \emph{roots}, the elements in $\Phi^+$ are called \emph{positive roots}, and the elements in $\Delta$ are called \emph{simple roots}.

\begin{theorem} [Theorem 1.9 of \cite{ClustII}]
  Let $\mathcal{A}$ be a Cluster algebra of finite type and let $\Phi_{\geq -1} = \Phi^+ \cup -\Delta$ be the set of \emph{almost positive roots} of the root system of the associated Cartan type given by the positive roots together with the simple negative roots. There exists a unique bijection between almost positive roots and the cluster variables for $\mathcal{A}$ for which the simple negative root $-\alpha$ is mapped to $x_\alpha$ and, for positive roots,
  $$\sum_{\alpha \in \Delta}n_\alpha \alpha \mapsto \frac{P_\alpha}{\prod x_\alpha^{n_\alpha}},$$
  with $P_\alpha$ having nonzero constant term. Here, $x_{\alpha_i}$ stands for $x_i$ for an appropriate ordering $\Delta = \{ \alpha_0, \ldots, \alpha_{n-1} \}$.
\end{theorem}

This connection in the finite types can be used in the cluster algebra package as follows:
  \sageex{

    \sagepromt for f in ClusterSeed(['A',2]).variable\_class():

    \sagedots print f, f.almost\_positive\_root()

    \begin{align*}
      x_0                       &\qquad -\alpha_1\\
      x_1                       &\qquad -\alpha_2\\
      (x_1 + 1)/x_0             &\qquad \alpha_1\\
      (x_0 + 1)/x_1             &\qquad \alpha_2\\
      (x_0 + x_1 + 1)/(x_0 x_1) &\qquad \alpha_1 + \alpha_2
    \end{align*}
    
    \sagepromt f
    
    \sageret{$(x_0 + x_1 + 1)/(x_0 x_1)$}
    
    \sagepromt root = f.almost\_positive\_root(); root

    \sageret{$\alpha_1 + \alpha_2$}
    
    \sagepromt root.parent()
    
    \sageret{Root lattice of the Root system of type ['A', 2]}
  }  

\subsection{Generalized associahedra} In this section, we will define generalized associahedra and describe how they can be realized as polytopal complexes in finite types. We will see then how these polytopal complexes are implemented in \sage. Generalized associahedra beyond finite type are not yet feasible as the needed tools to deal with infinite types are not yet developed. We start with the definition of generalized associahedra (not necessarily of finite type).

\begin{definition}[Generalized associahedron]
 The \emph{generalized associahedron} associated to a cluster algebra $\mathcal{A}$ can be defined as the exchange graph of the mutation class of all cluster seeds for $\mathcal{A}$. This is the unoriented graph with vertices given by the set of all cluster seeds, and with an edge joining two clusters if they can be obtained from each other by a mutation.
\end{definition}
Generalized associahedra reduce in classical types to known constructions, see e.g.~\cite[Section~12]{ClustII}. By \cite[Theorem~1.12]{ClustII}, a cluster seed of finite type is uniquely determined by its cluster, and two seeds are obtained from each other by a mutation if and only their clusters differ by exactly one cluster variable, see Theorem~\ref{th:exchangerelation}. In finite types, there exist realizations as polytopal complexes, see \cite{Polytope}. Let $S_+,S_-$ be the bipartition of the simple reflections $S = \{ s_\alpha : \alpha \in \Delta \}$ corresponding to the simple roots in $\Delta$. This means that $S_+$ and $S_-$ are chosen in such a way that the reflections in each pairwise commute. Observe that the fact that all quivers of finite type are bipartite ensures that such bipartitions always exist. Define two piecewise linear operators $\tau_+$ and $\tau_-$ on $V$ by
$$
  \tau_\epsilon(\beta) = 
    \begin{cases}
      \beta & \text{if } \beta = -\alpha \text{ for } s_\alpha \in S_{-\epsilon} \\
      \prod_{ s \in S_\epsilon } s (\beta) & \text{otherwise,}
    \end{cases}
$$
and let
$$\rho^\vee = \frac{1}{2} \sum_{\beta \in \Phi^+} \beta^\vee \in V^*.$$
In \cite[Theorem~1.1]{Polytope}, it is shown that every $\langle \tau_+,\tau_- \rangle$-orbit in $\Phi_{\geq -1}$ intersects $-\Delta$. Moreover, $\alpha_i,\alpha_j \in -\Delta$ lie in the same orbit if and only if $\alpha_i = -\omega_0(\alpha_j)$ where $\omega_0$ is the (unique) \emph{longest element} in $W$. Thus, the coefficients $[ \rho^\vee, \alpha_i^\vee ]$ and $[ \rho^\vee, \alpha_j^\vee ]$ coincide; for $\beta \in \Phi_{\geq -1}$, set $c_\beta$ to be this coefficient. After identifying $\varphi$ with the $n$-tuple $( \langle \varphi, \alpha_i \rangle )_{0 \leq i < n}$, define the half-space
$$H^+(\beta) := \{ \varphi \in \mathbb{R}^n : \langle \varphi, \beta \rangle \leq c_\beta \}$$
to obtain the polytopal realization of the generalized associahedron by
$$ \operatorname{Ass}(\Phi) = \bigcap_{\beta \in \Phi^+} H^+(\beta) \subseteq \mathbb{R}^n.$$
The operators $\tau_+$ and $\tau_-$ are implemented in \sage\ as operators for the root space.

\sageex{
  \sagepromt S = RootSystem(['A',2]).root\_space()

  \sagepromt tau\_plus, tau\_minus = S.tau\_plus\_minus()
  
  \sagepromt for beta in S.almost\_positive\_roots():
  
  \sagedots  \hspace{10pt} print beta, tau\_plus(beta), tau\_minus(beta)
  
  \sagedots  \hspace{10pt} print

  \sageret{
    \begin{align*}
      -\alpha_1,\quad \alpha_1, \quad -\alpha_1 \\
      \alpha_1, \quad-\alpha_1, \quad \alpha_1 + \alpha_2 \\
      \alpha_1 + \alpha_2, \quad \alpha_2, \quad \alpha_1 \\
      -\alpha_2, \quad -\alpha_2, \quad \alpha_2 \\
      \alpha_2, \quad \alpha_1 + \alpha_2, \quad -\alpha_2
    \end{align*}
}}
  \sageex{

  \sagepromt AssoA2 = Associahedron(['A',2]); AssoA2
  
  \sageret{The generalized associahedron of type ['A', 2] \\ \hspace{85pt} having 2 dimensions and 5 vertices}

  \sagepromt AssoB2 = Associahedron(['B',2]); AssoB2

  \sageret{The generalized associahedron of type ['B', 2] \\ \hspace{85pt} having 2 dimensions and 6 vertices}
  
  \sagepromt AssoC2 = Associahedron(['C',2]); AssoC2
  
  \sageret{The generalized associahedron of type ['C', 2] \\ \hspace{85pt} having 2 dimensions and 6 vertices}
  
  \sagepromt AssoG2 = Associahedron(['G',2]); AssoG2

  \sageret{The generalized associahedron of type ['G', 2] \\ \hspace{85pt} having 2 dimensions and 8 vertices}

  \sagepromt AssoA2.show(); AssoB2.show(); AssoC2.show(); AssoG2.show()
  
  \sageret{
    \includegraphics[width=2in]{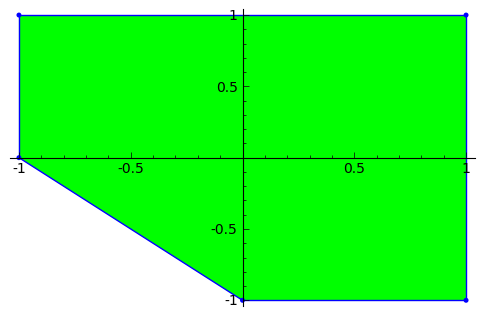} \quad
    \includegraphics[width=2in]{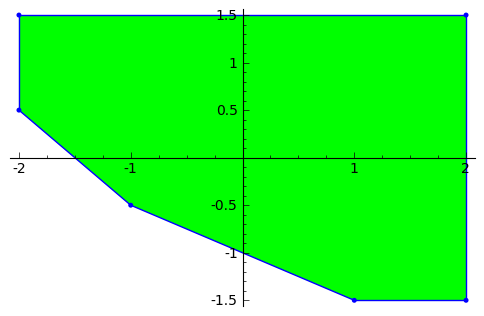}\\[10pt]
    \includegraphics[width=2in]{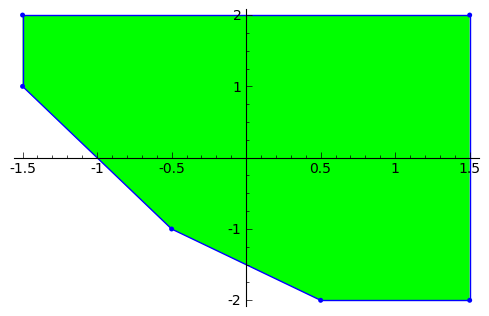} \quad
    \includegraphics[width=2in]{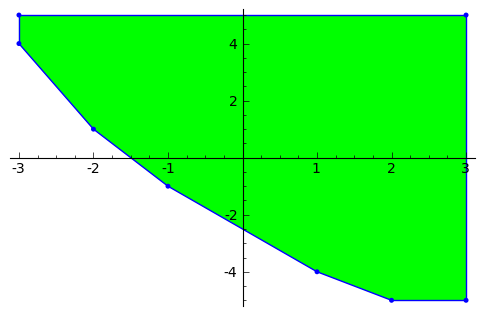}
  }
  
  \sagepromt AssoA3 = Associahedron(['A',3]); AssoA3
  
  \sageret{The generalized associahedron of type ['A', 3] \\ \hspace{85pt} having 3 dimensions and 14 vertices}
  
  \sagepromt AssoB3 = Associahedron(['B',3]); AssoB3
  
  \sageret{The generalized associahedron of type ['B', 3] \\ \hspace{85pt} having 3 dimensions and 20 vertices}
  
  \sagepromt AssoA3.show(); AssoB3.show()
  \sageret{
    \includegraphics[width=2.4in]{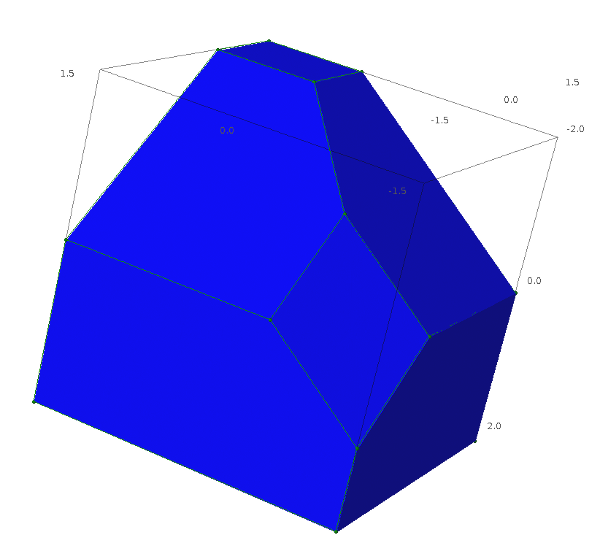}
    \includegraphics[width=2.4in]{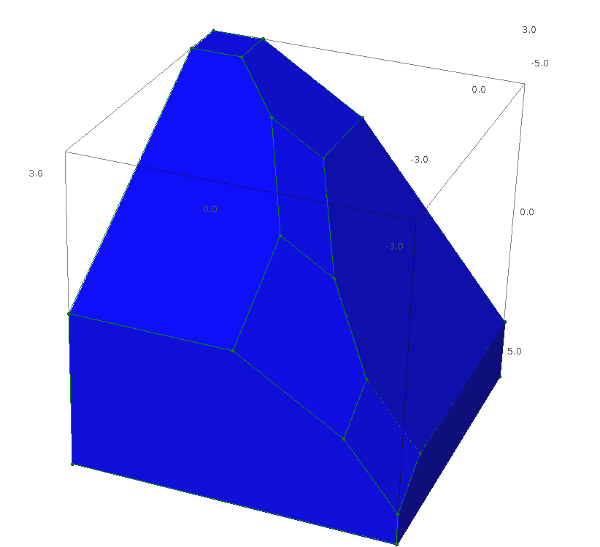}
  }
}
The associahedron of type $A_3$ has $14$ vertices ($13$ of which are visible, the $14$th is the origin, which corresponds to the cluster $\{ -\alpha_1,-\alpha_2,-\alpha_3\}$). As well, the $9$ facets corresponds to the almost positive roots, where the hyperplane $x_i = c_{-\alpha_i}$ correspond to the simple negative root $-\alpha_i$. Every vertex corresponds to exactly $3$ hyperplanes, and in type $B_3$, we have $20$ vertices and $12$ facets, as desired.

\subsection{The cluster complex} As with associahedra, we will define the cluster complex in general and then discuss the implementation for finite types.

\begin{definition}[Cluster complex]
 The \emph{cluster complex} associated to a cluster algebra $\mathcal{A}$ can be defined to be the simplicial complex with vertices being the cluster variables for $\mathcal{A}$ and with facets being the clusters.
\end{definition}

As we have seen, cluster variables in finite types are in bijection with almost positive roots. We use this description in the implementation of the cluster complex.

  \sageex{

  \sagepromt ClusterComplex(['A',2])

  \sageret{Simplicial complex with 5 vertices and 5 facets}

  \sagepromt ClusterComplex(['A',3])

  \sageret{Simplicial complex with 9 vertices and 14 facets}

  \sagepromt Delta = ClusterComplex(['B',3]); Delta

  \sageret{Simplicial complex with 12 vertices and 20 facets}
  }

In the following example, we see how we can use other \sage\ packages to further study objects we work with. As the cluster complex is a simplicial complex, there now exists various possible methods. For example, we can compute its \emph{homology},

  \sageex{

  \sagepromt Delta.homology()

  \sageret{$\{0: 0, 1: 0, 2: \mathbb{Z}\}$}

  }

  This is as expected, as this simplicial complex is the boundary complex of a triangulated polytope, and thus shellable and Cohen-Macaulay.

  \section{Methods and attributes}

\label{glossary}

In this section, we describe the different classes defined in this package, and list their attributes and methods. For the \lq\lq key\rq\rq\ methods, we also give descriptions of the algorithms.

In general, attribute names start with an underscore to emphasize that they should not be used directly but only through appropriate methods. As an example, a cluster seed has an attribute {\tt \_M} in which its exchange matrix is stored and a method {\tt b\_matrix} which is used to get the exchange matrix. The difference is that the method returns a copy of its exchange matrix, so it is safe to work with this matrix and to modify it without accidentally modifying the seed itself.

\sageex{
\sagepromt S = ClusterSeed(['A',3]);

\sagepromt M1 = S.\_M; M2 = S.b\_matrix();

\sagepromt M1 == M2

\sageret{True}

\sagepromt M1 is M2

\sageret{False}
}

\subsection{Skew-symmetrizable matrices}\label{sec:skew-sym}

We briefly want to describe the algorithm used to determine whether a square matrix $B$ is skew-symmetrizable, which also determines the associated diagonal matrix $D$ in the affirmative case. It was written by F.~Block, F.~Saliola, and C.~Stump during the \sage\ days 20.5 at the Fields Institute, Toronto, Canada, in May 2010.
\begin{algorithm}
 Let $B = (b_{ij})_{1 \leq i,j \leq n}$ be the input square matrix of dimension $n$, and let $D = (d_i)_{1 \leq i \leq n}$ be the diagonal matrix with positive coefficients we want to construct. We use the equivalent description of skew-symmetrizablility given by the property
 $$d_i b_{ij} = -d_j b_{ji} \text{ for all } i,j.$$
 \begin{enumerate}
  \item Check if $b_{ii} = 0$ for all $i$. If this is not the case, return \False,
  \item let $k$ be the smallest integer such that $d_k$ is not yet determined,
  \item set $d_k = 1$,
  \item for $i \in \{1,\ldots,n\}$ such that $b_{ik} \neq 0$ and $d_i$ is not yet determined, do
    \begin{enumerate}
      \item set $d_i = -d_k b_{ki} / b_{ik}$,
      \item if $d_i \leq 0$ return \False.
      \item if $\operatorname{any}\big( d_i b_{ij} \neq -d_j b_{ji}\big)$ for $j$ such that $d_j$ is already determined, return \False.
    \end{enumerate}
  \item repeat step \emph{(4)} with $k$ given by all integers for which $d_i$ was set since we passed step \emph{(3)} the last time,
  \item if $D$ is not yet completely determined, {\tt goto} step \emph{(2)},
  \item return $D$.
 \end{enumerate}

\end{algorithm}

\subsection{QuiverMutationType} \label{QMT} For coding reasons, we distinguish between the classes {\tt QuiverMutationType\_Irreducible} and {\tt QuiverMutationType\_Reducible}, but we refer here to both as {\tt QuiverMutationType}. Objects of those types are unique, i.e., there exists only one object of a given quiver mutation type.

\sageex{
\sagepromt mut\_type1 = QuiverMutationType('A',3)

\sagepromt mut\_type2 = QuiverMutationType('A',3)

\sagepromt mut\_type1 is mut\_type2

\sageret{True}
}

\noindent All the data for quiver mutation types is hard-coded. In particular, this concerns the graphs and digraphs, and the class size.

To construct a quiver mutation type, the function {\tt QuiverMutationType} is called. An irreducible quiver mutation type takes $3$ parameters, the {\tt letter}, the {\tt rank} or {\tt bi\_rank}, and the {\tt twist}, see the description below. Those calls are best explained in examples. Observe that the call arguments can be also wrapped into a list or tuple. We suppress the output whenever the output coincide with the input.
\begin{itemize}
  \item finite types

    \sageex{
            \sagepromt QuiverMutationType('A',1);
            
            \sagepromt QuiverMutationType('A',5);
            
            \sagepromt QuiverMutationType('B',2);
            
            \sagepromt QuiverMutationType('B',5);
            
            \sagepromt QuiverMutationType('C',2)
            
            \sageret{['B', 2]}
                        \sagepromt QuiverMutationType('C',5);
            
            \sagepromt QuiverMutationType('D',2)
            
            \sageret{[ ['A', 1], ['A', 1] ]}
            
            \sagepromt QuiverMutationType('D',3)
            
            \sageret{['A', 3]}
            
            \sagepromt QuiverMutationType('D',4);
            
            \sagepromt QuiverMutationType('E',6);
            
            \sagepromt QuiverMutationType('E',7);
            
            \sagepromt QuiverMutationType('E',8);
            
            \sagepromt QuiverMutationType('F',4);
            
            \sagepromt QuiverMutationType('G',2);
    }

  \item affine types

    \sageex{
            \sagepromt QuiverMutationType('A',(1,1),1);

            \sagepromt QuiverMutationType('A',(2,4),1);

            \sagepromt QuiverMutationType('BB',1,1)

            \sageret{['A', [1, 1], 1]}

            \sagepromt QuiverMutationType('BB',2,1);

            \sagepromt QuiverMutationType('BB',4,1);

            \sagepromt QuiverMutationType('CC',1,1)

            \sageret{['A', [1, 1], 1]}
            
            \sagepromt QuiverMutationType('CC',2,1);
            
            \sagepromt QuiverMutationType('CC',4,1);
            
            \sagepromt QuiverMutationType('BC',1,1);
            
            \sagepromt QuiverMutationType('BC',5,1);
            
            \sagepromt QuiverMutationType('BD',3,1);
            
            \sagepromt QuiverMutationType('BD',5,1);
            
            \sagepromt QuiverMutationType('CD',3,1);
            
            \sagepromt QuiverMutationType('CD',5,1);
            
            \sagepromt QuiverMutationType('D',4,1);
            
            \sagepromt QuiverMutationType('D',6,1);
            
            \sagepromt QuiverMutationType('E',6,1);
            
            \sagepromt QuiverMutationType('E',7,1);
            
            \sagepromt QuiverMutationType('E',8,1);
            
            \sagepromt QuiverMutationType('F',4,1);
            
            \sagepromt QuiverMutationType('F',4,-1);
            
            \sagepromt QuiverMutationType('G',2,1);
            
            \sagepromt QuiverMutationType('G',2,-1);
    }
    \item hyperbolic types

          \sageex{
            \sagepromt QuiverMutationType('E',6,[1,1]);

            \sagepromt QuiverMutationType('E',7,[1,1]);
            
            \sagepromt QuiverMutationType('E',8,[1,1]);
    }

    \item mutation-finite types
      \begin{itemize}
        \item rank $2$

          \sageex{
            \sagepromt QuiverMutationType('R2',(1,1),2)

            \sageret{['A', 2]}
            
            \sagepromt QuiverMutationType('R2',(1,2),2)
            
            \sageret{['B', 2]}
            
            \sagepromt QuiverMutationType('R2',(1,3),2)
            
            \sageret{['G', 2]}
            
            \sagepromt QuiverMutationType('R2',(1,4),2)
            
            \sageret{['BC', 1, 1]}

            \sagepromt QuiverMutationType('R2',(1,5),2);

            \sagepromt QuiverMutationType('R2',(2,2),2)

            \sageret{['A', [1, 1], 1]}

            \sagepromt QuiverMutationType('R2',(3,5),2);
          }
        \item exceptional types

          \sageex{
            \sagepromt QuiverMutationType('V',4,2);
            
            \sagepromt QuiverMutationType('W',4,2);
            
            \sagepromt QuiverMutationType('W',4,-2);
            
            \sagepromt QuiverMutationType('X',6,2);
            
            \sagepromt QuiverMutationType('X',7,2):
            
            \sagepromt QuiverMutationType('Y',6,2);
            
            \sagepromt QuiverMutationType('Z',6,2);
            
            \sagepromt QuiverMutationType('Z',6,-2);
          }
      \end{itemize}
    \item mutation-infinite types
      \begin{itemize}
        \item infinite type $E$

          \sageex{
            \sagepromt QuiverMutationType('E',9,3)
            
            \sageret{['E', 8, 1]}
            
            \sagepromt QuiverMutationType('E',10,3);
            
            \sagepromt QuiverMutationType('E',12,3);
            
            \sagepromt QuiverMutationType('AE',(1,1),3);
            
            \sagepromt QuiverMutationType('AE',(1,4),3);
            
            \sagepromt QuiverMutationType('BE',5,3);
            
            \sagepromt QuiverMutationType('CE',5,3);
            
            \sagepromt QuiverMutationType('DE',6,3);
          }
        \item Grassmannian types -- the second parameter $(a,b)$ must satisfy $1\leq a<b$ and one obtains a grid graph of width $a-1$ and height $b-a-1$

          \sageex{
            \sagepromt QuiverMutationType('GR',(2,4),3)
            
            \sageret{['A', 1]}
            
            \sagepromt QuiverMutationType('GR',(2,6),3)
            
            \sageret{['A', 3]}
            
            \sagepromt QuiverMutationType('GR',(3,6),3)
            
            \sageret{['D', 4]}
            
            \sagepromt QuiverMutationType('GR',(3,7),3)
            
            \sageret{['E', 6]}
            
            \sagepromt QuiverMutationType('GR',(3,8),3)
            
            \sageret{['E', 8]}

            \sagepromt QuiverMutationType('GR',(3,9),3)
            
            \sageret{['E', 8, [1,1]]}
            
            \sagepromt QuiverMutationType('GR',(3,10),3);
          }
        \item triangular types -- the second parameter gives the size of the graph

          \sageex{
            \sagepromt QuiverMutationType('TR',2,3)
            
            \sageret{['A', 3]}
            
            \sagepromt QuiverMutationType('TR',3,3)
            
            \sageret{['D', 6]}
            
            \sagepromt QuiverMutationType('TR',4,3)
            
            \sageret{['E', 8, [1, 1]]}
            
            \sagepromt QuiverMutationType('TR',5,3);
          }
        \item type $T$ -- the second parameter gives the lengths of the three legs

          \sageex{
            \sagepromt QuiverMutationType('T',(1,1,1),3)
            
            \sageret{['A', 1]}
            
            \sagepromt QuiverMutationType('T',(1,1,4),3)
            
            \sageret{['A', 4]}
            
            \sagepromt QuiverMutationType('T',(1,4,4),3)
            
            \sageret{['A', 7]}
            
            \sagepromt QuiverMutationType('T',(2,2,2),3)
            
            \sageret{['D', 4]}
            
            \sagepromt QuiverMutationType('T',(2,2,4),3)
            
            \sageret{['D', 6]}
            
            \sagepromt QuiverMutationType('T',(2,3,3),3)
            
            \sageret{['E', 6]}
            
            \sagepromt QuiverMutationType('T',(2,3,4),3)
            
            \sageret{['E', 7]}
            
            \sagepromt QuiverMutationType('T',(2,3,5),3)
            
            \sageret{['E', 8]}
            
            \sagepromt QuiverMutationType('T',(2,3,6),3)
            
            \sageret{['E', 8, 1]}
            
            \sagepromt QuiverMutationType('T',(2,3,7),3)
            
            \sageret{['E', 10, 3]}
            
            \sagepromt QuiverMutationType('T',(3,3,3),3)
            
            \sageret{['E', 6, 1]}
            
            \sagepromt QuiverMutationType('T',(3,3,4),3);
          }
      \end{itemize}
    \item reducible types

          \sageex{
            \sagepromt QuiverMutationType(['A',3],['B',4])
            
            \sageret{[ ['A', 3], ['B', 4] ]}
          }
\end{itemize}

\noindent As described in Section~\ref{sec:Kacclassification}, one can use also Kac's classification types~\cite{Kac}.

\begin{remark} \rm Most of the above types have already been explained as Dynkin diagrams, appear in Kac's list, or in the classification work of Derksen-Owen \cite{DerkOwen}, and Felikson-Shapiro-Tumarkin \cite{FeSTu, FeSTuII}.  The exceptions to these are the triangular seeds, Grassmannian seeds, and the ``T'' seeds.  The first two of these describe a certain family of quivers that have certain shapes (as triangles and grids, respectively) and correspond to certain coordinate rings of geometric objects. (See Examples 4.4 and 4.6 of \cite{Keller2} or the source papers \cite{ClustIII} and \cite{Scott}.) The ``T'' family are those which correspond to ``Dynkin diagrams'' of the shape of a $T$ with a certain number of vertices on each arm and one central vertex.  

\sageex{

\sagepromt ClusterSeed(['TR',5,3]).show()

\sagepromt ClusterSeed(['GR',[5,11],3]).show()

\sagepromt ClusterSeed(['T',[4,4,5],3]).show()

\sageret{\includegraphics[width=2.2in]{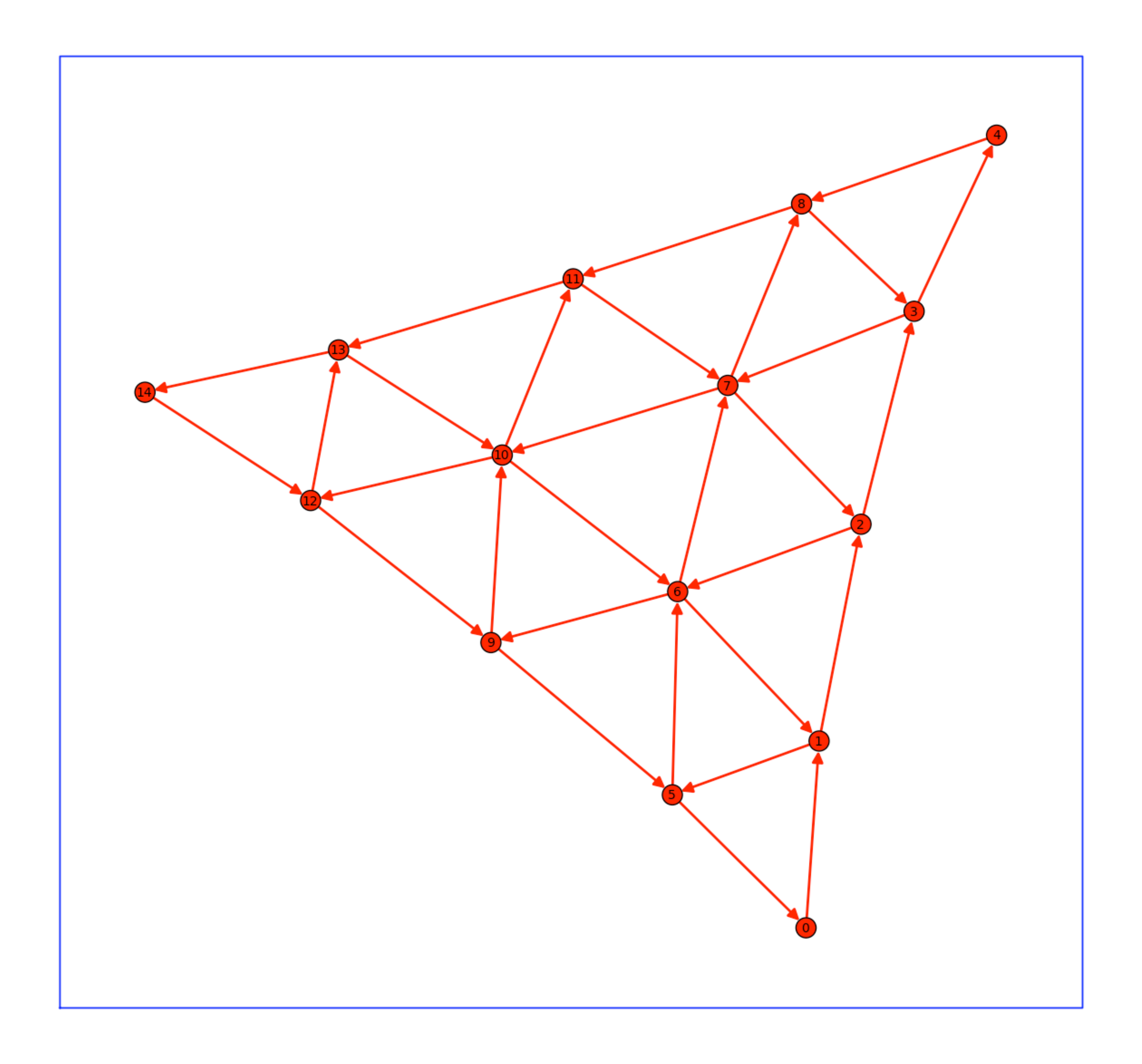}
\includegraphics[width=2.2in]{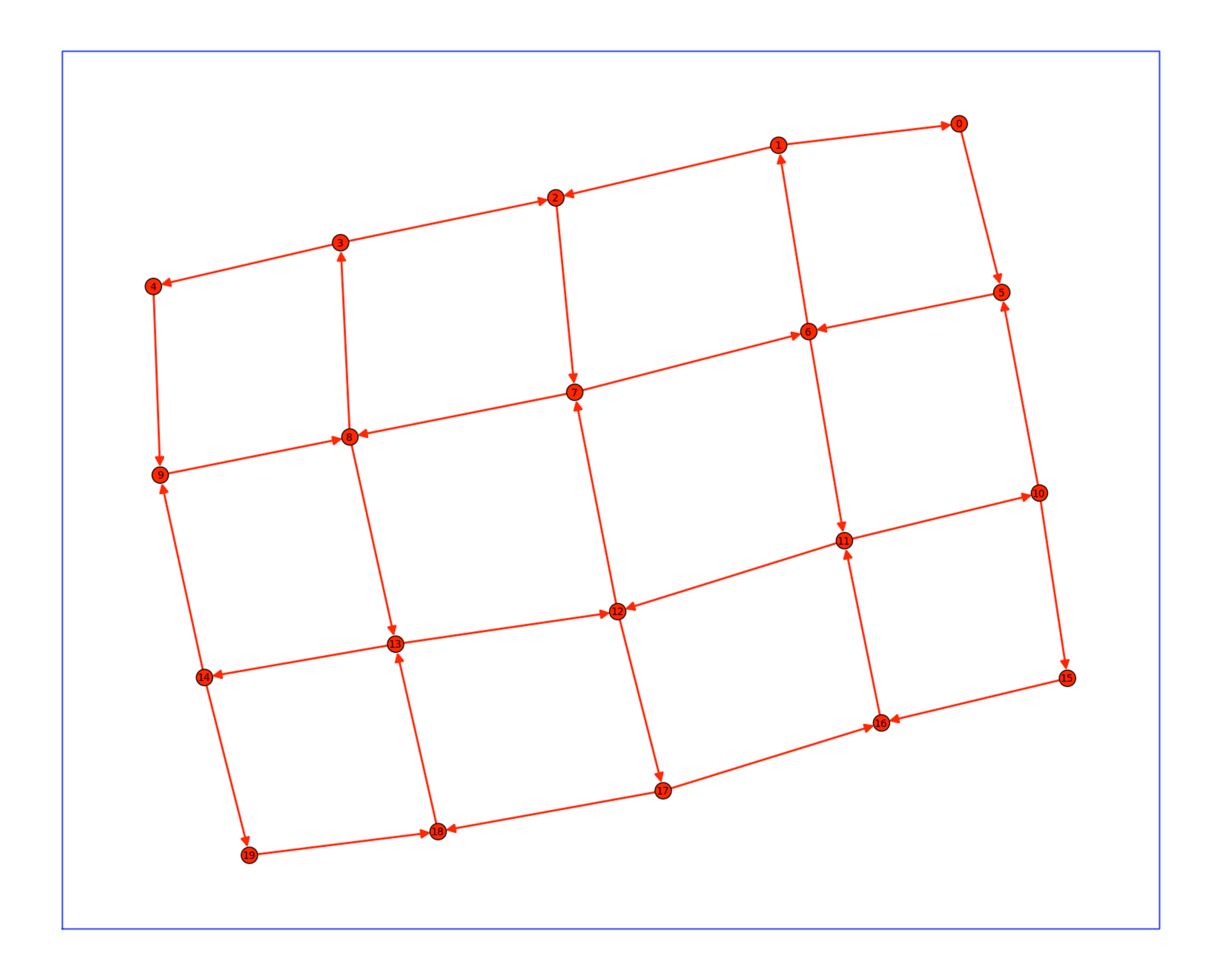}
\includegraphics[width=2.5in]{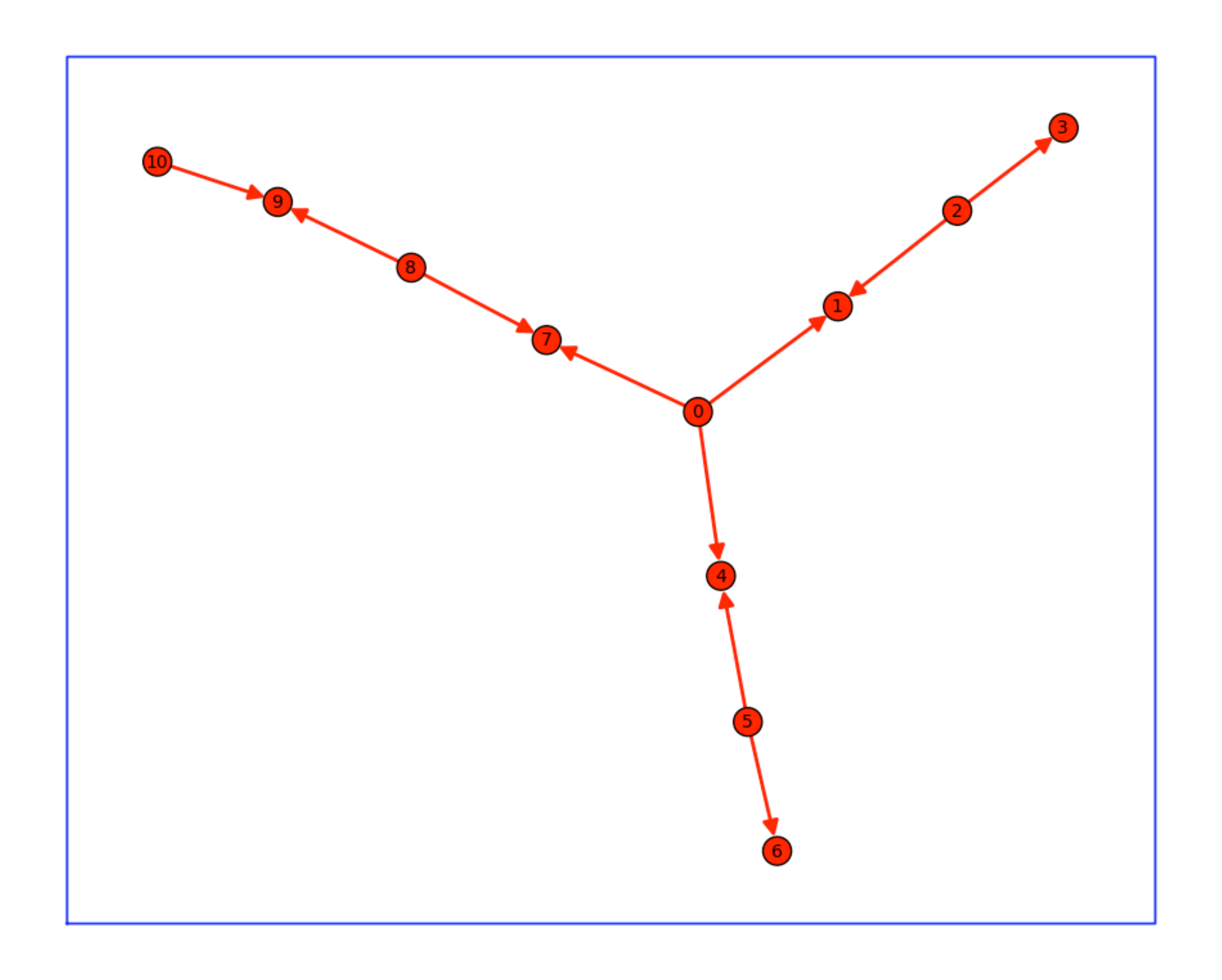}
}

}

We also illustrate a self-dual and two dual non-simply laced exceptional mutation-finite cases here too.

\sageex{
\sagepromt ClusterSeed(['X',6,2]).show()

\sagepromt S = ClusterSeed(['W',4,2]); S.show()

\sagepromt S = ClusterSeed(['W',4,-2]); S.show()
            
}            
            
\sageret{
\includegraphics[width=1.96in]{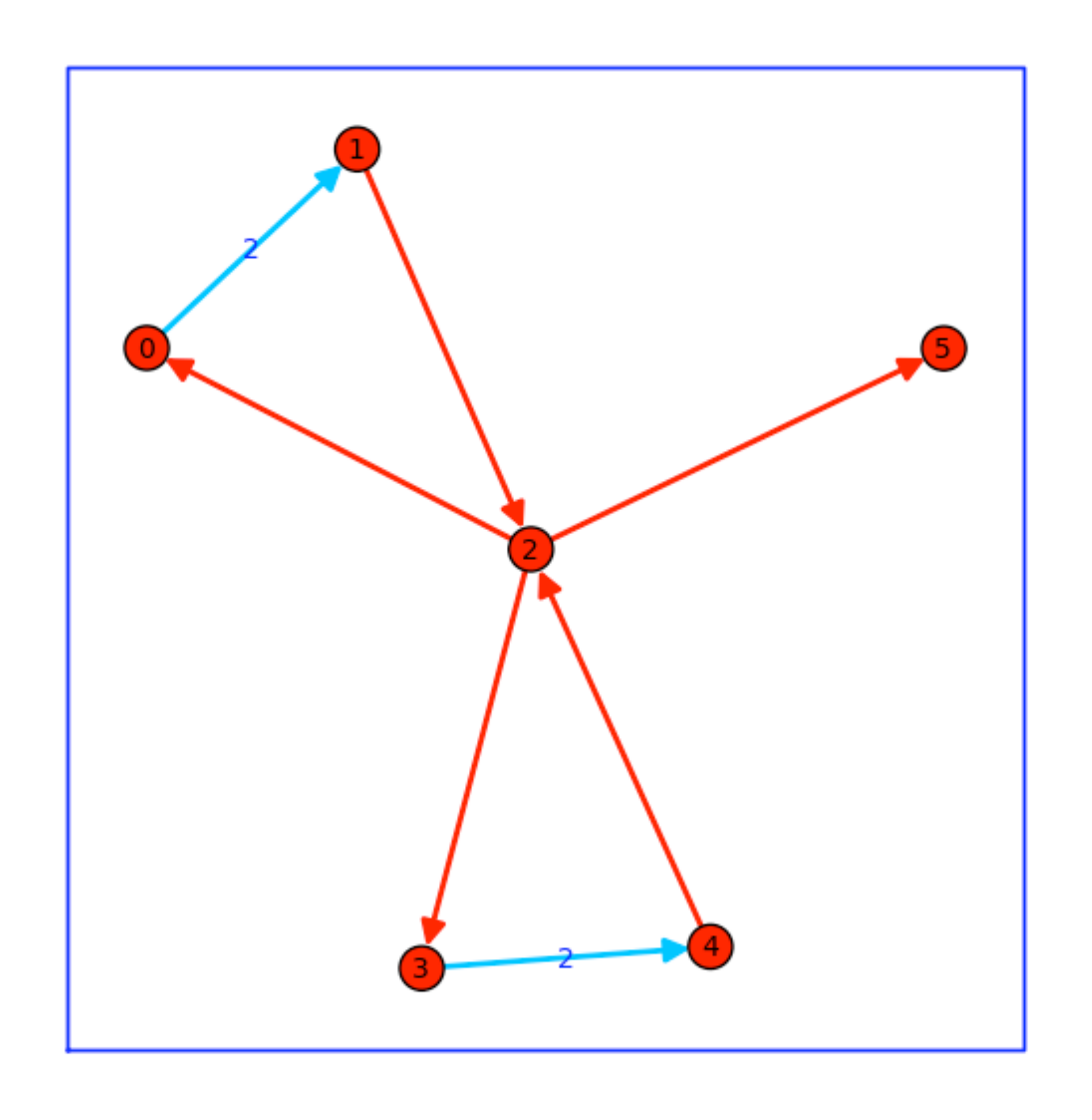}\hspace{3em}
\includegraphics[width=1in]{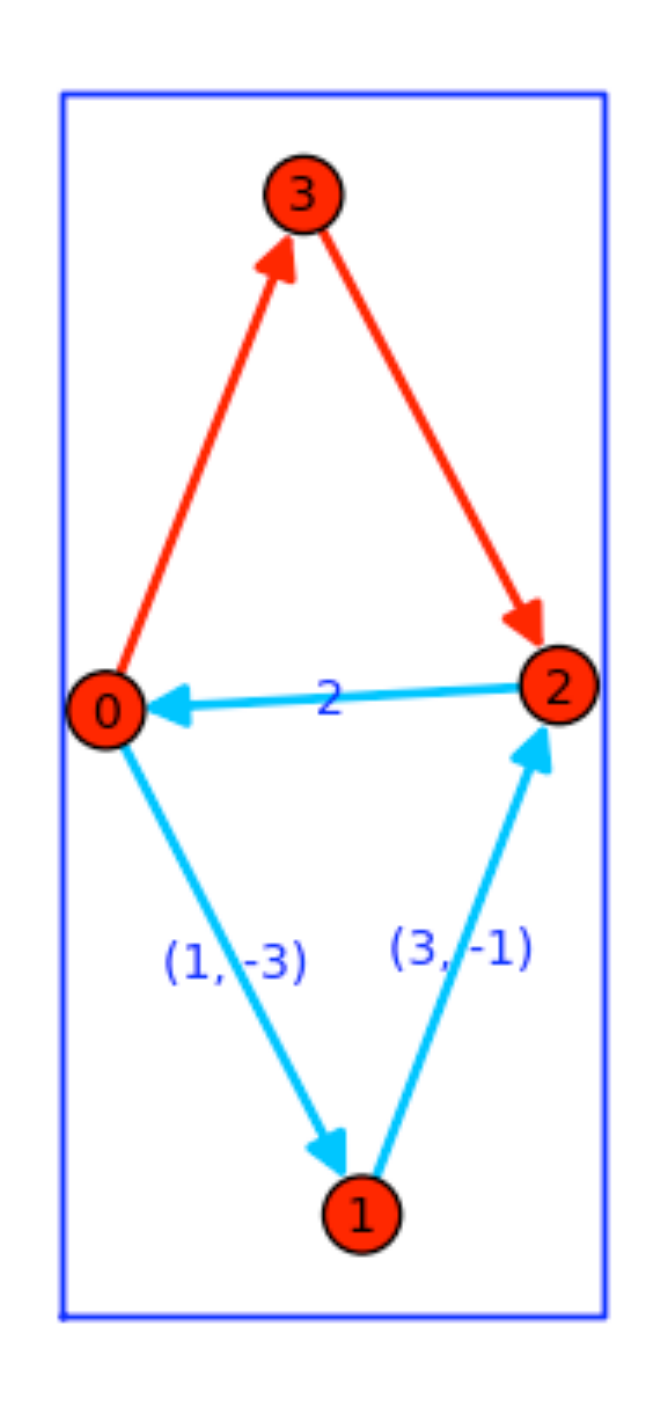} \hspace{3em}
\includegraphics[width=1in]{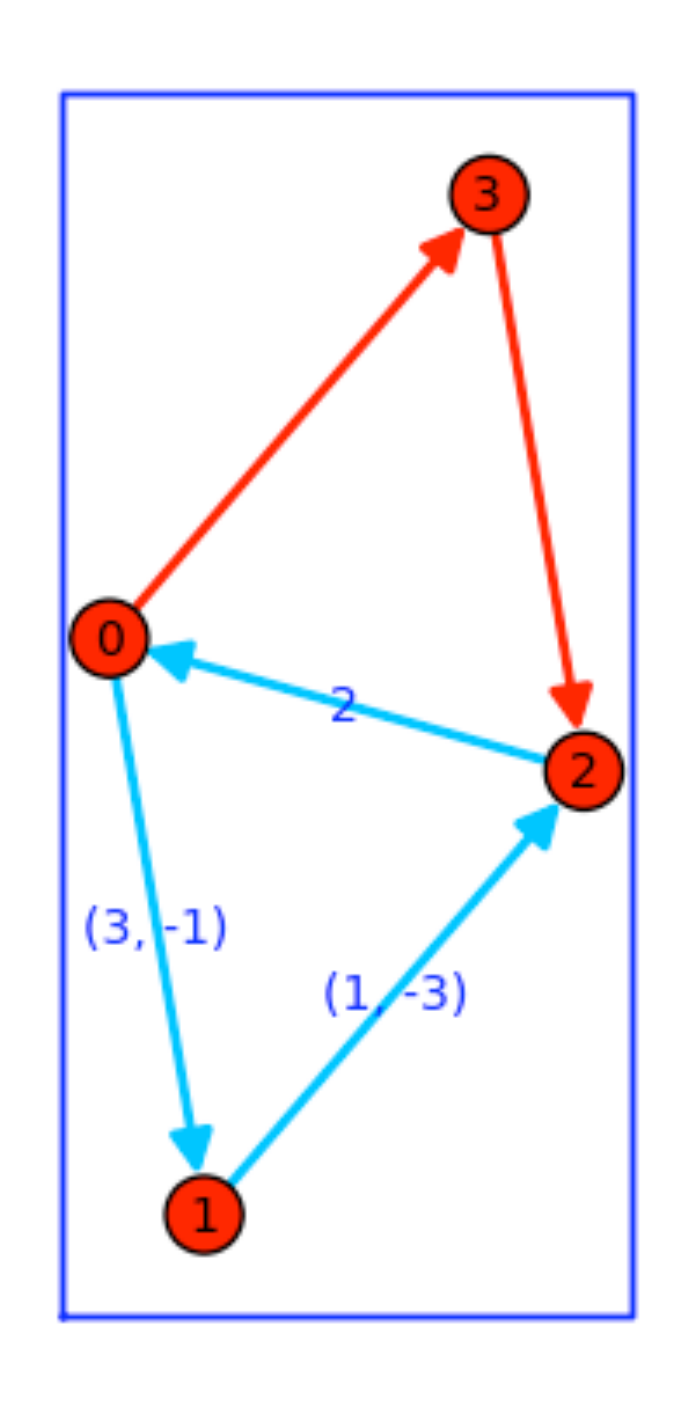}}

\end{remark}

The attributes of {\tt QuiverMutationType} are given by
\begin{itemize}
  \item {\tt \_letter}\\
  The string containing the letter(s) of the classification type.\\

  \item {\tt \_rank}\\
  The number of vertices in the standard quiver.\\

  \item {\tt \_bi\_rank}\\ 
  Is \None\, except for affine type $A$, where it denotes $[a,b]$ with $a+b$ being the rank and $a \leq b$ are the number of edges in the acyclic orientation of the standard quiver.\\
  
  \item {\tt \_twist}\\  
  Depends on the type of the classification type, and can be one of the following:
  \begin{itemize}
    \item \None\ for finite types,
    \item {\tt 1} for affine types,
    \item {\tt [1,1]} for elliptic types,
    \item {\tt 2} for finite mutation types which are not finite or elliptic,
    \item {\tt 3} for infinite mutation types.\\
  \end{itemize}
  
  \item {\tt \_graph}\\  
  Graph representing the underlying graph of the standard quiver.\\
  
  \item {\tt \_digraph}\\  
  Digraph representing the underlying graph of the standard quiver.\\
  
  \item {\tt \_description}\\  
  The string representation of the mutation class.\\
  
  \item {\tt \_info}\\  
  Dictionary containing the keys
  \begin{itemize}
    \item irreducible,
    \item finite,
    \item affine,
    \item elliptic,
    \item simply\_laced,
    \item mutation\_finite, and
    \item irreducible\_components.
  \end{itemize}
  The values are \True\ or \False, except for {\tt irreducible\_components} which is a list containing the irreducible components.
\end{itemize}

The methods of {\tt QuiverMutationType} are given by
\begin{itemize}
  \item {\tt \_\_eq\_\_(self,other)}\\  
  Returns \True, iff self and other represent the same quiver mutation type. As quiver mutation types are unique (i.e., there exists at most one object representing a given quiver mutation type), this method simply returns {\tt self is other}.\\

  \item {\tt \_repr\_(self)}\\ 
  Returns the string representation of self.\\

  \item {\tt plot(self, circular=False, directed=True)}\\  
  Returns a random or circular, directed or undirected plot of self.\\
  
  \item {\tt show(self, circular=False, directed=True)}\\
  Shows the plot of self.\\
  
  \item {\tt rank(self)}\\  
  Returns the rank (i.e., the number of vertices) of self.\\
  
  \item {\tt coxeter\_diagram(self)}\\  
  Returns the Coxeter diagram of self
  \sageex{
    
    \sagepromt QuiverMutationType(['A',5]).coxeter\_diagram()
    
    \sageret{Coxeter diagram of rank 5}

    \sagepromt QuiverMutationType(['A',3],['B',3]).coxeter\_diagram()
    
    \sageret{Coxeter diagram of rank 8}
  }
  
  \item {\tt b\_matrix(self)}\\  
  Returns the exchange matrix of self
  \sageex{

    \sagepromt QuiverMutationType(['A',5]).b\_matrix()

    \sageret{$\left(\begin{array}{rrrrr}0 & 1 & 0 & 0 & 0 \\-1 & 0 & -1 & 0 & 0 \\0 & 1 & 0 & 1 & 0 \\0 & 0 & -1 & 0 & -1 \\0 & 0 & 0 & 1 & 0\end{array}\right)$}

    \sagepromt QuiverMutationType(['A',3],['B',3]).b\_matrix()

    \sageret{$\left(\begin{array}{rrrrrr}0 & 1 & 0 & 0 & 0 & 0 \\-1 & 0 & -1 & 0 & 0 & 0 \\0 & 1 & 0 & 0 & 0 & 0 \\0 & 0 & 0 & 0 & 1 & 0 \\0 & 0 & 0 & -1 & 0 & -1 \\0 & 0 & 0 & 0 & 2 & 0\end{array}\right)$}
  }
  
  \item {\tt standard\_quiver(self)}\\  
  Returns the standard quiver of self.\\
  
  \item {\tt cartan\_matrix(self)}\\
  Returns the Cartan matrix of self which is obtained from its exchange matrix by replacing the positive entries by negative, and replace the $0$'s on the main diagonal by $2$'s.
  \sageex{

    \sagepromt QuiverMutationType('A',5).cartan\_matrix()

    \sageret{$\left(\begin{array}{rrrrr}2 & -1 & 0 & 0 & 0 \\-1 & 2 & -1 & 0 & 0 \\0 & -1 & 2 & -1 & 0 \\0 & 0 & -1 & 2 & -1 \\0 & 0 & 0 & -1 & 2\end{array}\right)$}

    \sagepromt QuiverMutationType(['A',3],['B',3]).cartan\_matrix()

    \sageret{$\left(\begin{array}{rrrrrr}2 & -1 & 0 & 0 & 0 & 0 \\-1 & 2 & -1 & 0 & 0 & 0 \\0 & -1 & 2 & 0 & 0 & 0 \\0 & 0 & 0 & 2 & -1 & 0 \\0 & 0 & 0 & -1 & 2 & -1 \\0 & 0 & 0 & 0 & -2 & 2\end{array}\right)$}
  }
  
  \item {\tt class\_size(self)}\\  
  Returns the number of quivers which are mutation-equivalent to self, up to isomorphism ({\bf Warning:} several class sizes are only conjectured, see Section~\ref{sec:classsizes}).
  \sageex{
  
  \sagepromt QuiverMutationType(['A',22],['BD',16,1]).class\_size()

  \sageret{4257164518523691840}

  \sagepromt QuiverMutationType(['GR',[4,9],3]).class\_size()

  \sageret{$\infty$}

  }

  \item {\tt dual(self)}\\  
  Returns the dual quiver mutation type of self.
  \sageex{
    
    \sagepromt QuiverMutationType('A',4).dual()
    
    \sageret{['A', 4]}
    
    \sagepromt QuiverMutationType('B',4).dual()
    
    \sageret{['C', 4]}
  }
  
  \item {\tt is\_irreducible(self)}\\  
  Returns \True, iff self is irreducible.
  \sageex{
  
  \sagepromt QuiverMutationType('A',4).is\_irreducible()
  
  \sageret{True}
  
  \sagepromt QuiverMutationType(['A',3],['B',3]).is\_irreducible()
  
  \sageret{False}
  }
  
  \item {\tt is\_mutation\_finite(self)}\\  
  Returns \True, iff self is of finite mutation type.
  \sageex{
  
  \sagepromt QuiverMutationType(['GR',[4,8],3]).is\_mutation\_finite()

  \sageret{True}

  \sagepromt QuiverMutationType(['GR',[4,9],3]).is\_mutation\_finite()

  \sageret{False}
  }

  \item {\tt is\_simply\_laced(self)}\\  
  Returns \True, iff self is simply-laced.\\
  
  \item {\tt is\_finite(self)}\\  
  Returns \True, iff self is of finite type.\\
  
  \item {\tt is\_affine(self)}\\  
  Returns \True, iff self is of affine type.\\
  
  \item {\tt is\_elliptic(self)}\\  
  Returns \True, iff self is of elliptic type.\\
  
  \item {\tt irreducible\_components(self)}\\  
  Returns a tuple containing the irreducible components of self.
  \sageex{

    \sagepromt QuiverMutationType('A',5).irreducible\_components()

    \sageret{(['A', 4],)}

    \sagepromt QuiverMutationType(['A',3],['B',3]).irreducible\_components()

    \sageret{(['A', 3], ['B', 3])}
  }
  
  \item {\tt properties(self)}\\
  Prints all properties of self. See Section~\ref{sec:typeclassification} for examples.
\end{itemize}

\subsection{Quiver} The next class we want to describe it the class {\tt Quiver}. It allows numerous ways to construct a quiver, several examples were described in Section~\ref{sec:quivers}.
  \begin{itemize}
    \item {\tt QuiverMutationType}
    \item {\tt list} or {\tt tuple} representing a quiver mutation type
    \item {\tt ClusterSeed}
    \item {\tt matrix}: a skew-symmetrizable matrix which represents the exchange matrix
    \item {\tt Quiver}
    \item {\tt DiGraph}: the digraph must represent a quiver
    \item {\tt list} of {\tt tuples} representing the edge list of a digraph for a quiver
  \end{itemize}
The attributes of {\tt Quiver} are given by
\begin{itemize}
  \item {\tt \_M}\\
  The exchange matrix of self.\\

  \item {\tt \_n}\\
  The number of cluster variables (which is the number of columns in the exchange matrix).\\
  
  \item {\tt \_m}\\
  The number of frozen variables (which is the number of rows minus the number of columns in the exchange matrix).\\
  
  \item {\tt \_description}\\
  The string representation of self.\\
  
  \item {\tt \_mutation\_type}\\
  The mutation type of self, if known, {\tt None} otherwise.
\end{itemize}
The methods of {\tt Quiver} are given by
\begin{itemize}
  \item {\tt \_\_init\_\_(self, data, frozen=0)}\\
  Frozen sets the later vertices to be frozen
  \sageex{

  \sagepromt Q1 = Quiver([(0,1),(1,2),(2,3)]); Q1

  \sageret{Quiver on 4 vertices}
  
  \sagepromt Q2 = Quiver([(0,1),(1,2),(2,3)],frozen=1); Q2
  
  \sageret{Quiver on 4 vertices with 1 frozen vertex}

  \sagepromt Q1.show()

  \sagepromt Q2.show()
  
  \sageret{\includegraphics[height=2in]{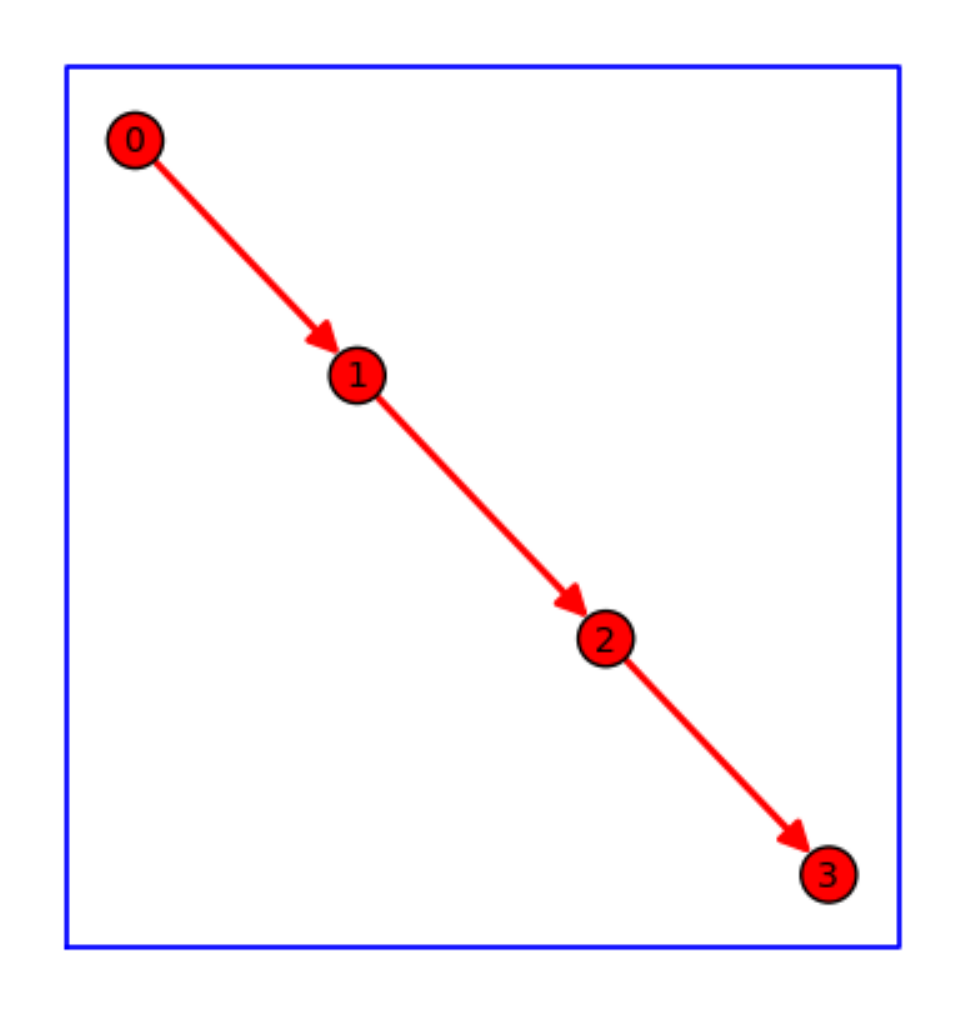} \hspace{2em}
  \includegraphics[height=2in]{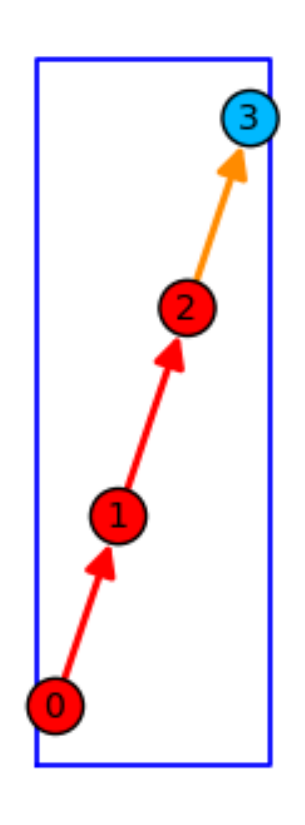}}
  }
  
  \item {\tt \_\_eq\_\_(self, other)}\\
  Returns \True, iff the b-matrices of self and other coincide
  \sageex{
  
  \sagepromt Q = Quiver(['A',5])
  
  \sagepromt T = Q.mutate( 2, inplace=False )
  
  \sagepromt Q.\_\_eq\_\_( T )
  
  \sageret{False}
  
  \sagepromt T.mutate( 2 )
  
  \sagepromt Q.\_\_eq\_\_( T )
  
  \sageret{True}
  }
  
  \item {\tt \_repr\_(self)}\\
  Returns the string representation of self
  \sageex{
  
  \sagepromt Q = Quiver(['A',5])
  
  \sagepromt Q.\_repr\_()
  
  \sageret{"Quiver on 5 vertices of type ['A', 5]"}
	}
	
  \item {\tt plot(self, circular=False, directed=True, mark=None)}\\
  Returns a random/circular and directed/undirected plot of self with a given vertex marked.\\
  
  \item {\tt show(self, fig\_size=1, circular=False, directed=True, mark=None)}\\
  Shows the plot of self.\\
  
  \item {\tt interact(self, fig\_size=1, circular=True)}\\
  Starts an interactive mode, as shown in Figure~\ref{fig:interactivemode} at the end of Section \ref{sec:quivers}.\\

  \item {\tt save\_image(self, filename, circular=False)}\\
  Saves the plot of self to filename.\\
  
  \item {\tt b\_matrix(self)}\\
  Returns the exchange matrix of self
  \sageex{
  
  \sagepromt Quiver(['A',4]).b\_matrix()
  
  $$\left(\begin{array}{rrrr}
    0 & 1 & 0 & 0 \\
    -1 & 0 & -1 & 0 \\
    0 & 1 & 0 & 1 \\
    0 & 0 & -1 & 0
    \end{array}\right)$$
    
  \sagepromt Quiver(['B',4]).b\_matrix()
    
  $$\left(\begin{array}{rrrr}
    0 & 1 & 0 & 0 \\
    -1 & 0 & -1 & 0 \\
    0 & 1 & 0 & 1 \\
    0 & 0 & -2 & 0
    \end{array}\right)$$
  
  \sagepromt Quiver(['D',4]).b\_matrix()
  
  $$\left(\begin{array}{rrrr}
    0 & 1 & 0 & 0 \\
    -1 & 0 & -1 & -1 \\
    0 & 1 & 0 & 0 \\
    0 & 1 & 0 & 0
    \end{array}\right)$$
  
  \sagepromt Quiver(QuiverMutationType([['A',2],['B',2]])).b\_matrix()
  
  $$\left(\begin{array}{rrrr}
    0 & 1 & 0 & 0 \\
    -1 & 0 & 0 & 0 \\
    0 & 0 & 0 & 1 \\
    0 & 0 & -2 & 0
    \end{array}\right)$$
  }

  \item {\tt digraph(self)}\\
  Returns the underlying digraph of self
  \sageex{
  
  \sagepromt Quiver(['A',4]).digraph()

  \sageret{Digraph on 4 vertices}

  }

  \item {\tt n(self)}\\
  Returns the number of free vertices of self
  \sageex{
  
  \sagepromt Q = Quiver([(0,1),(1,2),(2,3)],frozen=1)
  
  \sagepromt Q.n()
  
  \sageret{3}
  }
  
  \item {\tt m(self)}\\
  Returns the number of frozen vertices of self
  \sageex{
  
  \sagepromt Q = Quiver([(0,1),(1,2),(2,3)],frozen=1)
  
  \sagepromt Q.m()
  
  \sageret{1}
  }

  \item {\tt canonical\_label(self, certify=False)}\\
  Returns an isomorphic quiver with canonical vertex labeling. This is based on the canonical labeling of digraphs using the corresponding method for digraphs by R.L.~Miller based on \cite{McKay1981}. If certify is \True, a dictionary of the relabeling is also returned
  \sageex{
  
  \sagepromt Quiver(['A',4]).canonical\_label(certify=True)

	\sageret{(Quiver on 4 vertices of type ['A', 4], $\{0:0, 1:3, 2:1, 3:2\}$)}
  }

  \item {\tt is\_acyclic(self)}\\
  Returns \True, iff self is acyclic.\\
  
  \item {\tt is\_bipartite(self, return\_bipartition=False)}\\
  Returns \True, iff self is bipartite, if return\_bipartition is \True, the bipartition is returned
  \sageex{

  \sagepromt Quiver(['A',4]).is\_bipartite(return\_bipartition=True)

  \sageret{(set([0, 2]), set([1, 3]))}
  }

  \item {\tt principal\_restriction(self)}\\
  Returns the principal restricting of self. This is obtained from self by deleting all frozen variables.\\
  
  \item {\tt principal\_extension(self)}\\
  Returns the principal extension of self. This can be used only for seeds without frozen variables. Returns a new seed with exchange matrix of size $2n \times n$ given by the exchange matrix of self of size $n \times n$ with an additional identity matrix added below.\\
  
  \item {\tt mutate(self, data, inplace=True)}\\
  Mutates at a vertex or at a list of vertices, if inplace is \True, self is modified, otherwise a new quiver is returned.\\
  
  \item {\tt mutation\_sequence(self, sequence, show\_sequence=False,}\\ \hspace*{1em}{\tt fig\_size=1.2)}\\
  Returns a list of quivers obtained from a sequence of mutations. If the parameter show\_sequence is \True, the sequence is shown with a given fig\_size.\\
  
  \item {\tt reorient(self,data)}\\
  Reorients self with respect to the given total order, or with respect to an iterator of edges in self to be reverted.  \\
 
{\bf Warning:} This often will change the mutation class of the quiver except if the quiver is a tree (see Theorem \ref{th:tree}).
\\
 
  \item {\tt mutation\_class\_iter(self, depth=infinity,}\\ \hspace*{1em}{\tt show\_depth=False, return\_paths=False,} \\ \hspace*{2em}{\tt data\_type='quiver', up\_to\_equivalence=True,}\\ \hspace*{3em}{\tt only\_sink\_source=False)}\\
  Returns an iterator which goes through the mutation class of self depending on several parameters
    \begin{itemize}
      \item {\tt depth}: integer, only quivers with distance at most depth from self are returned
      \item {\tt show\_depth}: if \True, the actual depth of the mutation is shown
      \item {\tt return\_paths}: if \True, a shortest path of mutation sequences from self to the given quiver is returned as well
      \item {\tt data\_type}: can be one of the following:
      \sageex{\sageret{quiver, matrix, digraph, dig6, path}}
      \item {\tt up\_to\_equivalence}: if True, only quivers up to equivalence are considered
      \item {\tt sink\_source}: if True, only mutations at sinks and sources are applied\\
    \end{itemize}

  \item {\tt mutation\_class(self, depth=infinity,}\\ \hspace*{1em}{\tt show\_depth=False, return\_paths=False,}\\ \hspace*{2em}{\tt data\_type='quiver', up\_to\_equivalence=True,}\\ \hspace*{3em}{\tt only\_sink\_source=False)}\\
  Returns a list of all quivers in the corresponding iterator.\\
  
  \item {\tt group\_of\_mutations(self)}\\
  Returns the group of mutations of self. {\bf Warning:} The permutation group is very big! This group differs for quivers and for cluster seeds, as different cluster seeds may have the same exchange matrix and thus the same quiver. This group is defined to be the group of permutation given as follows. The ground set is the mutation class of self without taking equivalence of quivers into account, and the group is generated by the $n$ involutions on this set given by mutation at the $n$ different vertices. Observe that the analogous operation on the mutation class up to equivalence does not give a group (this can be easily checked in type $A_3$). Basically nothing is known about this group.  
  \sageex{
  
  \sagepromt Q = Quiver(['A',2])
  
  \sagepromt Q.group\_of\_mutations()
  
  \sageret{Permutation Group with generators [(1,2)]}
  
  \sagepromt Q = Quiver(['A',3])
  
  \sagepromt Q.group\_of\_mutations()
  
  \sageret{Permutation Group with generators [(1,2)(3,4)(5,9)(6,7)(8,12)(10,11)(13,14), (1,3)(2,5)(4,6)(7,14)(8,11)(9,13)(10,12), (1,4)(2,3)(5,10)(6,8)(7,13)(9,14)(11,12)]}

  \sagepromt Q = Quiver(['B',2])
  
  \sagepromt Q.group\_of\_mutations()
  
  \sageret{Permutation Group with generators [(1,2)]}
  
  \sagepromt Q = Quiver(['B',3])
  
  \sagepromt Q.group\_of\_mutations()
  
  \sageret{Permutation Group with generators [(1,2)(3,4)(5,6)(7,10)(8,9), (1,3)(2,6)(4,5)(7,9)(8,10), (1,4)(2,3)(5,7)(6,8)(9,10)]}
  
  \sagepromt Q = Quiver(['A',1])
  
  \sagepromt Q.group\_of\_mutations().cardinality()
  
  \sageret{1}
  
  \sagepromt Q = Quiver(['A',2])
  
  \sagepromt Q.group\_of\_mutations().cardinality()
  
  \sageret{2}

  \sagepromt Q = Quiver(['A',3])
  
  \sagepromt Q.group\_of\_mutations().cardinality()
  
  \sageret{322560}
  }

  \item {\tt is\_finite(self)}\\
  Returns \True, iff self is of finite type. This is done by checking if it is mutation-equivalent to a quiver of finite type.\\
  
  \item {\tt is\_mutation\_finite(self, nr\_of\_checks=None, return\_path=False)}\\
  Returns \True, iff self if of finite mutation type. {\bf Warning:} The algorithm is non-deterministic and uses random mutations in various directions. Might theoretically result in a wrong \True\ return. The number of checks can be set, the default is $1000$ times the number of vertices of self. If {\tt return\_path} is \True, then a path to a non-mutation-finite quiver is returned, if found.\\
  
  \item {\tt mutation\_type(self)}\\
  Returns the mutation type of self if it can be determined.
    \begin{itemize}
      \item First, it is checked if self is mutation-equivalent to a quiver of a classical type using the descriptions of the classification types,
      \item then, it is checked if self is contained in an exceptional mutation class which are hard-coded,
      \item if it was not possible to determine the mutation type, it is checked if self is mutation-finite or infinite
    \end{itemize}
    {\bf Warning:} The algorithm to determine quivers of mutation type $\tilde D_n$ (which is {\tt ['D',n,1]}) is not yet implemented!
\end{itemize}

\subsection{ClusterSeed} The constructor of the class {\tt ClusterSeed} allows the same input as the class {\tt Quiver} to construct a cluster seed. Moreover, many attributes and methods for cluster seeds and for quivers coincide. Often, the cluster seed simply calls the quiver method.
  \begin{itemize}
    \item {\tt QuiverMutationType}
    \item {\tt list} or {\tt tuple} representing a quiver mutation type
    \item {\tt ClusterSeed}
    \item {\tt matrix}: a skew-symmetrizable matrix which represents the exchange matrix
    \item {\tt Quiver}
    \item {\tt DiGraph}: the digraph must represent a quiver
    \item {\tt list} of {\tt tuples} representing the edge list of a digraph for a quiver
  \end{itemize}
The attributes of {\tt ClusterSeed} are given by
\begin{itemize}
  \item {\tt \_M}\\
  The exchange matrix of self.\\
  
  \item {\tt \_cluster}\\
  The cluster as a list of cluster variables.\\
  
  \item {\tt \_n}\\
  The number of cluster variables (which is the number of columns in the exchange matrix).\\
  
  \item {\tt \_m}\\
  The number of frozen variables (which is the number of rows $-$ the number of columns in the exchange matrix).\\
  
  \item {\tt \_R}\\
  The base ring in which the cluster variables live.\\
  
  \item {\tt \_quiver}\\
  The quiver attached to self.\\
  
  \item {\tt \_description}\\
  The string representation of self.\\
  
  \item {\tt \_mutation\_type}\\
  The mutation type of self, if known, \None\ otherwise
\end{itemize}
The methods of {\tt ClusterSeed} are given by
\begin{itemize}
  \item {\tt \_\_init\_\_(self, data, frozen=0)}\\
  Frozen sets the later vertices to be frozen
  \sageex{

  \sagepromt S1 = ClusterSeed([(0,1),(1,2),(2,3)]); S1

  \sageret{A seed for a cluster algebra of rank 4}
  
  \sagepromt Q2 = Quiver([(0,1),(1,2),(2,3)],frozen=1); Q2
  
  \sageret{A seed for a cluster algebra of rank 3 with 1 frozen variable}

  \sagepromt Q1.b\_matrix(); Q2.b\_matrix()
  
  $$
  \left(\begin{array}{rrrr}0 & 1 & 0 & 0 \\-1 & 0 & 1 & 0 \\0 & -1 & 0 & 1 \\0 & 0 & -1 & 0\end{array}\right) \qquad
  \left(\begin{array}{rrr}0 & 1 & 0 \\-1 & 0 & 1 \\0 & -1 & 0 \\0 & 0 & -1\end{array}\right)
  $$
  }

  \item {\tt \_\_eq\_\_(self, other)}\\
  Returns \True, iff self and other have the same exchange matrix and the same cluster.\\
  
  \item {\tt \_repr\_(self)}\\
  Returns the string representation of self
  \sageex{
  
  \sagepromt S = ClusterSeed(['A',3]); S.\_repr\_()
  
  \sageret{"A seed for a cluster algebra of rank 3 of type ['A', 3]"}
  }
  
  \item {\tt plot(self, circular=False, mark=None)}\\
  Returns a random/circular plot of self with a given marked vertex. Calls the method for quivers.\\
  
  \item {\tt show(self, fig\_size=1, circular=False, mark=None)}\\
  Shows the plot of self.\\
  
  \item {\tt interact(self, fig\_size=1, circular=True)}\\
  Starts an interactive mode, as shown in Figure~\ref{fig:interactivemode} at the end of Section \ref{sec:quivers}.\\

  \item {\tt save\_image(self, filename, circular=False)}\\
  Saves a plot of self to filename.\\
  
  \item {\tt b\_matrix(self)}\\
  Returns the exchange matrix of self.\\
  
  \item {\tt cluster(self)}\\
  Returns the cluster of self
  \sageex{
  
  \sagepromt S = ClusterSeed(['A',3]); S.cluster()

  $$[x_0, x_1, x_2]$$
  
  \sagepromt S.mutate(0); S.cluster()

  $$\left[ \frac{x_1 + 1}{x_0}, x_1, x_2 \right]$$
  
  \sagepromt S.mutate(1); S.cluster()

  $$\left[\frac{x_1 + 1}{x_0}, \frac{x_0 x_2 + x_1 + 1}{x_0 x_1}, x_2\right]$$
  }
  
  \item {\tt ground\_field(self)}\\
  Returns the ground field in which the cluster variables of self live
  \sageex{
  
  \sagepromt S.ground\_field()
  
  \sageret{Fraction Field of Multivariate Polynomial Ring \\in x0, x1, x2 over Rational Field}
  }
  
  \item {\tt x(self,k)}\\
  Returns the $k$th initial cluster variable of self.\\
  
  \item {\tt y(self,k)}\\
  Returns the $k$th frozen variable of self.\\
  
  \item {\tt n(self)}\\
  Returns the number of cluster variables of self.\\
  
  \item {\tt m(self)}\\
  Returns the number of frozen variables of self.\\
  
  \item {\tt exchangeable\_variables(self)}\\
  Returns a {\tt list} of all cluster variables of self.\\
  
  \item {\tt frozen\_variables(self)}\\
  Returns a {\tt list} of all frozen variables of self.\\
  
  \item {\tt quiver(self)}\\
  Returns the {\tt Quiver} associated to self.\\
  
  \item {\tt is\_acyclic(self)}\\
  Returns \True, iff self is acyclic.\\
  
  \item {\tt is\_bipartite(self, return\_bipartition=False)}\\
  Returns \True, iff self is bipartite, if return\_bipartition is \True, the bipartition is returned
  \sageex{
  
  \sagepromt ClusterSeed(['A',4]).is\_bipartite(return\_bipartition=True)

  \sageret{(set([0, 2]), set([1, 3]))}
  }
  \item {\tt mutate(self, sequence, inplace=True)}\\
  Mutates at an index or at a list of indices, if inplace is \True, self is modified, otherwise a new cluster seed is returned.\\
  
  \item {\tt mutation\_sequence(self, sequence, show\_sequence=False,}\\ \hspace*{1em}{\tt fig\_size=1.2,return\_output='seed')}\\
  Returns a list depending on {\tt return\_output} obtained from a sequence of mutations. If show\_sequence is \True, the sequence is shown with a given fig\_size. The possible outputs are
  \begin{itemize}
   \item {\tt 'seed'}: a list of cluster seeds is returned
   \item {\tt 'matrix'}: a list of exchange matrices is returned
   \item {\tt 'var'}: a list of cluster variables is returned\\
  \end{itemize}
  
  \item {\tt principal\_restriction(self)}\\
  Returns the principal restriction of self. This is obtained from self by deleting all frozen variables.\\
  
  \item {\tt principal\_extension(self)}\\
  Returns the principal extension of self. This can be used only for seeds without frozen variables. Returns a new seed with exchange matrix of size $2n \times n$ given by the exchange matrix of self of size $n \times n$ with an additional identity matrix added below.\\
  
  \item {\tt reorient(self,data)}\\
  Reorients self by reorienting the corresponding quiver. Calls the method for quivers.\\
  
  \item {\tt set\_cluster(self, cluster)}\\
  Sets the set of clusters of self to {\tt cluster}.\\
  
  \item {\tt reset\_cluster(self)}\\
  Sets the set of clusters of self back to the initial cluster.\\
  
  \item {\tt mutation\_class\_iter(self, depth=infinity,}\\ \hspace*{1em}{\tt show\_depth=False, return\_paths=False,} \\ \hspace*{2em}{\tt up\_to\_equivalence=True, only\_sink\_source=False)}\\
  Returns an iterator which goes through the mutation class of self depending on several parameters
    \begin{itemize}
      \item {\tt depth}: integer, only quivers with distance at most depth from self are returned
      \item {\tt show\_depth}: if \True, the actual depth of the mutation is shown
      \item {\tt return\_paths}: if \True, a shortest path of mutation sequences from self to the given quiver is returned as well
      \item {\tt up\_to\_equivalence}: if \True, only quivers up to equivalence are considered
      \item {\tt only\_sink\_source}: if \True, only mutations at sinks and sources are applied\\
    \end{itemize}
  
  \item {\tt mutation\_class(self, depth=infinity,}\\ \hspace*{1em}{\tt show\_depth=False, return\_paths=False,}\\ \hspace*{2em}{\tt up\_to\_equivalence=True, only\_sink\_source=False)}\\
  Returns a list of all quivers in the corresponding iterator\\

  \item {\tt cluster\_class\_iter(self, depth=infinity, show\_depth=False,\\
  \hspace*{1em}up\_to\_equivalence=True)}\\
  Returns an iterator through all clusters mutation-equivalent to self up to a given depth. Moreover, it is possible to show the actual depth together with several parameters, or to output clusters as labeled seeds.\\

  \item {\tt cluster\_class(self, depth=infinity, show\_depth=False,\\
  \hspace*{1em}up\_to\_equivalence=True)}\\
  Returns a list of all clusters mutation-equivalent to self up to a given depth. Moreover, it is possible to show the actual depth together with several parameters, or to output cluster as labeled seeds.\\
  
  \item {\tt b\_matrix\_class\_iter(self, depth=infinity,\\
  \hspace*{1em}up\_to\_equivalence=True)}\\
  Returns an iterator through all matrices mutation-equivalent to self up to a given depth, and up to permutation of rows and columns unless specified otherwise. \\
  
  \item {\tt b\_matrix\_class(self, depth=infinity, up\_to\_equivalence=True)}\\
  Returns a list of all matrices mutation-equivalent to self up to a given depth, and up to permutation of rows and columns unless specified otherwise. \\

  \item {\tt variable\_class\_iter(self, depth=infinity,\\
  \hspace*{1em}ignore\_bipartite\_belt=False)}\\
  Returns an iterator through all variables obtained from self by mutations up to a given depth. {\bf Warning:} If at some point a bipartite seed is reached, another algorithm is used unless the parameter {\tt ignore\_bipartite\_belt} is set to be {\tt True}. See the description in Section \ref{sec:varclass}.\\

  \item {\tt variable\_class(self, depth=infinity,\\
  \hspace*{1em}ignore\_bipartite\_belt=False)}\\
  Returns a list of all variables obtained from self by mutations up to a given depth. {\bf Warning:} If at some point a bipartite seed is reached, another algorithm is used unless the parameter {\tt ignore\_bipartite\_belt} is set to be {\tt True}. See the description in Section \ref{sec:varclass}.\\

  \item {\tt group\_of\_mutations(self)}\\
  Returns the group of mutations of self. {\bf Warning:} The permutation group is very big! This group differs for quivers and for cluster seeds, as different cluster seeds may have the same exchange matrix and thus the same quiver. This group is defined to be the group of permutation given as follows. The ground set is the mutation class of self without taking equivalence of seeds into account, and the group is generated by the $n$ involutions on this set given by mutation at the $n$ different vertices. Observe that the analogous operation on the mutation class up to equivalence does not give a group (this can be easily checked in type $A_3$). Basically nothing is known about this group.  
  \sageex{
  
  \sagepromt S = ClusterSeed(['A',2])
  
  \sagepromt S.group\_of\_mutations()
  
  \sageret{Permutation Group with generators [(1,2)(3,4)(5,6)(7,9)(8,10), (1,3)(2,5)(4,7)(6,8)(9,10)]}

  \sagepromt S = ClusterSeed(['B',2])
  
  \sagepromt S.group\_of\_mutations()
  
  \sageret{Permutation Group with generators \\[0pt] [(1,2)(3,4)(5,6), (1,3)(2,5)(4,6)]}

  \sagepromt Q = ClusterSeed(['A',1])
  
  \sagepromt Q.group\_of\_mutations().cardinality()
  
  \sageret{2}

  \sagepromt Q = ClusterSeed(['A',2])
  
  \sagepromt Q.group\_of\_mutations().cardinality()
  
  \sageret{10}

  \sagepromt Q = ClusterSeed(['A',3])
  
  \sagepromt Q.group\_of\_mutations().cardinality()
  
  \sageret{705438720}
  }

  \item {\tt is\_finite(self)}\\
  Returns \True, iff self is of finite type. Calls the method for the quiver of self.\\
  
  \item {\tt is\_mutation\_finite(self, nr\_of\_checks=None, return\_path=False)}\\
  Returns \True, iff self is of finite mutation type. Calls the method for the quiver of self.\\
  
  \item {\tt mutation\_type(self)}\\
  Returns the mutation type of self, if possible. Calls the method for the quiver of self.
\end{itemize}

\subsection{ClusterVariable}

By definition, a cluster variable is an element in the field of rational function in $n$ variables\footnote{Moreover, by Theorem \ref{th:Laurent}, they are actually multivariate Laurent polynomials in $n$ variables, although for the moment we do not use this functionality.}. The class {\tt ClusterVariable} provides two extra features for cluster variables:
\begin{enumerate}
 \item The connection to almost positive roots in finite types (positive roots are not yet provided in \sage\ for affine types).
 \item An ordering for cluster variables which is inspired by its connection to almost positive roots:
  \begin{itemize}
   \item They are ordered first by total degree of the denominator (in particular, the variables in the initial seed come first in natural order),
   \item If the degree is equal and positive, they are ordered lexicographically with $x_0 > x_1 > \ldots > x_{n-1}$.
  \end{itemize}
\end{enumerate}
  \sageex{

    \sagepromt for f in ClusterSeed(['A',2]).variable\_class():

    \sagedots \hspace{10pt} print f, f.almost\_positive\_root()

    \begin{align*}
      x_0                       &\qquad -\alpha_1\\
      x_1                       &\qquad -\alpha_2\\
      (x_1 + 1)/x_0             &\qquad \alpha_1\\
      (x_0 + 1)/x_1             &\qquad \alpha_2\\
      (x_0 + x_1 + 1)/(x_0 x_1) &\qquad \alpha_1 + \alpha_2
    \end{align*}
  }  
  
  \noindent Two further examples of the ordering can be found in Section~\ref{sec:Kacclassification}.  It is planned to include more functionalities for the cluster variable class in the future.

\end{document}